\documentclass[11pt]{article}

\usepackage[sectionbib]{natbib}
\usepackage{algorithm}
\usepackage{algorithmic}
\usepackage{booktabs}
\usepackage{lineno}
\usepackage{graphicx}
\usepackage{color,xcolor}
\usepackage{subfigure} 
\usepackage{geometry}
\usepackage{multirow}

\usepackage[colorlinks=true,linkcolor=blue,citecolor=blue]{hyperref}%

\usepackage{amsmath}
\usepackage{amssymb}
\usepackage{amsthm}

\usepackage[utf8]{inputenc} 
\usepackage{hyperref}       
\usepackage{url}            
\usepackage{booktabs}       
\usepackage{amsfonts}       
\usepackage{nicefrac}       
\usepackage{microtype}      

\usepackage{natbib}

\usepackage{array,epsfig,rotating}
\usepackage{tabularx}
\usepackage{subfigure}
\usepackage{graphicx}
\usepackage{verbatim}
\usepackage{arydshln}
\usepackage{authblk}

\newtheorem{theorem}{Theorem}
\newtheorem{lemma}{Lemma}

\newtheorem{prop}{Proposition}

\newtheorem{remark}{Remark}

\newtheorem{assumption}{Assumption}

\allowdisplaybreaks

\def\hat{\widehat}
\def\tilde{\widetilde}

\newcommand{\aaa}{\boldsymbol \alpha}
\newcommand{\bbb}{\boldsymbol \beta}
\newcommand{\rrr}{\boldsymbol \gamma}

\newcommand{\bw}{\boldsymbol w}

\newcommand{\bA}{\mathbf A}

\newcommand{\bp}{{\boldsymbol p}}

\newcommand{\bX}{\mathbf X}
\newcommand{\bY}{\mathbf Y}
\newcommand{\bZ}{\mathbf Z}

\newcommand{\bz}{\mathbf z}

\newcommand{\bN}{\mathbf N}

\newcommand{\bI}{\mathbf I}

\newcommand{\bmu}{\boldsymbol \mu}
\newcommand{\bpi}{\boldsymbol \pi}

\newcommand{\ME}{\mathbb E}
\newcommand{\MP}{\mathbb P}
\newcommand{\MU}{\mathbb U}
\newcommand{\Var}{\operatorname{Var}}

\newcommand\independent{\protect\mathpalette{\protect\independenT}{\perp}}
\def\independenT#1#2{\mathrel{\rlap{$#1#2$}\mkern2mu{#1#2}}}

\def\spacingset#1{\renewcommand{\baselinestretch}%
{#1}\small\normalsize} \spacingset{1}
\spacingset{1}

\usepackage{float}

\usepackage{mathtools}

\usepackage{dsfont}

\begin{document}

\title{Causal inference with high-dimensional discrete covariates}

\author{Zhenghao Zeng$^1$, Sivaraman Balakrishnan$^2$, Yanjun Han$^3$, Edward H. Kennedy$^2$  \\ 
    $^1$Graduate School of Business, Stanford University \\ 
    $^2$Department of Statistics and Data Science, Carnegie Mellon University \\ 
    $^3$Courant Institute of Mathematical Sciences and Center for Data Science, New York University
    }

\date{}

\maketitle

\begin{abstract}
When estimating causal effects from observational studies, researchers often need to adjust for many covariates to deconfound the non-causal relationship between exposure and outcome, and often many covariates are discrete. The behavior of commonly used estimators in the presence of many discrete covariates is not well understood, since standard approaches often employ structural assumptions such as smoothness, which do not apply in discrete settings. In this work, we study estimation of causal effects in a model where the covariates required for confounding adjustment are discrete but high-dimensional, meaning the number of categories $d$ can be comparable to or even larger than sample size $n$. Specifically, we show the mean squared error of commonly used regression, weighting and doubly robust estimators is bounded by $\frac{d^2}{n^2}+\frac{1}{n}$. We then prove that the minimax lower bound for estimating the average treatment effect is of order $\frac{d^2}{n^2 \log^2 n}+\frac{1}{n}$, which characterizes the fundamental difficulty of causal effect estimation in the high-dimensional discrete setting, and shows the estimators mentioned above are rate-optimal up to log factors. Finally we consider two other kinds of structure that can be exploited: effect homogeneity, and prior knowledge of the covariate distribution. We propose new estimators that enjoy faster convergence rates here, of order $\frac{d}{n^2} + \frac{1}{n}$, thus achieving consistency in a broader regime. The results are illustrated empirically via simulation studies and a real data example.
\end{abstract}


\def\spacingset#1{\renewcommand{\baselinestretch}%
{#1}\small\normalsize} \spacingset{1}
\spacingset{1}

\section{Introduction}
To draw causal conclusions from observational studies, researchers typically have to measure and adjust for many covariates (e.g., any confounders that could affect  treatment assignment or the outcome of interest). Numerous approaches exist for estimating causal effects from such data. Common methods include outcome modeling \citep{rubin1979using, hernan2010causal}, inverse propensity score weighting \citep{rosenbaum1983central, hahn1998role, hirano2003efficient}, and semiparametric methods (i.e., doubly robust estimation, or augmented inverse probability weighting, or targeted or double/debiased machine learning) \citep{van2003unified, chernozhukov2018double, kennedy2022semiparametric}. Notably, the latter semiparametric methods are consistent for the average treatment effect (ATE) even when either the outcome regression or propensity score is misspecified, yielding a robust tool for estimating ATEs under possible model misspecification. More generally, the bias of these methods is a product of errors in outcome regression and propensity score estimation, thus allowing parametric rates of convergence (and semiparametric efficiency) even when nuisance functions are estimated at slower nonparametric rates. Sample splitting or cross-fitting \citep{robins2008higher, zheng2010asymptotic, chernozhukov2018double} are often used to prevent overfitting the nuisance functions and avoid empirical process conditions \citep{kennedy2022semiparametric}. In the classic setting where the covariate dimension is held fixed, while the sample size grows to infinity, the aforementioned methods can enjoy  favorable properties (e.g., $\sqrt{n}$-consistency and asymptotic normality) and performance in applications. However, when the number of covariates needed for confounding adjustment is large, some new and important problems  arise.

Indeed, recently there has been much work in the causal inference literature focused on understanding the behavior and properties of various estimators in the high-dimensional setting. \cite{yadlowsky2022explaining} 
derived the excess variance of commonly used estimators due to first-stage nuisance estimation with high-dimensional covariates by assuming a linear outcome model and known propensity score. They also illustrated the inflation of variance in simulations, showing even the doubly robust estimator may not achieve the efficiency bound when the covariates included are high-dimensional. \cite{jiang2022new} established a novel central limit theorem for the doubly robust estimator, when the outcome regression model and propensity score follow high-dimensional generalized linear models, without assumptions on sparsity, but under a stylized setting where covariates are normally distributed and the covariate dimension $d$ is of the same order as sample size $n$. They also showed estimates obtained by permuting the folds in cross-fitting are asymptotically correlated in the high-dimensional regime. \cite{celentano2023challenges} proposed a novel debiased method for missing data models in the $n < d$ setting, where ordinary least square estimation is not feasible. However, to help overcome the curse of dimensionality and achieve non-trivial rates of convergence, most of the work in this literature combines the debiased machine learning/doubly robust techniques with additional structural assumptions on the nuisance functions, such as the aforementioned linearity, or smoothness \citep{robins2009quadratic, kennedy2023towards, semenova2021debiased}, or sparsity \citep{belloni2017program, chernozhukov2018double, athey2018impact, bradic2019minimax}.   We refer to \cite{maathuis2009estimating, lin2013agnostic, zhao2016topics, chakrabortty2018inference, ma2019robust, lei2021regression, antonelli2022causal, tang2023ultra, du2024causal} and others for more relevant discussion on causal inference with high-dimensional data. 

Our work enriches the high-dimensional causal literature under a different structural assumption, i.e., that the covariates are discrete. This setting is of interest for several reasons. First, it can help uncover new phenomena in more general high-dimensional regimes where the dimension $d$ can be comparable to or larger than sample size $n$, and motivate further exploration. 
The discrete covariate setting can also be viewed as an interesting base case, with crucial implications for continuous data, as has been seen in semiparametric efficiency bounds \citep{chamberlain1992efficiency},  bandit problems \citep{kleinberg2004nearly}, and hypothesis testing  \citep{balakrishnan2019hypothesis}, for example. Further, discrete high-dimensional covariates often arise in practice, for example  in applied medical and health policy research, e.g., when adjusting for International Classification of Diseases (ICD) codes \citep{world2004international}. ICD codes include over 100,000 indicators, so their cartesian product, together with other common demographic covariates (e.g., sex, race, education level), induces a huge number of categories, certainly corresponding to a high-dimensional regime. The structured discrete data, including graphs, texts and images, is also ubiquitous in the machine learning literature, which commonly has high-dimensional representations. It is often important to take these non-numerical variables into account when estimating treatment effects \citep{yao2019estimation, keith2020text}. Our work establishes theoretical foundations for adjusting for high-dimensional discrete covariates and helps practitioners evaluate the reliability of their estimates given the sample size and covariate dimension.

Our work is also related to \citet{paninski2003estimation, valiant2010clt, jiao2015minimax, wu2016minimax, wu2019chebyshev}, and others, who considered estimating simpler functionals (e.g., entropy, support size) in high-dimensional discrete models. In this work, important breakthroughs have been made in terms of both improved estimation methodologies and minimax fundamental limits. Using tools in polynomial approximation theory \citep{timan2014theory}, it has been established that a best-polynomial-approximation estimator usually enjoys the so-called ``effective sample size enlargement" property, meaning its behavior with $n$ samples resembles that of an MLE (the naive plug-in estimator) with $n \log n$ samples. Lower bounds are established by a moment matching technique. 
Upper bounds and lower bounds using moment matching usually coincide here since, in the optimization sense, moment matching is the dual problem of best polynomial approximation where strong duality holds. The readers are referred to \cite{wu2016minimax, luenberger1997optimization} for more detailed discussion on this duality phenomenon.  Estimation of causal effects under the discrete covariate framework has remained largely unexplored, but our work helps bridge this gap.

Specifically, in this paper we study the treatment effect estimation problem when adjusting for discrete covariates with possibly more categories than samples. Our four main contributions are summarized as follows:
\begin{enumerate}
    \item First, we provide finite-sample bounds on the mean squared error of commonly used regression, weighting and doubly robust estimators of ATE estimation. Our results imply $d/n \rightarrow 0$ is a sufficient and necessary condition for these estimators to be consistent under positivity assumption, where $d$ is the number of categories and $n$ is the sample size. 
    \item We then characterize the fundamental limits of ATE estimation under a positivity assumption. The minimax lower bound is (in terms of mean squared error) of order $\frac{d^2}{n^2 \log^2 n} + \frac{1}{n}$, which shows that commonly used estimators are minimax optimal up to log factors. Moreover, $d=o(n \log n)$ is a necessary condition for consistent ATE estimation.
    \item Next, we explore the role of effect homogeneity in ATE estimation, showing that faster rates are achievable when the treatment effects are more homogeneous across categories. In fact, here consistent estimation is possible if $d/n^2 \rightarrow 0$. 
    \item Finally we study estimation of treatment effects given prior knowledge of the covariate distribution. We first provide a negative result, showing the covariate distribution may not help improve the rate for ATEs. We then consider a variance-weighted average treatment effect, and derive its faster estimation rate, achieved by a second-order estimator. As in the homogeneous effects setting, our second-order estimator is consistent here if $d/n^2 \rightarrow 0$.
\end{enumerate}
The structure of our paper is as follows: After introducing the data-generating process, causal assumptions and notation in Section \ref{sec:setup}, in Section \ref{sec:ate-plugin} we study the properties of commonly used regression, weighting and doubly robust estimators in estimating the average treatment effects (ATE). We first show the numerical equivalence between these three estimators and the plug-in estimator. 
We then provide finite-sample bounds on their mean squared error and characterize the sufficient and necessary conditions for them to be consistent in Section \ref{sec:ate-upper}. In Section \ref{sec:minimax} we study the minimax lower bound and sample complexity of ATE estimation in the high-dimensional setting, based on the moment matching method borrowed from the theoretical computer science literature. Next, we consider two additional structures that one can exploit in the high-dimensional problem: effect homogeneity (Section \ref{sec:homogeneity}) and prior knowledge of covariate distribution (Section \ref{sec:cova-distribution}). We propose novel estimators that can take advantage of these structures and enjoy faster convergence rates, making consistent estimation possible in the regime where conventional estimators examined in Section \ref{sec:ate-plugin} fail to be consistent. Finally in Section \ref{sec:simulation} we perform numerical experiments to verify our results in previous sections. All the proofs and additional complementary results are presented in the appendix. To the best of our knowledge, our work is the first in the literature that directly analyzes the theoretical properties of different estimators of causal effects, provides a minimax lower bound for ATE estimation, and explores additional structures that we can exploit to achieve faster convergence rates, in the high-dimensional discrete setting. 

\section{Setup and Notation}\label{sec:setup}
In this section, we first introduce the data-generating process in the covariate setting and characterize the distributions of several counting statistics that play an important role in estimating causal effects. Identification assumptions of causal estimands and additional notation are further introduced with discussion.

\subsection{Data Generating Process}\label{sec:DGP}
Suppose we observe $n$ i.i.d. copies of $\bZ =(X,A,Y)$ where $X$ is the discrete covariate with $d$ categories, $A \in\{0,1\}$ is the binary treatment and $Y \in\{0,1\}$ is the binary outcome. On the population level, the covariate $X$ has a categorical distribution on $[d]=\{1,\dots,d\}$ with
\[
\MP(X=k) = p_k, \,1 \leq k \leq d.
\]
Let $\bp=(p_1,\dots,p_d)$ be the probability vector. In real applications one may observe multiple discrete covariates. For example, we may observe $K$ binary covariates and it is easy to see that one could encode these binary variables as one categorical variable with $d=2^K$ categories.
Hence we will assume a single discrete covariate $X$ in our problem. Given the covariate $X=k$, the treatment $A$ follows a Bernoulli distribution with parameter $\pi_k$
\[
A \mid X=k\, \sim \text{Bernoulli} (\pi_k)  \text{ with } \pi_k  =\mathbb{P}(A=1 \mid X=k).
\]
Conditioned on the covariate $X=k$ and the treatment $A=a$, $Y$ has a Bernoulli distribution with parameter $\mu_{ak}$
\[
Y\mid A=a, X=k \, \sim \text{Bernoulli} (\mu_{a k}) \text{ with } \mu_{a k}  =\mathbb{P}(Y=1 \mid X=k, A=a).
\]
$\pi_k$ and $ \mu_{ak}$ are the propensity score and regression functions, respectively. We further denote $q_{a k}=\mathbb{P}(X=k, A=a, Y=1)=p_k\pi_k^a\left(1-\pi_k\right)^{1-a}\mu_{a k}$ and $w_k=\mathbb{P}(X=k, A=1)= p_k \pi_k$. Based on a sample of size $n$, define the following empirical average estimators of the parameters:
\[
n \widehat{q}_{a k}= \, \sum_{i=1}^n I\left(X_i=k, A_i=a, Y_i=1\right) \sim \operatorname{Bin}\left(n, q_{a k}\right),
\]
\[
n \hat{w}_k =\,   \sum_{i=1}^n I\left(X_i=k, A_i=1\right) \sim \operatorname{Bin}\left(n, w_k\right),
\]
\[
n \hat{p}_k =\,  \sum_{i=1}^n I\left(X_i=k\right) \sim \operatorname{Bin}\left(n, p_k\right)
\]
for each $1 \leq k \leq d$. The corresponding empirical average estimates of the propensity score $\pi_k$ and regression function $\mu_{ak}$ are 
\begin{equation}\label{eq:mle-nuisance}
    \begin{aligned}
    \hat{\pi}_{k} = &\, \frac{\hat{w}_k}{\hat{p}_k} = \frac{\#\{i: X_i=k, A_i=1\}}{\#\{i: X_i=k\}},  \\
    \hat{\mu}_{1k} = &\, \frac{\hat{q}_{1k}}{\hat{w}_k} = \frac{\#\{i: X_i=k, A_i=1, Y_i=1\}}{\#\{i: X_i=k,A_i=1\}}, \\
    \hat{\mu}_{0k} = &\, \frac{\hat{q}_{0k}}{\hat{p}_k-\hat{w}_k} = \frac{\#\{i: X_i=k, A_i=0, Y_i=1\}}{\#\{i: X_i=k,A_i=0\}},
    \end{aligned}
\end{equation}
where we define $0/0=0$ whenever both the numerator and denominator are zero. This may happen when no individual in $k$-th category is assigned to treatment ($\hat{w}_k = 0$) and hence the response $Y$ under treatment is unavailable ($\hat{q}_{1k} = 0$). Let $\bX^n = (X_1 ,\dots, X_n)$ and $\bA^n = (A_1, \dots, A_n)$ be the collection of covariates and treatments in the sample, respectively. According to the sampling schemes we have
\[
(n\hat{p}_1, \dots,n\hat{p}_d ) \, \sim \text{Multinomial} (n, p_1, \dots, p_d),
\]
\[
n \widehat{w}_k=n \widehat{p}_k \widehat{\pi}_k \mid \bX^n \sim \operatorname{Bin}\left(n \widehat{p}_k, \pi_k\right),
\]
\[
\begin{aligned}
   &\,n \widehat{q}_{a k}= n \widehat{p}_k \left[a \widehat{\pi}_k+(1-a)\left(1-\widehat{\pi}_k\right)\right]  \widehat{\mu}_{a k} \mid \bX^n, \bA^n \\
   &\,\sim \operatorname{Bin}\left(n \widehat{p}_k \left[a \widehat{\pi}_k+(1-a)\left(1-\widehat{\pi}_k\right)\right] , \mu_{a k}\right). 
\end{aligned}
\]
Moreover, we have $\widehat{w}_k \independent \widehat{w}_{\ell} \mid \bX^n$ for $k \neq \ell$ since conditioned on the number of samples in each category, the treatment assignment within different categories $X=k,\ell (k \neq \ell)$ proceeds independently. Similarly, conditioned on the number of samples in each category and treatment assignment, the outcomes within different categories are conditionally independent, i.e., $\widehat{q}_{a k} \independent \widehat{q}_{a \ell} \mid \bX^n, \bA^n$ for $k \neq \ell$. 

\subsection{Causal Assumptions and Other Notation}

To properly define and identify causal estimands of interests, we rely on the potential outcome framework \citep{rubin1974estimating, splawa1990application} and additional identification assumptions to connect counterfactual outcomes with the observed data. We use the random variable $Y^a$ to denote the potential/counterfactual outcome we would have observed
had a subject received treatment $A = a$, which may be contrary to the observation $Y$. The average treatment effect (ATE) $\psi$ is defined as 
\[
\psi = \ME\left[Y^1-Y^0\right].
\]
The following identification assumptions are often imposed to identify $\psi$ as a functional of the observed data distribution $\MP$.
\begin{assumption}\label{as-consistency}
  Consistency: 
$Y = Y^a \text{\, if \,} A=a.$
\end{assumption}

\begin{assumption}\label{as-positivity}
Positivity: For any $k \in [d]$, $\pi_k \in [\epsilon, 1-\epsilon]$ for some constant $\epsilon \in (0,1/2)$.
\end{assumption}

\begin{assumption}\label{as-exchangeability}
No unmeasured confounding: 
$Y^a \independent A \mid X$ for $a=0,1$.
\end{assumption}
Under Assumptions \ref{as-consistency}--\ref{as-exchangeability}, the ATE can be identified as
\begin{equation}\label{eq:ATE}
    \begin{aligned}
      \psi = &\, \ME\left[\ME(Y\mid X,A=1)-\ME(Y\mid X,A=0)\right]  \\
      = &\, \sum_{k=1}^d p_k(\mu_{1k}- \mu_{0k}) =\sum_{k=1}^d p_k\left(\frac{q_{1k}}{w_k}- \frac{q_{0k}}{p_k-w_k}\right).  
    \end{aligned} 
\end{equation}
We refer the readers to \cite{hernan2010causal} for detailed discussion on these identification assumptions. 
From this point forward, $\psi$ will denote the functional of observed data distribution in \eqref{eq:ATE}, which is equal to ATE when Assumptions \ref{as-consistency}--\ref{as-exchangeability} hold. If Assumption \ref{as-consistency}--\ref{as-exchangeability} are violated, the functional $\psi$ should only be viewed as the expected difference in the regression functions between the treatment and control group, which may not represent a causal effect. Nonetheless, all our results still hold for $\psi$ in \eqref{eq:ATE} under only the positivity assumption. 

In this paper we will consider the following model in which the number of categories for the covariate is at most $d$:
\begin{equation}\label{eq:model-D}
   \mathcal{D}(\epsilon) = \left\{ \sum_{k=1}^d p_k = 1, 0 \leq p_k\leq 1  , \epsilon \leq \pi_k \leq 1-\epsilon, 0\leq \mu_{0k}, \mu_{1k} \leq 1 ,\, \forall  k \in [d]  \right \}.
\end{equation}
We will always assume positivity in addition to the basic bounds on model parameters. The binary nature of the outcome variable $Y$ is assumed primarily for simplicity and is not essential to our analysis. Our results on rates remain valid as long as $\mathbb{E}[Y^2 \mid X = k, A = a] \leq C$ for some constant $C > 0$. This condition holds if, within each covariate level $X = k$ and treatment arm $A = a$, the outcome $Y$ has bounded variance. For example, this is satisfied when $Y$ is bounded or sub-Gaussian with a uniformly bounded sub-Gaussian norm given $X=k$ and $A=a$. 

\begin{remark}
It is worth noting that in the high-dimensional regime where $d$ can be large (the regime we are mainly interested in), the positivity assumption \ref{as-positivity} can impose additional restriction on the observed distribution \citep{d2021overlap}. However, when positivity is violated, the ATE is not identified and may not be an appropriate estimand to focus on. This together with other possible issues due to violation of positivity are tangential to our main points and thus we may proceed assuming positivity holds. Some of our results characterize the rate in terms of $\epsilon$ explicitly (e.g., Theorem~\ref{thm:ate-var}) and still provide meaningful implications under weak positivity assumption where $\epsilon = \epsilon_n $ shrinks to zero. 
\end{remark}

For a (possibly random) function $f$ of the observation $\bZ = (X,A,Y)$, we use $\MP_n[f(\bZ)]$ to denote the sample average $\frac{1}{n} \sum_{i=1}^n f(\bZ_i)$, $\MP[f(\bZ)] = \int f(\bz) d\MP(\bz)$, and $\|f\|_2$ to denote the $L_2$-norm $\left[\int f^2(\bz) d \MP(\bz)\right]^{1/2}$, where all expectations are only taken with respect to the randomness of $\bZ$. For a bivariate function $g$ on $(\bZ_1,\bZ_2)$, let $\MU_n[g(\bZ_1,\bZ_2)]$ denote the second-order U-statistic measure $\frac{1}{n(n-1)} \sum_{i \neq j} g(\bZ_i,\bZ_j)$. For two sequences $\{a_n\}$ and $ \{b_n\}$, we use $a_n \lesssim b_n$ to denote that there exists a constant $C>0$ such that $a_n \leq C b_n$ when $n$ is sufficiently large, and $a_n \gtrsim b_n$ means $b_n \lesssim a_n$. We use $a\vee b$ and $a \wedge b$ to denote the maximum and minimum of $a$ and $b$, respectively.

\section{Average Treatment Effects}\label{sec:ate-plugin}
In this section, we focus on the theoretical properties of commonly used estimators of ATE
\[
\psi = \ME\left[\ME(Y\mid X,A=1)-\ME(Y\mid X,A=0)\right] = \sum_{k=1}^d p_k(\mu_{1k}- \mu_{0k}) .
\]
A simple plug-in-style estimator based on empirical average estimates of model parameters $(p_k,\pi_k,\mu_{1k}, \mu_{0k})$ in  \eqref{eq:mle-nuisance} is 
\begin{equation}\label{eq:plugin}
    \hat{\psi} = \sum_{k=1}^d \hat{p}_k \left(\hat{\mu}_{1k}-\hat{\mu}_{0k} \right)= \sum_{k=1}^d \hat{p}_k\left(\frac{\hat{q}_{1k}}{\hat{w}_k}- \frac{\hat{q}_{0k}}{\hat{p}_k-\hat{w}_k}\right), 
\end{equation}
We also consider three popular estimators in the literature to estimate ATE, namely the outcome regression estimator $\hat{\psi}_{\text{reg}}$ \citep{rubin1979using}, inverse probability weighting estimator $\hat{\psi}_{\text{ipw}}$ \citep{rosenbaum1983central, hahn1998role} and doubly robust estimator $\hat{\psi}_{\text{dr}}$ \citep{robins1994estimation, scharfstein1999adjusting} defined as follows:
\begin{equation}\label{eq:ate-ests}
    \begin{aligned}
       \hat{\psi}_{\text{reg}} = &\, \MP_n\left[\hat{\mu}_{1X} - \hat{\mu}_{0X} \right], \\
       \hat{\psi}_{\text{ipw}} = &\, \MP_n\left[\frac{AY}{\hat{\pi}_X} - \frac{(1-A)Y}{1-\hat{\pi}_X} \right], \\
       \hat{\psi}_{\text{dr}} = &\, \MP_n\left[\frac{A(Y-\hat{\mu}_{1X})}{\hat{\pi}_X} + \hat{\mu}_{1X} - \frac{(1-A)(Y-\hat{\mu}_{0X})}{1-\hat{\pi}_X}- \hat{\mu}_{0X}\right],
    \end{aligned}
\end{equation}
where again in the IPW and DR estimator we define $0/0=0$ whenever it occurs. In this paper, we will focus on the observational study setting where propensity score is unknown and needs to be estimated. In randomized experiments where the treatment process is known, the IPW estimator with known propensity scores is unbiased and $\sqrt{n}$-consistent. 

Surprisingly, all these three commonly used estimators of ATE are numerically equivalent to the simple plug-in estimator \eqref{eq:plugin} in the discrete covariate setting, if we use the same sample to estimate the nuisance parameters in \eqref{eq:mle-nuisance} and take average over in \eqref{eq:ate-ests}. This property unique to the discrete covariate setting is presented in regression coefficients estimation with missing covariates in \cite{wang2007numerical}. Recently similar equivalence in ATE estimation is shown in \cite{sloczynski2023covariate} when the correct parametric model for propensity score is specified. The discrete covariate setting can be viewed as a special case where the design matrix only contains indicators specifying the category membership of each sample. We restate their result in the discrete covariate setting and emphasize that such equivalence still holds in the high-dimensional setting where some categories may not have treated/untreated samples observed, as long as we define $0/0 = 0$ whenever it appears.

\begin{prop}\label{prop:ests-equivalence}
        Supposed $\mathcal{D}=\{(X_i,A_i,Y_i), 1 \leq i \leq n\}$ is a sample of size $n$ and $X \in [d]$ is discrete. If the nuisance estimators are the empirical averages defined in equation \eqref{eq:mle-nuisance} using $\mathcal{D}$ and we take average over $\mathcal{D}$ in \eqref{eq:ate-ests}. Then we have
    \[
    \hat{\psi}_{\text{reg}} = \hat{\psi}_{\text{ipw}} = \hat{\psi}_{\text{dr}} = \hat{\psi}.
    \]
    So the regression, weighting and doubly robust estimators are numerically equivalent in the discrete covariate setting.
\end{prop}
The numerical equivalence in Proposition \ref{prop:ests-equivalence} does not necessarily hold when the covariates have continuous components and smoothing is applied to construct estimates of propensity scores and regression functions. 
One explanation is that in the discrete case, the arguments used to show $\ME [{\mu}_{1X} ] = \ME[AY/\pi_X]=\psi_1$ (See e.g., Chapter 2 in \cite{hernan2010causal}) can be applied to the empirical distribution $\MP_n$ (i.e. replace every expectation and conditional expectation in the argument with sample averages) to show $\MP_n [\hat{\mu}_{1X}] = \MP_n[AY/\hat{\pi}_X]=\hat{\psi}$.

Proposition \ref{prop:ests-equivalence} shows that we only need to consider the properties of plug-in-style estimator $\hat{\psi}$ and all the results hold for estimators in \eqref{eq:ate-ests} as well. In the following discussion, we will first derive the upper bound on the mean squared error of $\hat{\psi}$ in Section \ref{sec:ate-upper}. Then we characterize the minimax lower bound in Section \ref{sec:minimax}. 

\subsection{Upper Bound }\label{sec:ate-upper}

In this section, we study the behavior of $\hat{\psi}$ in a potentially high-dimensional regime and characterize its mean squared error in terms of $(n,d)$. In the low-dimensional regime where $d$ is fixed, the empirical average estimators of nuisance functions \eqref{eq:mle-nuisance} have parametric convergence rates. One can expect this ideal property to propagate to plug-in-style estimator $\hat{\psi}$, which crucially depends on the empirical average estimates. We provide a central limit theorem of $\hat{\psi}$ in the appendix when $d$ is fixed to avoid distracting the readers from the high-dimensional regime of main interest. In the following discussions, our main results hold non-asymptotically for all $n,d$ sufficiently large. We will focus on the point estimation theory in this work and characterize the condition under which $\hat{\psi}$ is consistent. Let $\psi_a = \ME[Y^a] = \ME[\ME(Y\mid X,A=a)]= \sum_{k=1}^d p_k \mu_{ak}$ be the mean of potential outcome $Y^a$ and $\hat{\psi}_a = \sum_{k=1}^d \hat{p}_k \hat{\mu}_{ak}$ be the corresponding plug-in-style estimator. 
We first summarize the exact bias of $\hat{\psi}_a$ in the following proposition.

\begin{prop}\label{thm:ate-exact-bias}
The exact bias of $\hat{\psi}_a$ is
\begin{equation*}
    \begin{aligned}
        \ME\left[\hat{\psi}_1 - \psi_1\right] =&\,  -\sum_{k=1}^d \mu_{1k}p_k(1-\pi_k)(1-p_k \pi_k)^{n-1}, \\
        \ME\left[\hat{\psi}_0 - \psi_0\right] =&\,  -\sum_{k=1}^d \mu_{0k}p_k\pi_k(1-p_k +p_k\pi_k)^{n-1}. \\
    \end{aligned}
\end{equation*}
Hence the exact bias of $\hat{\psi}=\hat{\psi}_1 - \hat{\psi}_0$ is
\[
\begin{aligned}
  &\,\ME\left[\hat{\psi} - \psi\right] \\
  = &\,- \sum_{k=1}^d \left[ \mu_{1k}p_k(1-\pi_k)(1-p_k \pi_k)^{n-1} -  \mu_{0k}p_k\pi_k(1-p_k +p_k\pi_k)^{n-1}\right] 
\end{aligned}
\]
\end{prop}

From the proof of Proposition \ref{thm:ate-exact-bias}, the bias of $\hat{\psi}_a$ comes from categories with no subjects receiving treatment $A=a$, and hence it is not possible to obtain an unbiased estimator of regression functions $\mu_{ak}$ within these categories. It is unclear from the exact bias how fast $d$ can grow with $n$ while still having vanishing bias. To make the dependency of bias on $(n,d)$ explicit, we maximize the bias over model $\mathcal{D}(\epsilon)$ and characterize the worst-case bias in the next proposition. In the following discussion of this section, we will mainly consider the functional $\psi_1$, with the understanding that analogous arguments can be applied to $\psi_0$ to obtain the same rate/property. 

\begin{prop}\label{prop:ate-worst-bias}
The worst-case bias of $\hat{\psi}_1$ is
\[
\sup_{\MP \in \mathcal{D(\epsilon)}}\big |\ME_\MP[\hat{\psi}_1 - \psi_1]\big| = \sup_{\bp} \sum_{k=1}^d p_k (1-\epsilon)(1-\epsilon p_k )^{n-1}.
\]
As a consequence, we have the following bounds on the worst-case bias:
\[
\frac{1-\epsilon}{2e}\left(\frac{d-1}{\epsilon n}\wedge 1\right) \leq \sup_{\MP \in \mathcal{D(\epsilon)}} \big|\ME_\MP[\hat{\psi}_1 - \psi_1]\big| \leq \left( \frac{1-\epsilon}{\epsilon} \right) \frac{d}{n}.
\]
\end{prop}
Proposition \ref{prop:ate-worst-bias} shows that a sufficient and necessary condition for the bias to vanish as $n \rightarrow \infty$ in the worst case over model class $\mathcal{D}(\epsilon)$ is $d/n \rightarrow 0$. The following theorem further characterizes a bound on the variance of $\hat{\psi}_1$ and arrives at the MSE of the plug-in-style estimator. 

\begin{theorem}\label{thm:ate-var}
    The variance of $\hat{\psi}_1$ is upper bounded as
    \[
    \sup_{\MP \in \mathcal{D}(\epsilon)} \operatorname{Var}_{\MP}(\hat{\psi}_1) \leq \frac{C}{\epsilon n},
    \]
    where $C>0$ is an absolute constant. Hence the MSE is upper bounded as
    \[
    \sup_{\MP \in \mathcal{D}(\epsilon)} \ME_\MP\left[\left(\hat{\psi}_1 - \psi_1\right)^2\right] \leq \frac{(1-\epsilon)^2}{\epsilon^2}\frac{d^2}{n^2}  + \frac{C}{\epsilon n}.
    \]
\end{theorem}
Under positivity condition \ref{as-positivity}, the marginal probability for a subject to receive treatment is $\MP(A=1) \geq \epsilon$. Thus on average the total number of treated samples is lower bounded by $\epsilon n$. Since the denominator of the bound on variance is exactly $\epsilon n$, intuitively $\epsilon n$ acts as the effective sample size for estimating $\psi_1$. The results in Proposition \ref{prop:ate-worst-bias} and Theorem \ref{thm:ate-var} imply:
\begin{itemize}
    \item The plug-in estimator $\hat{\psi}_1$ is consistent in non-classical high-dimensional regimes where $d\rightarrow \infty$ as long as $d/n \rightarrow 0$ holds.
    \item $d/n \rightarrow 0$ is also a necessary condition for the worst-case bias to vanish as $n \rightarrow \infty$. Thus consistency of $\hat{\psi}_1$ in the high-dimensional regime $n \lesssim d$ is not achievable without further assumptions.
\end{itemize}
The intuition is that with $n$ samples and $d$ categories, there are $n/d$ samples within each category on average. Without further assumptions, one needs to consistently estimate the regression functions $\mu_{1k}$ in each category to achieve consistent estimation of $\psi_1 = \sum_k p_k \mu_{1k}$. And consistent estimation of $\mu_{1k}$ requires infinite samples assigned to each category, i.e. $n/d \rightarrow \infty$. 

\begin{remark}
    When the covariate distribution is sparse in the sense that the support size $s := |{k : p_k > 0}| \ll d$, we can replace $d$ with $s$ in the bound of Theorem \ref{thm:ate-var}, yielding a rate of $s^2/n^2 + 1/n$. Thus, consistency is achieved as long as $s/n \to 0$.
\end{remark}

\begin{remark}
    Note that in the estimators above (including $\hat{\psi}_{\text{reg}}, \hat{\psi}_{\text{ipw}}, \hat{\psi}_{\text{dr}}$, and $\hat{\psi}$), the same sample is used to estimate the nuisance functions and to compute the final estimator. Alternatively, one can adopt a sample splitting strategy, as in the double machine learning literature \citep{chernozhukov2018double}, where a separate sample is used for nuisance estimation. In this case, the numerical equivalence in Proposition \ref{prop:ests-equivalence} no longer holds. However, sample splitting does not necessarily improve the convergence rate or guarantee consistency in the high-dimensional regime where $n \lesssim d$. Specifically, assume the nuisance functions are estimated from an independent sample of size $n$ using empirical averages, which are natural estimators in the discrete covariate setting without additional structural assumptions. Following the proof of Theorem \ref{thm:cova-dist}, we have
\[
\|\hat{\mu}_{aX}-\mu_{aX}\|_2=O_{\MP}\left(\sqrt{\frac{d}{n}} \right), \, \|\hat{\pi}_{X}-\pi_{X}\|_2=O_{\MP}\left(\sqrt{\frac{d}{n}} \right),
\]
where for a function $g_X$, we define $\|g_X\|_2^2 = \sum_k p_k g_k^2$. The squared errors of the regression estimator $\hat{\psi}_{\text{reg}}$, the IPW estimator $\hat{\psi}_{\text{ipw}}$, and the doubly robust estimator $\hat{\psi}_{\text{dr}}$ are therefore dominated by $\|\hat{\mu}_{aX}-\mu_{aX}\|_2^2=O_{\mathbb{P}}(d/n)$, $\|\hat{\pi}_{X}-\pi_{X}\|_2^2=O_{\mathbb{P}}(d/n)$, and $\|\hat{\mu}_{aX}-\mu_{aX}\|_2^2\|\hat{\pi}_{X}-\pi_{X}\|_2^2=O_{\mathbb{P}}(d^2/n^2)$, respectively. As a result, consistency of these estimators still requires $d/n \to 0$, even when sample splitting is employed. 
\end{remark}

However, our results in this section do not preclude the existence of other consistent estimates. In order to conclusively determine whether consistent estimates exist in the high-dimensional regime, one needs to characterize the minimax lower bound for $\psi_1$.

\subsection{Minimax Lower Bound}\label{sec:minimax}



In Section \ref{sec:ate-upper} we showed that $d/n \rightarrow 0$ is a sufficient and necessary condition for the plug-in estimator $\hat{\psi}$ to be consistent. In this section, we study the existence of consistent estimators in the high-dimensional regime by considering the minimax lower bound for the mean of the regression function $\psi_1 = \ME[\ME(Y\mid X, A=1)]$, with the understanding that similar arguments show $\psi_0=\ME[\ME(Y\mid X, A=0)]$ and ATE share the same minimax rate as $\psi_1$. Recall the model we consider is:
\[
\mathcal{D}(\epsilon) = \left\{ 0 \leq p_k\leq 1  , \epsilon \leq \pi_k \leq 1-\epsilon, 0\leq \mu_{0k}, \mu_{1k} \leq 1 ,\, \forall k \in [d]  \right \},
\]
which corresponds to the setting where statisticians have knowledge on a set of possible values of $X$ but some categories may have zero proportion.
By Theorem \ref{thm:ate-var}, the MSE of plug-in-style estimator satisfies
\begin{equation}\label{eq:plugin-p}
\ME_\MP\left[\left(\hat{\psi}_1 - \psi_1\right)^2\right] \lesssim \frac{d^2}{n^2}  + \frac{1}{ n}, \forall \, \MP \in \mathcal{D} (\epsilon),  
\end{equation}
which shows a sufficient condition for the plug-in estimator to be consistent is $d/n \rightarrow 0$, i.e. we require sample size $n \gg d$ to achieve consistency. The following theorem characterizes the minimax lower bound of the estimation error of $\psi_1$ over model class $\mathcal{D} (\epsilon)$.


\begin{theorem}\label{thm:ate-minimax}
    In the regime $ d \lesssim n \log n$, we have
    \[
    R^*(d,n) := \inf_{\hat{\psi}_1} \sup_{\MP \in \mathcal{D} (\epsilon)}  \ME_{\MP}\left[\left(\hat{\psi}_1 - \psi_1\right)^2\right] \gtrsim \left( \frac{d}{n \log n}\right)^2 + \frac{1}{n},
    \]
    where the constant behind ``$\gtrsim$" depends on $\epsilon$.
\end{theorem}

The proof adopts the idea of moment matching commonly used in deriving minimax lower bounds for functionals of discrete distributions \citep{jiao2015minimax, wu2016minimax, wu2019chebyshev}. Specifically, we need to show that there exist two probability measures $\mu_0, \mu_1$ over the tuple $(p, \pi, \mu )$ satisfying
\begin{align}
a_0 := \ME_{\mu_0}[p] &=  \ME_{\mu_1}[p] \leq \frac{1}{d-1}, \label{eq:mean-constraint} \\
\ME_{\mu_0}\left[ p^i (p \pi)^j (p \pi \mu)^k \right] &= \ME_{\mu_1} \left[ p^i (p \pi)^j (p \pi \mu)^k \right], \quad \forall i, j ,k \geq 0, \,i+j+k \leq K, \label{eq:moment-matching} \\
 |\ME_{\mu_0}[p \mu] - \ME_{\mu_1}[p \mu]| &\gtrsim \frac{1}{n \log n},\label{eq:separation}
\end{align}
with $K \asymp \log n$. Once the measures $(\mu_0, \mu_1)$ are constructed, we construct two ``difficult" hypotheses $(H_0, H_1)$, where under the hypothesis $H_u$ with $u\in \{0,1\}$, let $(p_k, \pi_k, \mu_{1k})\overset{\text{i.i.d.}}{\sim} \mu_u$ for $ 1\leq k \leq d-1$, and $(p_{d}, \pi_{d}, \mu_{1d})=(1-(d-1)a_0, \epsilon, 0)$. On a high level, the mean constraint \eqref{eq:mean-constraint} ensures that $\bp$ is an approximate probability measure with high probability, the moment matching constraint \eqref{eq:moment-matching} ensures that $H_0$ and $H_1$ are statistically indistinguishable based on the observations $(\bX^n, \bA^n, \bY^n)$, and the separation condition \eqref{eq:separation} ensures that the value of $\psi_1 = \sum_{k=1}^d p_k \mu_{1k}$ is separated apart by an amount of $\Omega(d/(n\log n))$ under hypotheses $H_0$ and $H_1$. Based on these intuitions, the minimax lower bound in Theorem \ref{thm:ate-minimax} follows from the method of fuzzy hypothesis \citep{Tsybakov2009}.

To show the existence of measures $(\mu_0, \mu_1)$ satisfying \eqref{eq:mean-constraint}, \eqref{eq:moment-matching}, and \eqref{eq:separation}, by the duality of moment matching and best polynomial approximation (cf. e.g. \citep[Sec. 4.3]{lepski1999estimation}), the maximum separation in \eqref{eq:separation} subject to the moment constraint \eqref{eq:moment-matching} is characterized by
\begin{align}\label{eq:3D-approx-error}
\inf_{Q\in \mathbb{R}[x,y,z], \deg Q\le K} \sup_{(p,p\pi,p\pi\mu)\in D} |Q(p,p\pi,p\pi\mu) - p\mu|, 
\end{align}
where the approximation domain $D$ is a 3-dimensional polytope given by
\begin{align*}
D = \{(x,y,z): 0\le x\le 1, \epsilon x\le y\le (1-\epsilon)x, 0\le z\le y \}. 
\end{align*}
Characterizing the best 3D polynomial approximation error \eqref{eq:3D-approx-error} is generally very challenging \citep{rice1963tchebycheff}, and we lower bound \eqref{eq:3D-approx-error} via a proper one-dimensional subproblem, carefully chosen as
\begin{equation}\label{eq:1D-approx}
    \inf_{P \in  \text{span}\{1/x, 1, x, \dots, x^{K}\} } \max_{x \in [c/K^2,1]} \left | \frac{x}{x+c/K^2}-P(x) \right|. 
\end{equation}
The intuition behind the subproblem \eqref{eq:1D-approx} is a one-dimensional trajectory $x\mapsto (x,\epsilon x+c',\epsilon x)\in D$ (roughly) parallel to an edge of $D$, which is further motivated by the approximation-theoretic results in \citep[Chap. 12]{ditzian2012moduli} and \citep{totik2014polynomial}. Lower bounding \eqref{eq:1D-approx} using machineary in \citep{ditzian2012moduli} then resolves \eqref{eq:moment-matching} and \eqref{eq:separation}. As for the mean constraint \eqref{eq:mean-constraint}, we overcome it by a change-of-measure trick, which is the reason behind the additional basis $1/x$ and the constraint $x\ge c/K^2$ in \eqref{eq:1D-approx}. We defer the details to the appendix.

Theorem \ref{thm:ate-minimax} shows in the regime $ d \gtrsim n \log n$, a consistent estimator of $\psi_1$ does not exist over model class $\mathcal{D} (\epsilon)$. Compared with the upper bound in \eqref{eq:plugin-p}, we conclude that up to log-factors the plug-in estimator $\hat{\psi}_1$ is minimax optimal and one needs sample size at least of order $d$ to consistently estimate $\psi_1$ over model class $\mathcal{D}(\epsilon)$. In applications, if the observed number of categories of the discrete covariate $X$ is comparable to sample size $n$, then we should interpret the estimated ATE carefully since it is possible that our estimator is not consistent in that high-dimensional regime. 

It is worth noting that the lower bound in Theorem \ref{thm:ate-minimax} does not exactly match the upper bound provided by the plug-in-style estimator with the difference being a log-factor. It is possible that some estimator of $\psi_1$ based on polynomial approximation, which further reduces the bias, enjoys the ``effective sample size enlargement" property \citep{wu2016minimax, jiao2015minimax, wu2019chebyshev} and achieves the minimax lower bound.  We leave the exploration of such estimators to future investigation.

In the following two sections, we consider how effect homogeneity and prior knowledge of the covariate distribution can allow faster rates and consistent estimation of causal effects in the regime $n \lesssim d$ under certain scaling conditions on $n$ and $d$.

\section{Role of Effect Homogeneity}\label{sec:homogeneity}

In this section, we study the role of effect homogeneity in consistent estimation of treatment effects in the regime $n \lesssim d$. Let $\tau_k = \ME[Y|X=k, A=1] - \ME[Y|X=k, A=0] = \mu_{1k}-\mu_{0k}$ be the conditional average treatment effect (CATE) in the $k$-th category. In the high-dimensional regime where $d=d_n$ can grow with $n$, $\tau_k$'s should be viewed as a triangular array $\{\tau_{nk}, 1\leq k\leq d_n, n\geq 1\}$ and we will slightly abuse the notation to denote $\tau_k = \tau_{nk}$.
Homogeneous effects (i.e. when $\tau_k = \psi$ for all $k \in [d]$) can be helpful in terms of estimation. Intuitively, since the treatment effects are the same within each level of covariate, there is no need to consistently estimate CATE $\tau_k$ for all possible $d$ categories. One only needs to estimate the CATE within a few categories accurately and due to effect homogeneity, these estimated CATEs generalize to other levels of covariate, which potentially reduces the sample size required for consistent estimation. In this work, we adopt a novel form of approximate effect homogeneity via parameter capturing the extent of heterogeneity. Mathematically, denote
\begin{equation}\label{eq:homogeneity}
    \sigma_n := \max_{1 \leq k \leq d_n} |\tau_k -\psi|
\end{equation}
as the maximal effect heterogeneity. Note that $\sigma_n=0$ corresponds to the constant conditional average treatment effects and $\sigma_n=2$ imposes no restriction on the model class. Our parameterization can interpolate between these two extremes as $\sigma_n$ varies. 

\subsection{Upper Bound}

The estimator we propose to take advantage of the effect homogeneity is 
\begin{equation}\label{eq:homo-est}
    \hat{\tau} =  \frac{\sum_{k=1}^d \hat{t}_k \hat{\tau}_k}{\sum_{k=1}^d \hat{t}_k}
\end{equation}
where $\hat{t}_k = \hat{p}_kI( 0 < \hat{\pi}_k <1) $ is the indicator of whether both treated and untreated samples are observed in the $k$-th level and $\hat{\tau}_k = \hat{\mu}_{1k}-\hat{\mu}_{0k}$. The idea is to only estimate the CATE $\tau_k$ within those categories with both treated and untreated units and take a weighted average over such categories.  We restrict our attention to these categories so that unbiased estimation of $(\mu_{1k}, \mu_{0k})$ (and hence $\tau_k$) is possible. We again define $0/0=0$ if $\sum_{k=1}^d \hat{t}_k = 0$, i.e. there is no such category that contains both treated and untreated units. Clearly on the event $\left\{\sum_{k=1}^d \hat{t}_k = 0\right\}$, we cannot obtain any information on treatment effects from $\hat{\tau}$. Under appropriate scaling conditions, the probability of this ``adverse" event converges to $0$ even in the case $n \lesssim d$, as summarized in the following lemma.

\begin{lemma}\label{lemma:nocollision}
The chance of having every category consist of either all treated or all untreated units (i.e., no collisions of any treated and untreated units at any category) is upper-bounded as 
\[
\MP\left(\sum_{k=1}^d \hat{t}_k = 0\right) \leq 2 \exp\left(-C (\epsilon) \frac{n^2}{n \vee d} \right)
\]
where $C(\epsilon)> 0$ is a constant depending on $\epsilon$.
\end{lemma}
This lemma is not only critical to our analysis of $\hat{\tau}$ but also of independent interests in the literature of applied probability. It is closely related to the birthday problem \citep{clevenson1991majorization} and can be viewed as the occupancy problem \citep{wendl2003collision, nakata2014number} under a different sampling scheme. In the classic occupancy problem, the number of treated and untreated units are fixed first and then they are assigned to different categories of $X$ with probability vectors $\bp^{\text{T}}$ and $\bp^{\text{C}}$ separately. While in our setting people's covariates are first sampled, following which their treatment assignments are determined. Importantly, Lemma \ref{lemma:nocollision} shows the probability that the denominator of $\hat{\tau}$ is $0$ vanishes even in the high-dimensional regime $n \lesssim d$, as long as $d/n^2 \rightarrow 0$. With Lemma \ref{lemma:nocollision} we can derive the following theorem bounding the MSE of $\hat{\tau}$.

\begin{theorem}\label{thm:homogeneity}
    The bias and variance of $\hat{\tau}$ are bounded as 
    \begin{equation*}
        |\ME[\hat{\tau}-\psi]| \leq \sigma_n + 2 \exp\left( 
-C(\epsilon)\frac{n^2}{n \vee d} \right),
    \end{equation*}
    \[
    \operatorname{Var}(\hat{\tau}) \lesssim  \sigma_n^2  + \frac{ d}{n^2} +\frac{1}{n},
    \]
    where $C(\epsilon)$ and the constant behind ``$\lesssim$" depend on $\epsilon$.
    As a consequence, the estimator $\hat{\tau}$ is consistent if $\sigma_n \rightarrow 0$ and $d/n^2 \rightarrow 0$ in the regime $n \lesssim d$. 
\end{theorem}


Theorem \ref{thm:homogeneity} shows that if asymptotic effect homogeneity holds (i.e. $\sigma_n \rightarrow 0$), then the estimator $\hat{\tau}$ is consistent in the regime $d/n^2 \rightarrow 0$. Compared with the rate condition $d/n \rightarrow0$ required for consistency by the plug-in estimator $\hat{\psi}$, $\hat{\tau}$ has a faster convergence rate and enables us to achieve consistency in a wider regime under asymptotic effect homogeneity.


\subsection{Minimax Lower Bound}
In this section, we establish a matching minimax lower bound under exact effect homogeneity. Specifically, we consider the model class
\[
\mathcal{H}(\epsilon):=\left\{p_k=\frac{1}{d},\ \epsilon \le \pi_k \le 1-\epsilon,\ 0 \le \mu_{0k},\mu_{1k}\le 1,\ \tau_k=\tau_{k'},\, \forall  k,k'\in[d]\right\}.
\]
By Theorem~\ref{thm:homogeneity}, the estimator $\hat{\tau}$ attains the upper bound
\[
\frac{d}{n^2}+\frac{1}{n}.
\]
The following theorem shows that this rate is minimax optimal over $\mathcal{H}(\epsilon)$.

\begin{theorem}\label{thm:homo-minimax}
    In the regime $d \lesssim n^2$, we have
    \[
    \inf_{\hat{\psi}} \sup_{\MP \in \mathcal{H}(\epsilon)}\ME_{\MP}\left[\left(\hat{\psi} - \psi\right)^2\right] \gtrsim  \frac{d}{n^2} + \frac{1}{n},
    \]
    where the constant behind ``$\gtrsim$" depends on $\epsilon$.
\end{theorem}

Theorem~\ref{thm:homo-minimax} implies that consistency over the model class $\mathcal{H}(\epsilon)$ cannot be achieved when $d \gtrsim n^2$. Intuitively, in this regime, with high probability, very few categories contain repeated observations, and in particular very few contain both treated and untreated units. Moreover, one can construct pairs of distributions whose one-observation marginals within each category are identical, so that categories observed only once carry no information for distinguishing them. As a result, only within-category collisions are informative, and when $d \gtrsim n^2$ there are too few such collisions to permit consistent estimation.

Theorem \ref{thm:homogeneity} establishes that $\sigma_n \to 0$ is a sufficient condition for consistency when $d/n^2 \to 0$. We regard $\sigma_n \rightarrow 0$ as an interpretable structural assumption capturing asymptotic effect homogeneity. In addition, the lower bound in Theorem \ref{thm:ate-minimax} suggests that some such restriction is necessary to extend consistency beyond the general $d=o(n\log n)$ regime. We leave a sharp characterization of the minimax dependence on $\sigma_n$ to future work.

\section{Prior Knowledge of Covariate Distribution}\label{sec:cova-distribution}

In this section, we consider a different structure that could be exploited in causal effects estimation: the covariate distribution is known or can be estimated at a fast rate. The role of the covariate distribution has been studied in both causal functional estimation (e.g., expected conditional covariance (ECC); see \citet{robins2008higher}) and function estimation (e.g., CATE; see \citet{kennedy2023towards}), where incorporating knowledge of the covariate density can lead to improved convergence rates. When the covariate is discrete, faster rate of ATE estimation is not achievable with information on covariate distribution, as shown in Section \ref{sec:cova-ate}. The fundamental difficulty is that the measure underlying ATE estimation is the product of covariate distribution and propensity score, whereas propensity score is unknown in our observational study setting. We thus switch attention to variance-weighted average treatment effects (WATE) in Section \ref{sec:cova-wate}, which itself also receives an increasing amount of attention in the literature. The underlying measure of WATE is the covariate distribution and a faster rate of estimation is achieved with covariate information.

\subsection{Average Treatment Effects}\label{sec:cova-ate}

In this section, we introduce a minimax lower bound for $\psi_1= \ME[\ME(Y|X,A=1)]$ under a uniform covariate distribution $\bp = (1/d,\cdots,1/d)$. We show that a significant improvement on the rate of ATE estimation, such as the $d/n^2$ rate discussed in Section \ref{sec:homogeneity}, cannot be attained. This highlights the limitations of leveraging covariate distribution knowledge in improving the efficiency of ATE estimation. 
Consider the following model class where the covariate distribution is uniform over $[d]$:
\[
\mathcal{D}^U(\epsilon) = \left\{ p_k = 1/d  , \epsilon \leq \pi_k \leq 1-\epsilon, 0\leq \mu_{0k}, \mu_{1k} \leq 1 ,\, \forall 1\leq k \leq d  \right \}.
\]
The minimax lower bound of $\psi_1$ over model class $\mathcal{D}^U(\epsilon)$ is summarized in the following theorem.
\begin{theorem}\label{thm:ate-minimax-cova}
    For any fixed constant $\beta \in (0,1)$, in the regime $n \lesssim d^{1-\beta}$ we have
    \[
    \inf_{\hat{\psi}_1} \sup_{\MP \in \mathcal{D}^U (\epsilon)}  \ME_{\MP}\left[\left(\hat{\psi}_1 - \psi_1\right)^2\right] \gtrsim 1,
    \]
    where the constant behind ``$\gtrsim$" depends on $\beta, \epsilon$. 
\end{theorem}

Similar to the proof of Theorem \ref{thm:ate-minimax}, the proof of Theorem \ref{thm:ate-minimax-cova} relies on the method of fuzzy hypotheses \citep{Tsybakov2009} and the characterization of the following best 2D polynomial approximation error: 
\begin{align}\label{eq:2D-approx-error}
\inf_{Q\in \mathbb{R}[x,y], \deg Q\le L} \max_{\pi\in [\epsilon, 1-\epsilon], \mu\in [0,1]} |Q(\pi, \pi\mu) - \mu|. 
\end{align}
Again, we lower bound the above quantity via a suitable $1$-dimensional subproblem, which simply sets $\mu = \epsilon/\pi$. The main difference from the proof of Theorem \ref{thm:ate-minimax} is that, instead of choosing $L \asymp \log n$, here we choose $L = O(1)$ to be a large constant. It turns out that the approximation error \eqref{eq:2D-approx-error} is lower bounded by a constant $c=c(\epsilon,L)>0$, and 
the two fuzzy hypotheses are statistically indistinguishable as long as $n\lesssim d^{1-\beta}$.

Theorem \ref{thm:ate-minimax-cova} establishes that the minimax rate for estimating $\psi_1$ is lower bounded away from zero within the regime $n\lesssim d^{1-\beta}$, provided the model class includes $\mathcal{D}^U(\epsilon)$ as a subclass. This result indicates that enhancing the estimation rate of $\psi_1$ by incorporating knowledge of the covariate distribution may not yield significant improvements when $n \lesssim d$. Specifically, achieving consistency in the scenario where $d/n^2 \rightarrow 0$—analogous to the situation of effect homogeneity discussed in Section \ref{sec:homogeneity}—is unattainable solely with insights into the covariate distribution.

The negative result in Theorem \ref{thm:ate-minimax-cova} can be explained by the fact that the underlying measure of ATE estimation is the product of covariate distribution and propensity score, as shown in analyzing its second-order estimator \citep{robins2009quadratic,zeng2023efficient}. In observational studies where the propensity score is unknown, it's hard to obtain nearly unbiased estimates of $p_k \mu_{1k}$ (the summand in $\psi_1$) with the observed data. One may write $p_k \mu_{1k} = p_k \pi_k \mu_{1k}/\pi_k$ with $p_k \pi_k \mu_{1k}$ unbiasedly estimable. The extra term $1/\pi_k$ could be approximated by a polynomial with possibly diverging degrees
\[
\frac{1}{\pi_k} \approx \sum_{j=0}^L (1-\pi_k)^j,
\]
where $L= L_n$ determines the approximation error. Since the product $(p_k \pi_k)^j$ can be estimated unbiasedly and $p_k$ is known, this strategy could reduce the bias. However, this approach suffers from a large variance, particularly in situations where $n \lesssim d$. In fact, our construction in the proof of Theorem \ref{thm:ate-minimax-cova} indeed replies on the difficulty of approximating $1/\pi_k$ with a polynomial of $\pi_k$, thereby highlighting the fundamental barrier in ATE estimation even with information on the covariate distribution.

\subsection{Variance-weighted Average Treatment Effects}\label{sec:cova-wate}

In this section we switch our attention to a popular variance-weighted average treatment effect estimand (WATE), defined as  
\[
\theta= \frac{\ME[\operatorname{Cov}(Y,A\mid X)]}{\ME[\operatorname{Var}(A\mid X)]} = \frac{\ME[\operatorname{Var}(A\mid X )\tau_X]}{\ME[\operatorname{Var}(A\mid X)]},
\]
where the weight $\operatorname{Var}(A\mid X)$ minimizes the asymptotic variance among all weighted average treatment effects with known weight under homoskedasticity \citep{crump2006moving}. Another interpretation of $\theta$ in terms of robustness is that if we assume a partial linear (homogeneous effect) model for the outcome and use E-estimator \citep{robins1992estimating} to estimate the parametric component, under model misspecification (i.e. the partial linear model is wrong and effect of treatment is heterogeneous), one can show the E-estimator converges to WATE \citep{vansteelandt2014regression}. The WATE has also been derived under different frameworks recently. \cite{zhou2022marginal} showed the marginal interventional effect under incremental propensity score intervention \citep{kennedy2019nonparametric} coincides with WATE. The expected conditional covariance also appears in conditional independence test \citep{shah2020hardness} and can be interpreted as a causal effect under a stochastic intervention \citep{diaz2023non}. In contrast to ATE, the underlying measure of $\theta$ is the covariate distribution alone, making improvements from knowledge of covariate distribution possible.


The strategy of estimating $\theta$ is to deal with the numerator and the denominator separately. Let $\eta = \ME[\operatorname{Cov}(Y,A\mid X)]$ and $\rho = \ME[\operatorname{Var}(A\mid X)]$. Further denote the regression function in $k$-th category as $\mu_k : = \ME[Y|X=k]$. The first-order and second-order influence functions of $\eta$ under a nonparametric model, if they exist, have form
\[
\varphi_1(\bZ) = (A-\pi_X)(Y-\mu_X),
\]
\[
\varphi_2(\bZ_1,\bZ_2) = -\frac{(A_1-\pi_{X_1})I(X_1=X_2)(Y_2-\mu_{X_2})}{p_{X_1}}.
\]
Similar forms of influence functions of $\rho$ can be obtained by replacing $(Y,\mu)$ with $(A,\pi)$. However, in general settings where covariates have continuous components, the second-order influence function does not exist for $\eta$ \citep{robins2009quadratic}. Intuitively, for two independent samples $X_1, X_2$ from a continuous distribution, we have $X_1 \neq X_2$ with probability $1$ and hence $\varphi_2(\bZ_1,\bZ_2)=0$. So the ``second-order influence function" is always $0$ when covariates have continuous components and cannot help us improve the estimator. On the other hand, when $X$ is discrete it is possible that $X_1, X_2$ fall into the same category and $\varphi_2(\bZ_1,\bZ_2) \neq 0$. Then the idea is to use second-order estimator based on $\varphi_1, \varphi_2$ to simultaneously correct for the first-order and second-order bias of plug-in-style estimator $\hat{\eta} = \MP_n[YA - \hat{\mu}_X \hat{\pi}_X]$. Consider the following second-order estimator of $\eta$ and $\rho$,
\begin{equation}\label{eq:high-order-eta}
    \hat{\eta} = \MP_n[(A-\hat{\pi}_X)(Y-\hat{\mu}_X)] - \MU_n \left[ \frac{(A_1-\hat{\pi}_{X_1})I(X_1=X_2)(Y_2-\hat{\mu}_{X_2})}{\hat{p}_{X_1}}\right],
\end{equation}
\begin{equation}\label{eq:high-order-rho}
    \hat{\rho} = \MP_n\left[(A-\hat{\pi}_X)^2\right] - \MU_n \left[ \frac{(A_1-\hat{\pi}_{X_1})I(X_1=X_2)(A_2-\hat{\pi}_{X_2})}{\hat{p}_{X_1}}\right],
\end{equation}
where $\MP_n,\MU_n$ are empirical and U-statistic measures. The following theorem establishes the estimation guarantees of second-order estimators when the nuisance functions are estimated from a separate independent sample agnostically, i.e., we do not specify the way to estimate $\mu$ and $\pi$.

\begin{theorem}\label{thm:cova-dist-agnostic}
Suppose the nuisance functions $\pi_k,\mu_k$ and covariate distribution $p_k$ are estimated from a separate independent sample $D$ as $\hat{\pi}_k,\hat{\mu}_k,  \hat{p}_k$ with $\hat{\mu}_k, \hat{\pi}_k \in [0,1]$. Then for second-order estimators in \eqref{eq:high-order-eta}--\eqref{eq:high-order-rho} we have
\[
\begin{aligned}
    |\ME\left[\hat{\eta}\mid D\right]-\eta| \leq&\,  \|\hat{\mu}-\mu\|_2 \|\hat{\pi}-\pi\|_2 \max_{k} \left| 1-\frac{p_k}{\hat{p}_k} \right|, \\
    \operatorname{Var}\left(\hat{\eta}\mid D\right) \lesssim &\, \frac{1}{n} \left(1+ \max_k \left| 1-\frac{p_k}{\hat{p}_k} \right|\right)^2 + \frac{d}{n^2}\left(1+ \max_k \left| 1-\frac{p_k}{\hat{p}_k} \right|\right)^2, \\
    |\ME\left[\hat{\rho}\mid D \right]-\rho| \leq&\,   \|\hat{\pi}-\pi\|_2^2 \max_{k} \left| 1-\frac{p_k}{\hat{p}_k} \right|, \\
    \operatorname{Var}\left(\hat{\rho}\mid D\right) \lesssim &\, \frac{1}{n} \left(1+ \max_k \left| 1-\frac{p_k}{\hat{p}_k} \right|\right)^2 + \frac{d}{n^2}\left(1+ \max_k \left| 1-\frac{p_k}{\hat{p}_k} \right|\right)^2,
\end{aligned}
\]
where recall for a function $f = f(\bZ)$ the $L_2$-norm is defined as $\|f\|_2^2 = \int f^2(\bz) d \MP (\bz)$ and the constant behind ``$\lesssim$" is an absolute constant.
\end{theorem}

Theorem \ref{thm:cova-dist-agnostic} shows the advantage of second-order estimators: after correcting for the first-order and second-order bias, the conditional bias of $\hat{\eta}$ depends on the product of estimation error of $\mu, \pi, p$ and is ``third-order small". In the following theorem, we parameterize the estimation rate of covariate distribution to show the consistency of second-order estimators in the high-dimensional regime $n \lesssim d$ with different choices of nuisance estimators $\hat{\pi}, \hat{\mu}$.

\begin{theorem}\label{thm:cova-dist}
    Suppose the nuisance functions $\pi_k,\mu_k$ and covariate distribution $p_k$ are estimated from a separate independent sample $D$ as $\hat{\pi}_k,\hat{\mu}_k,  \hat{p}_k$. Further assume the estimated covariate distribution $\hat{p}_k$ satisfies
    \[
    \max_{1 \leq k \leq  d} \left | 1-\frac{p_k}{\hat{p}_k} \right|\leq \xi_n < 1. 
    \]
    If nuisance functions $\pi_k, \mu_k$ are estimated from a sample of size $n$ using empirical averages, then the second-order estimators of $\eta$ and $\rho$ in \eqref{eq:high-order-eta}--\eqref{eq:high-order-rho} satisfy
    \[
    \ME\left[ (\hat{\eta}-\eta)^2 \right] \lesssim \xi_n^2 \frac{d\wedge n}{n} + \frac{d }{n^2}+ \frac{1}{n},
    \]
    \[
    \ME\left[ (\hat{\rho}-\rho)^2 \right] \lesssim \xi_n^2 \frac{d \wedge n}{n} + \frac{d}{n^2}+ \frac{1}{n}.
    \]
    If we set the nuisance functions $\hat{\pi}_k, \hat{\mu}_k$ as $0$, then the second-order estimators of $\eta$ and $\rho$ satisfy
    \[
    \ME\left[ (\hat{\eta}-\eta)^2 \right] \lesssim \xi_n^2  + \frac{d }{n^2}+ \frac{1}{n},
    \]
    \[
    \ME\left[ (\hat{\rho}-\rho)^2 \right] \lesssim \xi_n^2  + \frac{d}{n^2}+ \frac{1}{n}.
    \]
    where the constant behind ``$\lesssim$" is an absolute constant. 
    
    As a consequence, under positivity $\epsilon \leq \pi_k \leq 1-\epsilon$  and assume $\hat{\rho} \geq \epsilon(1-\epsilon)$, the estimator $\hat{\theta} = \hat{\eta}/\hat{\rho}$
    satisfies
    \[
    \begin{aligned}
        \ME\left[ (\hat{\theta}-\theta)^2 \right] \lesssim &\, \ME\left[ (\hat{\eta}-\eta)^2 \right] + \ME\left[ (\hat{\rho}-\rho)^2 \right]\\
        \lesssim &\, \xi_n^2  + \frac{d}{n^2}+ \frac{1}{n},
    \end{aligned}
    \]
    where the constant behind ``$\lesssim$" depends on $\epsilon$.
\end{theorem}

We note that Theorem \ref{thm:cova-dist} implies a faster convergence rate than $d^2/n^2 + 1/n$ when $\xi_n$ is small, which is hard to achieve for estimators $\hat{p}_k$'s constructed from a sample of size $n$ in the high-dimensional regime $n \lesssim d$. For example, in the uniform case $p_k = 1/d, 1\leq k \leq d$, one can use Chernoff bound to show that with probability $1-\delta$, 
\[
\max_{1\leq k \leq d}|d\hat{p}_k-1| \leq \sqrt{\frac{3d \log(2d/\delta)}{n}}.
\]
Thus we need a sample size of order $d \log d$ to guarantee small uniform estimation errors of $p_k$, which is not achievable in the case $n \lesssim d$. The improvements on the convergence rate of $\hat{\eta}$ should come from prior knowledge/assumptions on covariate distribution. The results imply that consistency is still possible in the regime $n \lesssim d$ when such information on covariate distribution is available. One example is that the covariate distribution is known to statisticians, then we have $\hat{p}_k = p_k,\xi_n = 0$ and the estimators are consistent as long as $d/n^2 \rightarrow 0$. Another useful setting is semi-supervised causal inference \citep{chakrabortty2022general, zhang2023semi, kallus2020role,zeng2024continuous}: there is a large number ($\gg n$) of individuals in the database who were not selected into the randomized trial/observational study, hence treatment and outcome are not available for them but one can use the database to estimate the covariate distribution very accurately thanks to the large sample size. Then $\xi_n$ is small and the estimation error of $\hat{\eta}$ will be small according to Theorem \ref{thm:cova-dist}.

When the covariate distribution can be estimated very well and $\xi_n$ is small, we can use arbitrary (bounded) nuisance estimators and still maintain consistency. As shown in Theorem \ref{thm:cova-dist}, simply setting $\hat{\pi}, \hat{\mu}$ as $0$ leads to consistent estimation as long as $\xi_n$ vanishes and $d/n^2 \rightarrow 0$. Using empirical average estimators has advantage over setting $\hat{\pi}, \hat{\mu}$ as $0$ only when $d \ll n $. The interpretation is that when $\xi_n$ is small, the conditional bias in Theorem \ref{thm:cova-dist-agnostic} is negligible regardless of the convergence rate of $\hat{\pi}, \hat{\mu}$. The conditional variance, whose order is independent of the nuisance estimation rate, is the dominant term in the MSE. 

We conclude this section with additional discussion on second-order estimation of $\psi_1 = \ME[\ME(Y\mid X,A=1)]$ in the discrete covariate setting. The forms of first-order and second-order influence functions are
\[
\phi_1(\bZ) = \frac{A(Y-\mu_{1X})}{\pi_X}+ \mu_{1X}, 
\]
\[
\phi_2(\bZ_1,\bZ_2) = -\left(\frac{A_1}{\pi_{X_1}}-1 \right)\frac{I(X_1=X_2)}{p_{X_1}\pi_{X_2}}A_2(Y_2-\mu_{1X_2}).
\]
The second-order estimator is then
\begin{equation}\label{eq:high-order-ate}
\MP_n\left[\hat{\phi}_1(\bZ)\right] + \MU_n\left[\hat{\phi}_2(\bZ_1,\bZ_2)\right].
\end{equation}
Assume we use a separate independent sample $D$ to estimate nuisance functions $\pi, \mu_1$ and covariate distribution $p$, the conditional bias of second-order estimator \eqref{eq:high-order-ate} is
\[
\ME \left[ (\hat{\mu}_{1X}-\mu_{1X})\left(1-\frac{\pi_X}{\hat{\pi}_X} \right)\left(1-\frac{\pi_X p_X}{\hat{\pi}_X \hat{p}_X} \right) \mid D\right].
\]
In order to obtain a faster rate by making use of information on covariate distribution, one also needs an accurate estimator of propensity score $\pi$ to ensure
\begin{equation}\label{eq:product-pi-p}
    \max_{1 \leq k \leq d}\bigg |1-\frac{\pi_k p_k}{\hat{\pi}_k \hat{p}_k} \bigg |
\end{equation}
is small. In observational studies where the treatment process is unknown, it is hard to estimate the product of propensity score and covariate distribution well and make \eqref{eq:product-pi-p} small in the regime $n \lesssim d$, due to the same reason discussed after Theorem \ref{thm:cova-dist}. This further demonstrates the difficulty in ATE estimation with solely information on covariate distribution. 

The fact that covariate distribution improves the estimation rates for ECC and WATE, but not for ATE, can be understood through the framework of \citet{robins2008higher}. In their analysis, the improved convergence rates achieved by higher-order estimators rely on knowledge of the so-called “$g$-function” (defined in their equation (3.3)), which is assumed to be known or sufficiently smooth. For the expected conditional covariance (ECC) $\eta$, this $g$-function simplifies to the covariate density. Hence, knowledge of the covariate distribution directly contributes to improved estimation rates. In contrast, for the covariate-adjusted mean in the treated group $\psi_1$, the $g$-function involves the product of the covariate density and the propensity score. Therefore, faster rates can only be achieved when information is available about this product, not the covariate distribution alone. In practice, WATE is also a useful and interpretable summary measure of treatment effects that is often important to estimate and report.

\section{Simulation Study}\label{sec:simulation}
In this section, we provide numerical results to verify the theoretical properties of the estimators discussed in Section \ref{sec:ate-plugin}--\ref{sec:cova-distribution}. Consider the following data generating process: $X$ has a categorical distribution with uniform probability $\bp=(1/d,\dots, 1/d)^{\top} \in \mathbb{R}^d$. For each category $X=k$, the propensity score $\pi_k = 1/2$ and the conditional means of outcome within treated and untreated groups are $\mu_{1k}=1/2, \mu_{0k}=1/4$, respectively. Hence the treatment effects are homogeneous across different levels (i.e. $\tau_k = 1/4$ for all $k \in [d]$) in our simulation setting and expected conditional covariance between treatment and outcome $\eta=1/16$ (considered in Section \ref{sec:cova-distribution}). For each estimator we generate data containing $n$ observations with $n \in \{100,500,1000,5000,10000 \}$, apply the proposed method, and compute the estimated RMSE. Such procedure is repeated $M=500$ times and the average RMSE is reported in each plot.

\begin{figure}[H]
    \centering
    \includegraphics[width=5in]{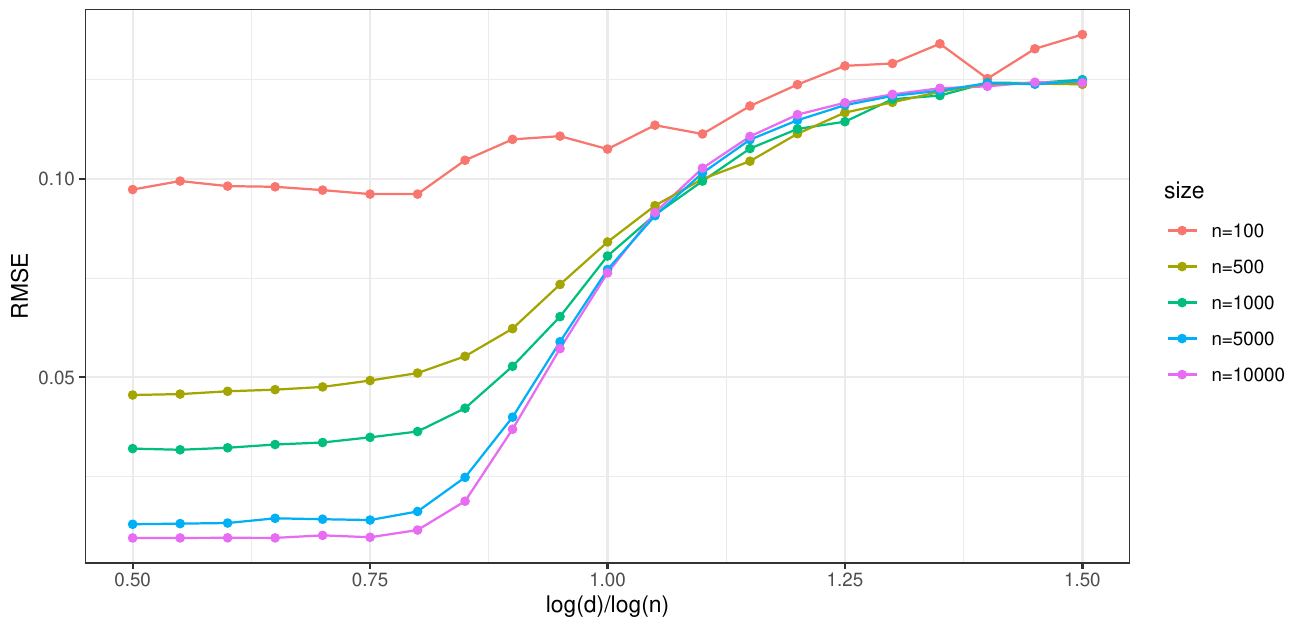}
    \caption{RMSE V.S. $\gamma$ for the plug-in estimator, where $\gamma$ controls the order of $d$.}
    \label{fig:rmse-plugin}
\end{figure}

We first consider the plug-in estimator discussed in Section \ref{sec:ate-plugin}. For each sample size $n$ let the number of categories $d=\lfloor n^{\gamma} \rfloor$ with $\gamma \in \{0.5, 0.55, 0.6,\\ \dots, 1.5\}$. The idea is to evaluate how large the order of $d$ can be to maintain consistency \citep{he2021phase}. The average treatment effect is $\psi=1/4$ and we estimate RMSE as
\[
\widehat{\text{RMSE}} =\left[ \frac{1}{M}\sum_{m=1}^M (\hat{\psi}^m-\psi)^2 \right]^{1/2},
\]
where $\hat{\psi}^m$ is the estimated ATE from $m$-th repetition. The relationship between RMSE and $\gamma$ is summarized in Figure \ref{fig:rmse-plugin}. For a fixed sample size $n$, as $\gamma$ increases (i.e. $d$ increases) the estimated RMSE also increases as expected. We see a clear phase transition in the plot: in the region $\gamma < 0.8$, the RMSE is quite stable of order $1/\sqrt{n}$; around $\gamma=1$ the RMSE increases drastically, indicating the plug-in estimator starts to have larger error. This corresponds to our theoretical results in Section \ref{sec:ate-plugin}: the plug-in estimator is consistent if and only if $d/n \rightarrow 0$.
Then we consider the estimator \eqref{eq:homo-est} proposed in Section \ref{sec:homogeneity} under effect homogeneity. Again for each fixed sample size $n$ we let $d=\lfloor n^{\gamma} \rfloor$ with $\gamma \in \{0.5, 0.55, 0.6, \dots, 2\}$. Note that the treatment effects are homogeneous in our setting (i.e. $\sigma_n = 0$) and we expect $\hat{\tau}$ in equation \eqref{eq:homo-est} to be consistent in a wider regime. 

\begin{figure}[H]
    \centering
    \includegraphics[width=5in]{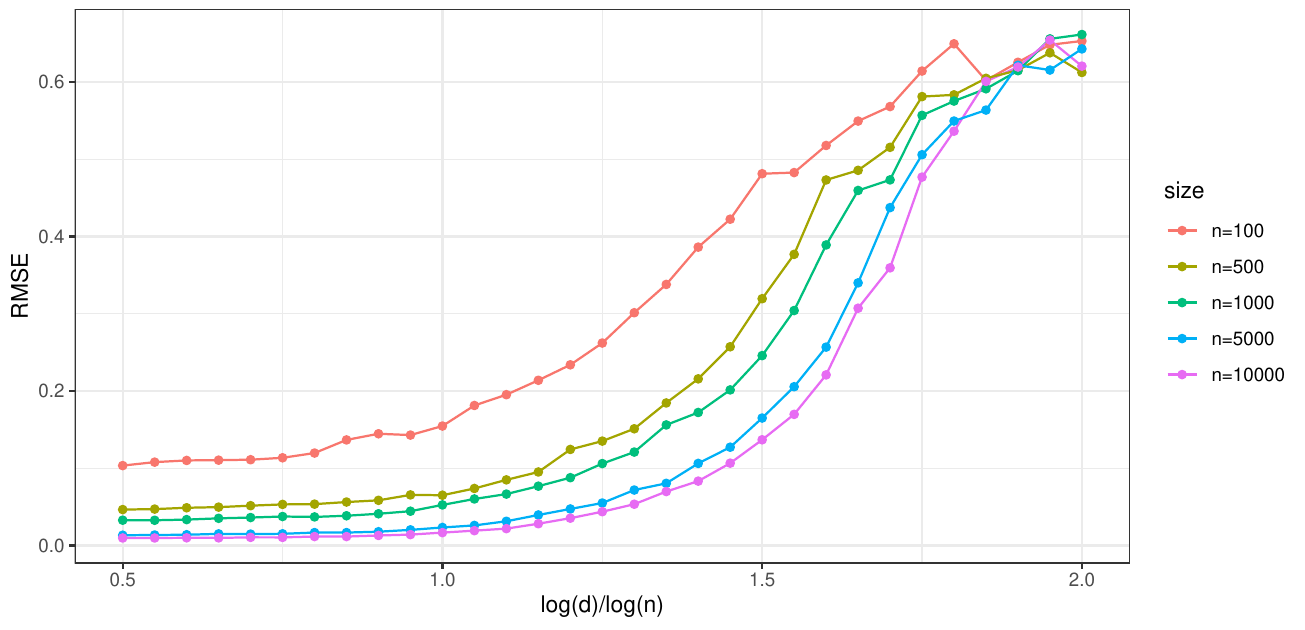}
    \caption{RMSE V.S. $\gamma$ for the estimator under effect homogeneity, where $\gamma$ controls the order of $d$.}
    \label{fig:rmse-homo}
\end{figure}

The relationship between RMSE and $\gamma$ is summarized in Figure \ref{fig:rmse-homo}. The phase transition occurs in the region $\gamma > 1$ (instead of at $\gamma=1$) and the estimator in \eqref{eq:homo-est} has a smaller error in the region $1 < \gamma < 1.25$ compared with the plug-in estimator when $n \geq 1000$. For instance, in the case $n=1000, \gamma = 1$ the plug-in estimator has RMSE around $0.08$ while the estimator under effect homogeneity has RMSE around $0.05$. This coincides with our theoretical results in Section \ref{sec:homogeneity} that the estimator in equation \eqref{eq:homo-est} has a faster rate and is consistent in a wider regime. 

\begin{figure}[H]
	\centering
	\subfigure[n=1000]{
		\begin{minipage}[t]{0.48\linewidth}
			\centering
			\includegraphics[width=2.5in]{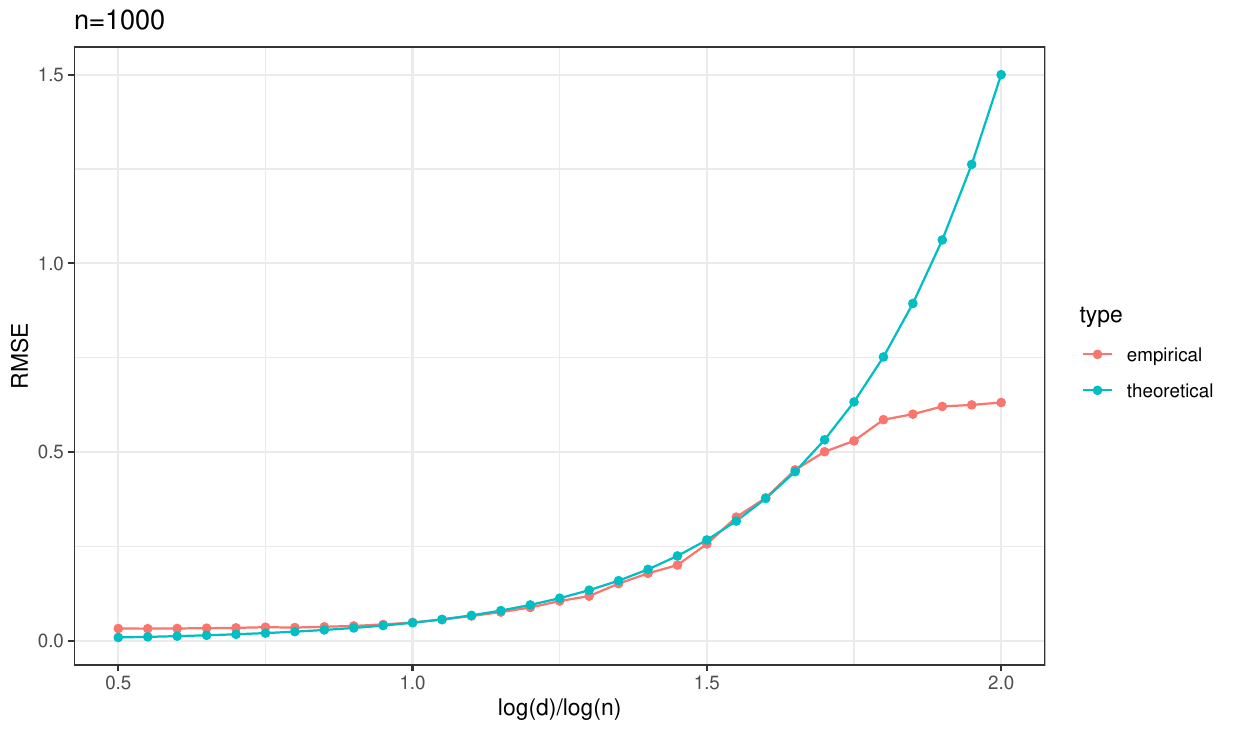}
	\end{minipage}}
	\subfigure[n=10000]{
		\begin{minipage}[t]{0.48\linewidth}
			\centering
			\includegraphics[width=2.5in]{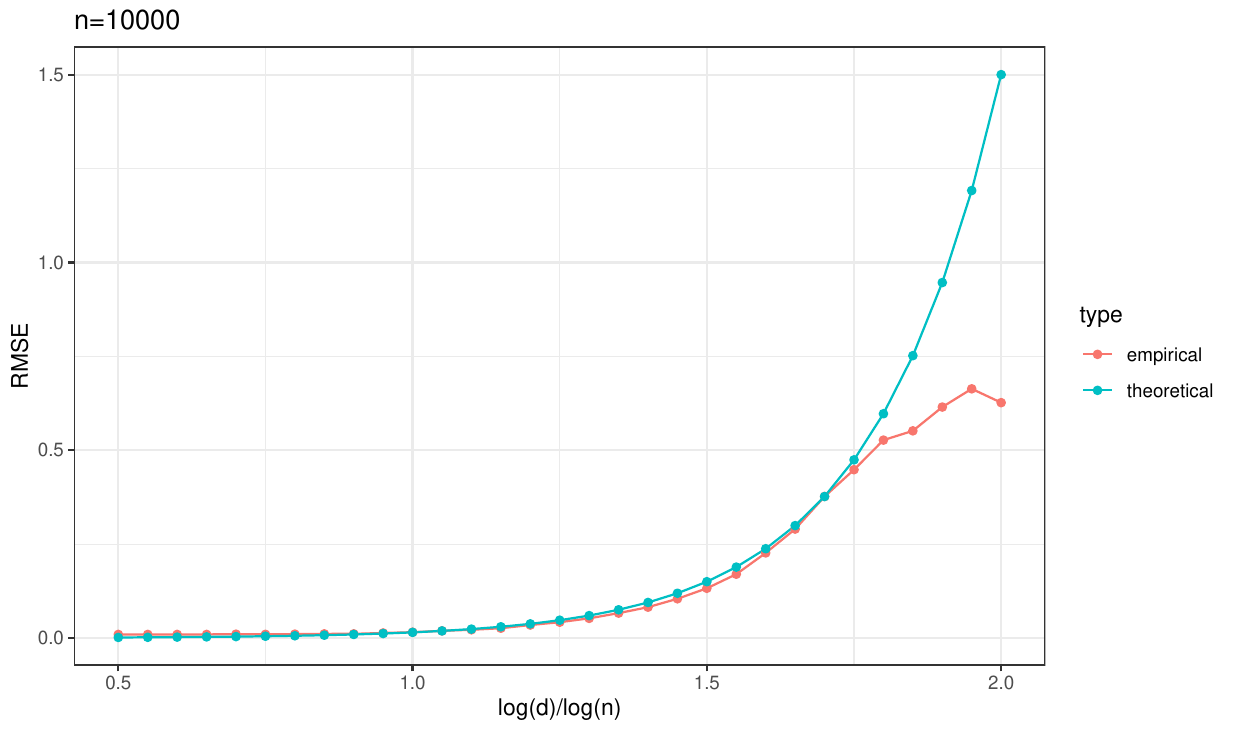}
	\end{minipage}}
	\centering
	\caption{Comparison of theoretical and empirical order of RMSE }
	\label{fig:compare-homo}
\end{figure}

To further understand the order of the RMSE, for $n \in \{1000, 10000\}$ we include the theoretical upper bound on RMSE, which is of order $C_1\sqrt{d/n^2} = C_1 n^{\gamma/2-1} $, in the plot as a benchmark. Here the constant $C$ is chosen as $1.5$ to fit the empirical RMSE curve. The results are summarized in Figure \ref{fig:compare-homo}. When $\gamma < 1.5$, the empirical RMSE fits the theoretical upper bound quite well. When $\gamma>1.5$ the empirical RMSE starts to deviate from the theoretical bound. In our experiments we found when $\gamma > 1.5$ and $d$ is large, the denominator in the estimator \eqref{eq:homo-est} is usually small (i.e., only a few categories have both treated and untreated samples) and the variation is large, which may explain the deviation of the empirical RMSE from the theoretical bound.

\begin{figure}[H]
    \centering
    \includegraphics[width=5in]{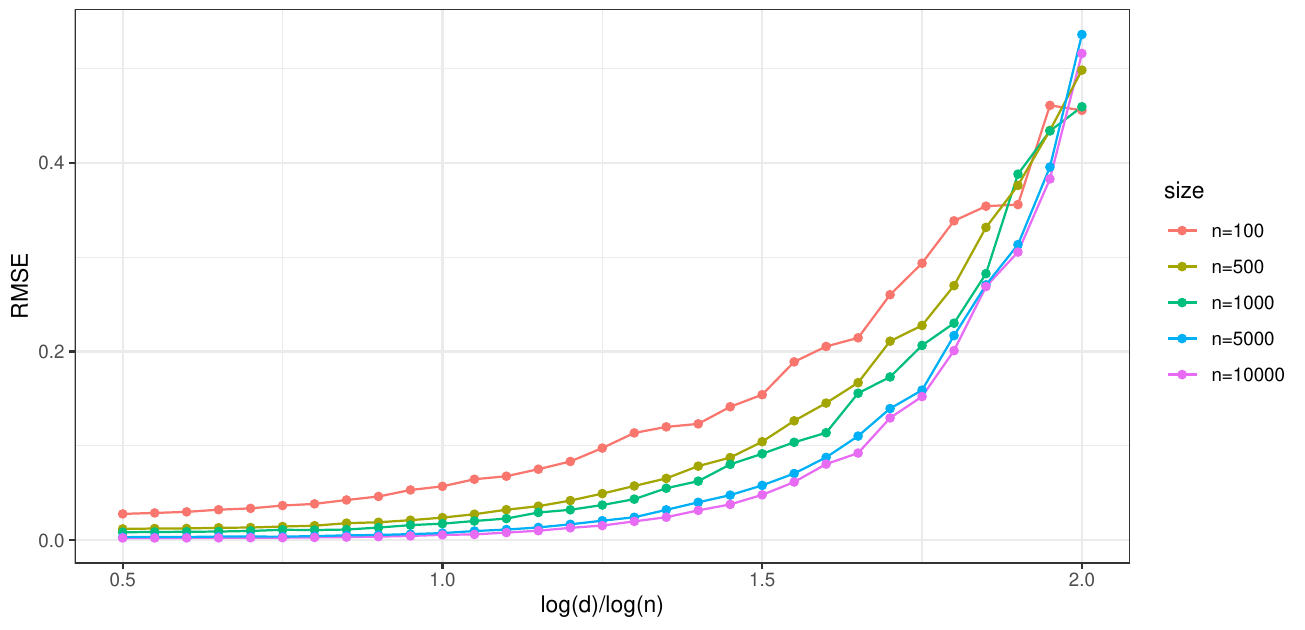}
    \caption{RMSE V.S. $\gamma$ for the second-order estimator \eqref{eq:high-order-eta} using true covariate distribution $\bp$, where $\gamma$ controls the order of $d$.}
    \label{fig:rmse-cova}
\end{figure}

\begin{figure}[H]
	\centering
	\subfigure[n=1000]{
		\begin{minipage}[t]{0.48\linewidth}
			\centering
			\includegraphics[width=2.5in]{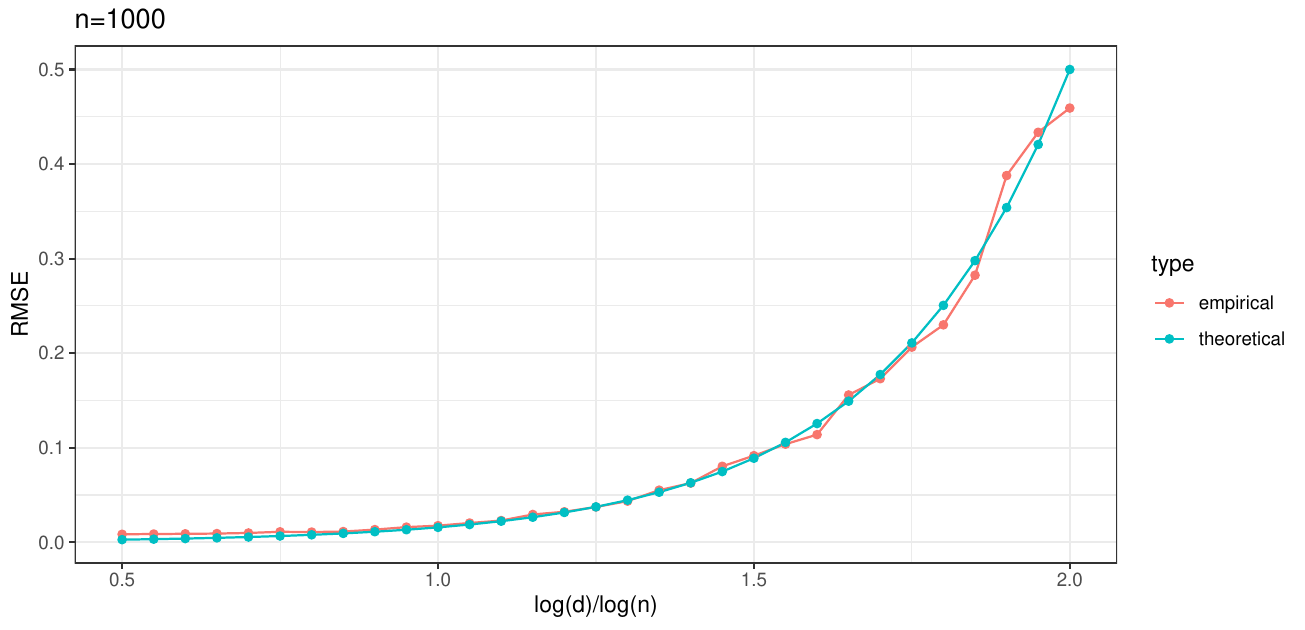}
	\end{minipage}}
	\subfigure[n=10000]{
		\begin{minipage}[t]{0.48\linewidth}
			\centering
			\includegraphics[width=2.5in]{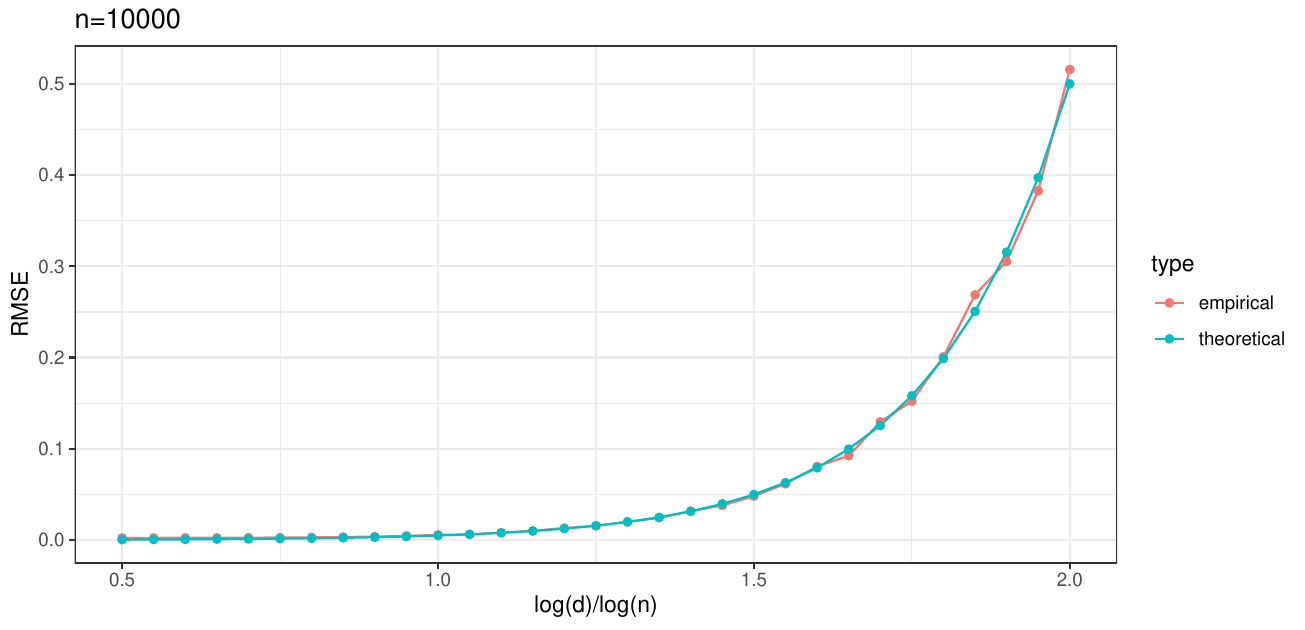}
	\end{minipage}}
	\centering
	\caption{Comparison of theoretical and empirical order of RMSE }
	\label{fig:compare-cova}
\end{figure}

Finally, we evaluate the performance of the second-order estimator in  \eqref{eq:high-order-eta} of expected conditional covariance $\eta$  (the estimator \eqref{eq:high-order-rho} of expected conditional variance of treatment is expected to have similar performance). In our setting $\eta=1/16$ and for each fixed sample size $n$, let $d=\lfloor n^{\gamma} \rfloor$ with $\gamma \in \{0.5, 0.55, 0.6, \dots, 2\}$. We set the estimates of probabilities $\hat{p}_k$'s as the true $p_k = 1/d$ and hence $\xi_n = 0$ in Theorem \ref{thm:cova-dist}. The estimates $\hat{\pi}_k, \hat{\mu}_k$'s are all set to $0$. The results are summarized in Figure \ref{fig:rmse-cova}. Similar to Figure \ref{fig:rmse-homo}, the estimated RMSE is quite stable in the region $\gamma < 1.25$ and the phase transition seems to happen in the region $\gamma > 1$.
We also plot the relationship between empirical RMSE and theoretical bound $C_2\sqrt{d/n^2} =C_2 n^{\gamma/2-1} $ with $C_2 = 0.5$ in Figure \ref{fig:compare-cova}. We see the empirical RMSE matches the theoretical bound very well.

\section{Real Data Analysis}
In this section, we employ the methods presented in Sections \ref{sec:ate-plugin} and \ref{sec:homogeneity}
to investigate the impact of 401(k) participation on net financial assets, utilizing data from the study by \citet{chernozhukov2004effects}. Since 401(k) eligibility is not randomly assigned, the effects must be assessed within an observational study framework. Following the strategy outlined in \citet{poterba1994401, poterba1995401, chernozhukov2018double}, 401(k) eligibility can be treated as exogenous after adjusting for confounding variables related to job choice that might be correlated with eligibility, such as income and education level.

The outcome variable $Y$ represents net financial assets, which includes IRA balances, 401(k) balances, checking accounts, U.S. savings bonds, other interest-earning accounts in banks and financial institutions, other interest-earning assets (such as personally held bonds), stocks, and mutual funds, less nonmortgage debt. The treatment variable $A$ is an indicator of 401(k) eligibility, where $A=1$ represents eligibility. The confounding variables $X$ include age, income, family size, years of education, marital status, two-earner status, defined benefit pension status, IRA participation, and home ownership. In our analysis, the variables age, income, family size, and years of education are discretized into categorical variables with four levels each. The sample size is $n=9915$, and the number of categories, considering a saturated model with all possible interactions, is $d=8196$.

The treatment effect estimated using the plug-in-style estimator discussed in Section \ref{sec:ate-plugin} is $7802$. The 95\% confidence intervals, based on bootstrap and asymptotic normality in the appendix, are $[5352, 9559]$ and $[6212, 10238]$, respectively. It is important to note that, given the sample size $n=9915$ and number of categories $d = 8196$, these confidence intervals may not achieve the nominal coverage probability  \citep{Karoui2018}; thus, they are presented for reference. Simulation results illustrating bootstrap-based variance estimation in high-dimensional settings with discrete covariates are provided in the appendix. These results indicate that 401(k) eligibility significantly increases the net financial assets. This finding aligns with the results of \citet{poterba1994401,poterba1995401,chernozhukov2018double}. If researchers are confident that the covariates account for all relevant confounders, our estimate can be interpreted causally; otherwise, it should be understood as the expected difference between the covariate-adjusted regression functions for the treatment and control groups. 

Additionally, we identified 993 categories in our data where positivity may be violated, meaning there is no overlap between treatment and control within those categories. The estimator in \eqref{eq:homo-est}, which only considers categories with both treated and control units, yields an estimate of $9149$, with a 95\% bootstrap confidence interval of  $[6164,11715]$. This confidence interval is comparable in length to that of the plug-in-style estimator. It is possible that the treatment effects are heterogeneous, suggesting that the faster convergence rate of the estimator in \eqref{eq:homo-est} under effect homogeneity may not be fully achievable with the given data.
\section{Discussion}

In this paper, we studied the treatment effects estimation problem in the context of high-dimensional discrete covariates. Theoretical properties of commonly used regression, weighting and doubly robust estimators are examined in this non-classic regime. We also evaluated the role of effect homogeneity and covariate distribution in treatment effects estimation and proposed estimators that can properly take advantage of these structures and achieve faster convergence rates. Finally, we explored the fundamental limits of treatment effects estimation and showed consistent estimation of ATE is a difficult task on high-dimensional data. The discrete covariate setting is not only an interesting base case but also informative for the general dataset with continuous components. We hope our work can help researchers appropriately understand and interpret the treatment effects estimated from datasets with many covariates.

There are several possible extensions for future work. In this paper, we borrowed the moment matching tools from the theoretical computer science literature \citep{jiao2015minimax, wu2016minimax, wu2019chebyshev} to study the minimax lower bound on ATE estimation. It would be interesting to explore how to construct estimators of ATE based on polynomial approximation theory and realize the “effective sample size enlargement" phenomenon, which is feasible in entropy and support size estimation. Moreover, the minimax lower bound for ATE we proved is an initial result that does not exactly match the upper bound. More effort is needed to come up with new construction and tighten the lower bound if possible. Examining how other structural assumptions can help us achieve faster estimation rates is also an interesting topic. For example, although the covariate may take many categories, the conditional means $\mu_{ak}$ may be constant over a small number of (known or unknown) groups:
\[
\mu_{ak} = \mu_{a,g(k)}, \qquad g(k) \in [G], \qquad G \ll d.
\]
When the grouping map $g(\cdot)$ is known, the problem reduces from estimating $2d$ stratum-specific means to estimating only $2G$ group-level means, and Theorem 1 then yields the corresponding rate $G^2/n^2 + 1/n$ (with $d$ replaced by $G$). When the grouping is unknown, one can encourage a small number of distinct values among the $\mu_{ak}$ by solving
\[
\min_{\mu_{a1},\dots,\mu_{ad}}
\sum_{k=1}^d \bigl(\hat{\mu}_{ak} - \mu_{ak} \bigr)^2
\;+\;
\lambda \sum_{k < \ell} \bigl|\mu_{ak} - \mu_{a\ell}\bigr|,
\]
motivated by the literature on convex clustering with fusion penalties
\citep{hocking2011clusterpath,chi2015splitting,tan2015statistical,radchenko2017convex,sun2021convex}.
The fusion term shrinks stratum-specific means toward one another and can therefore recover, or approximately recover, latent grouping structure when many of the $\mu_{ak}$ are equal or nearly equal. The resulting estimates of $\mu_{ak}$ can then be substituted into the plug-in ATE estimator considered in Section 3. Similarly, sparsity may also arise in the effect modifiers, in the sense that the category-specific treatment effects
    \[
    \tau_k = \mu_{1k} - \mu_{0k}
    \]
    are constant across a small number of groups. In this case, one can regularize the $\tau_k$ directly by solving
    \[
    \min_{\tau_1,\dots,\tau_d}
    \sum_{k=1}^d \bigl(\hat{\tau}_{k} - \tau_{k} \bigr)^2
    \;+\;
    \lambda \sum_{k < \ell} \bigl|\tau_{k} - \tau_{\ell}\bigr|,
    \]
    which encourages a small number of distinct values among the $\tau_k$.

Other extensions could involve the estimation of different causal estimands in a similar discrete setting, including time-varying treatment effects, generalizability and transportability, optimal treatment regimes, instrumental variable and more. Developing a comprehensive understanding on the properties of popular estimators and fundamental limits of different causal functionals in the discrete setting could be an interesting avenue left for future investigation.


\bibliographystyle{apalike}
\bibliography{refer}

\newpage
\appendix
\begin{center}
{\Large\bf Appendix}
\end{center}

\section{Asymptotic Normality of the Plug-in Estimator}\label{sec:CLT}

In this section, we provide a central limit theorem for the plug-in estimator when the number of categories $d$ is fixed as a constant, summarized in the following theorem. 

\begin{theorem}\label{thm:plugin-lowdim}
    Supposed $X \in [d]$ is discrete with $d$ fixed and the nuisance estimators are the empirical averages defined in \eqref{eq:mle-nuisance}. Then we have
    \[
    \sqrt{n}\left(\widehat{\psi}-\psi\right) \stackrel{d}{\rightarrow} N\left(0, \operatorname{Var}(\varphi(\bZ))\right),
    \]
    where 
    \[
    \varphi(\bZ) = \mu_{1X}-\mu_{0X}+\left (\frac{A}{\pi_X}-\frac{1-A}{1-\pi_X}\right )\left (Y-\mu_{AX}\right )
    \]
    is the first-order influence function of ATE $\psi$ under a nonparametric model.
\end{theorem}

Theorem \ref{thm:plugin-lowdim} implies that when the covariate is discrete with fixed dimension $d$, the plug-in-style estimator $\hat{\psi}$ is $\sqrt{n}$-consistent and asymptotically normal. Hence in low-dimensional problems, the plug-in estimator $\hat{\psi}$ enjoys appealing properties and we can construct confidence intervals and perform statistical tests on ATE based on Theorem \ref{thm:plugin-lowdim}. It is worth noting that asymptotic normality also holds in the regime $d=o(\sqrt{n})$. However, truncating the propensity score estimates at $\epsilon$ and $1-\epsilon$ is required to avoid instability induced by imprecise estimation of $\pi_k$ in the high-dimensional regime. Moreover, the empirical influence function $\hat{\varphi}$ belongs to a Donsker class when $d$ is fixed since it could be expressed as finite-dimensional parametric models (the dimension depends on $d$). As $d=d_n$ grows with $n$, they may not belong to a Donsker class and sample splitting is required to control the empirical process term.

\section{Additional Simulation Results}

In this section, we present additional simulation results on the asymptotic normality of the estimators studied in the main text, as well as on the accuracy of bootstrap-based variance estimation.

\subsection{Asymptotic Normality}

\begin{figure}[t]
  \centering
  \subfigure[$n=1000$]{
  \begin{minipage}[t]{0.48\linewidth}
\centering
    \includegraphics[width=\linewidth]{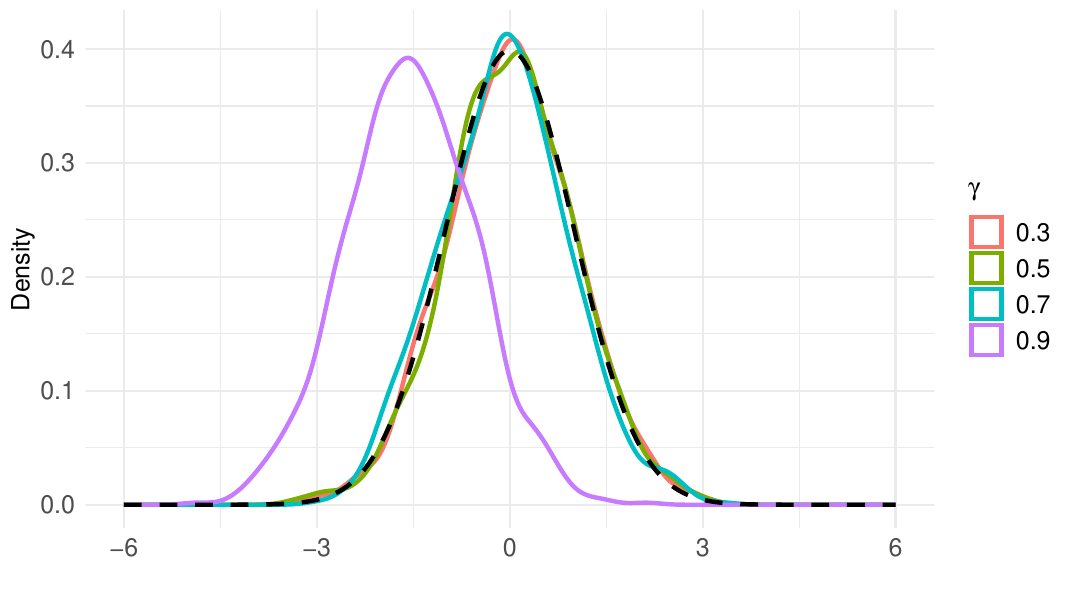}
  \end{minipage}}\hfill
  \subfigure[$n=10000$]{
  \begin{minipage}[t]{0.48\linewidth}
\centering
    \includegraphics[width=\linewidth]{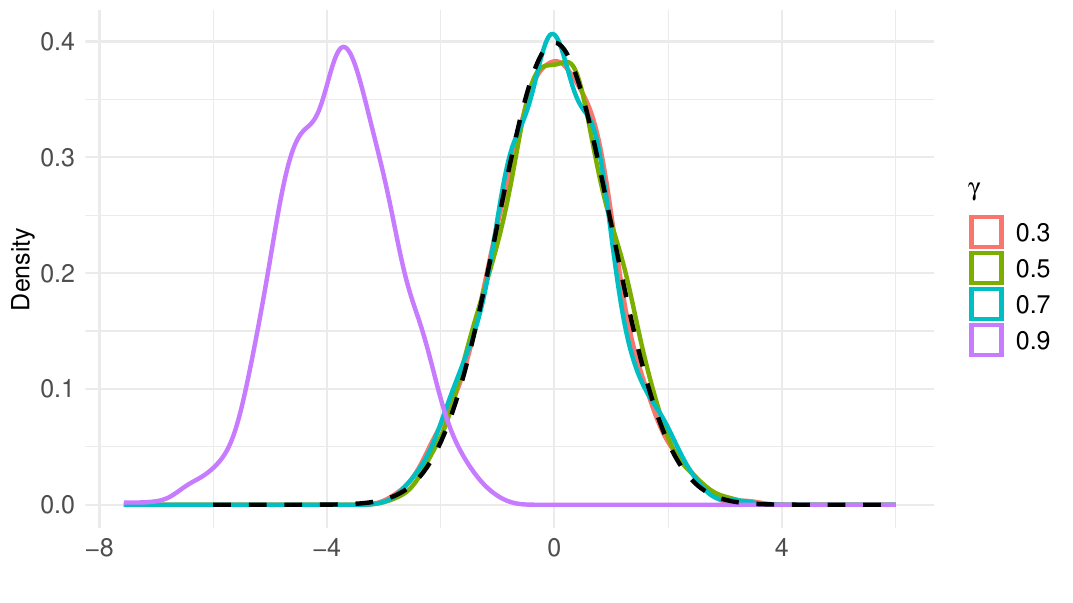}
  \end{minipage}}
  \caption{Density of standardized plug-in estimator $\hat{\psi}_1$. The black dashed line is the standard normal density.}
  \label{fig:simu-psi1}
\end{figure}

\begin{figure}[t]
  \centering
  \subfigure[$n=1000$]{
  \begin{minipage}[t]{0.48\linewidth}
\centering
    \includegraphics[width=\linewidth]{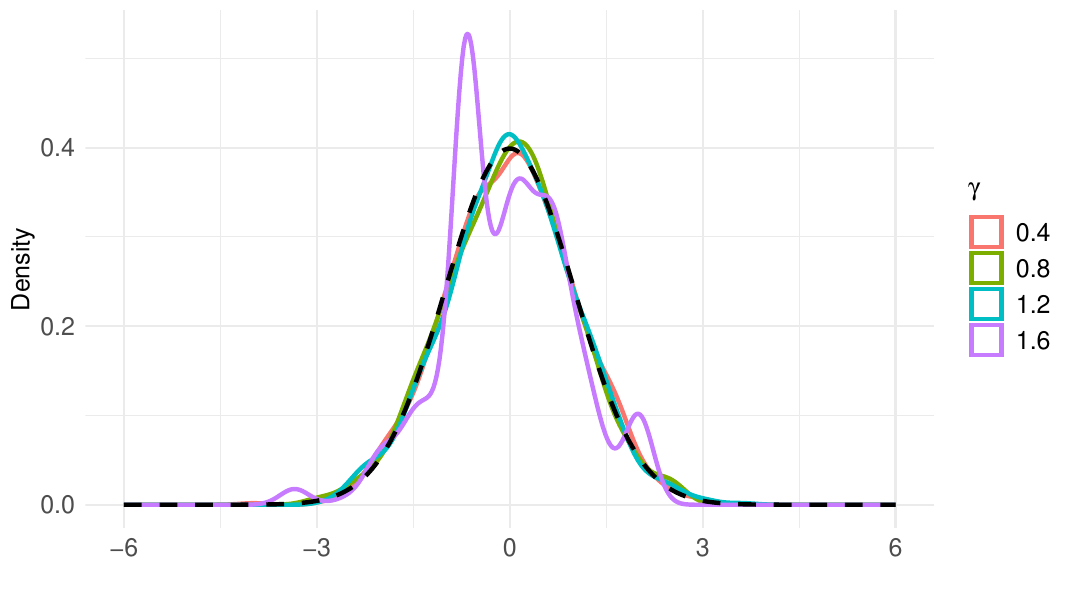}
  \end{minipage}}\hfill
  \subfigure[$n=10000$]{
  \begin{minipage}[t]{0.48\linewidth}
\centering
    \includegraphics[width=\linewidth]{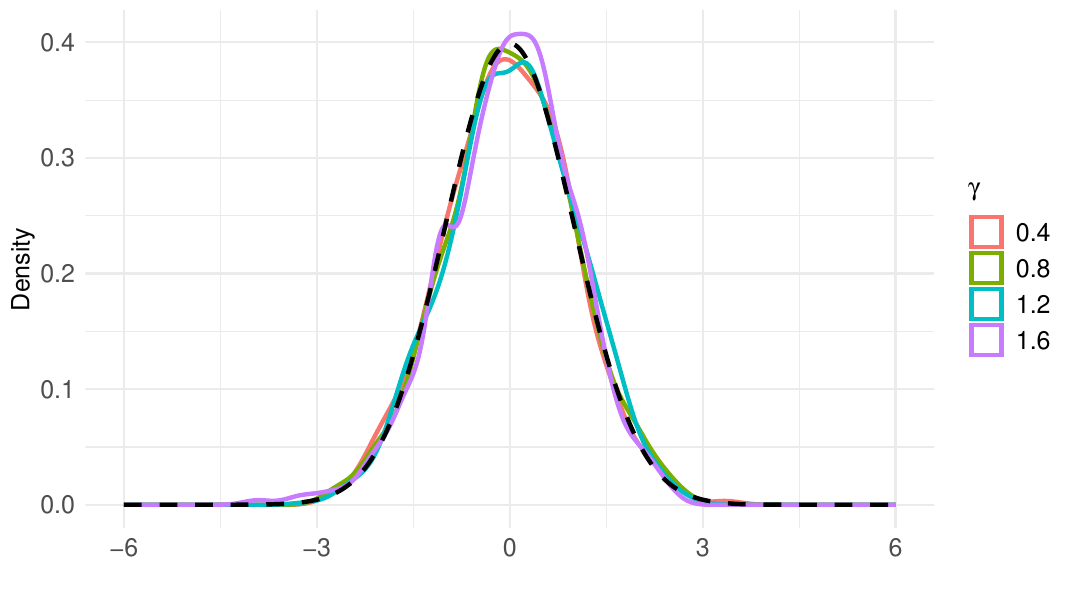}
  \end{minipage}}
  \caption{Density of the estimator $\hat{\tau}$ exploiting effect homogeneity. The black dashed line is the standard normal density.}
  \label{fig:simu-tau}
\end{figure}

\begin{figure}[t]
  \centering
  \subfigure[$n=1000$]{
  \begin{minipage}[t]{0.48\linewidth}
\centering
    \includegraphics[width=\linewidth]{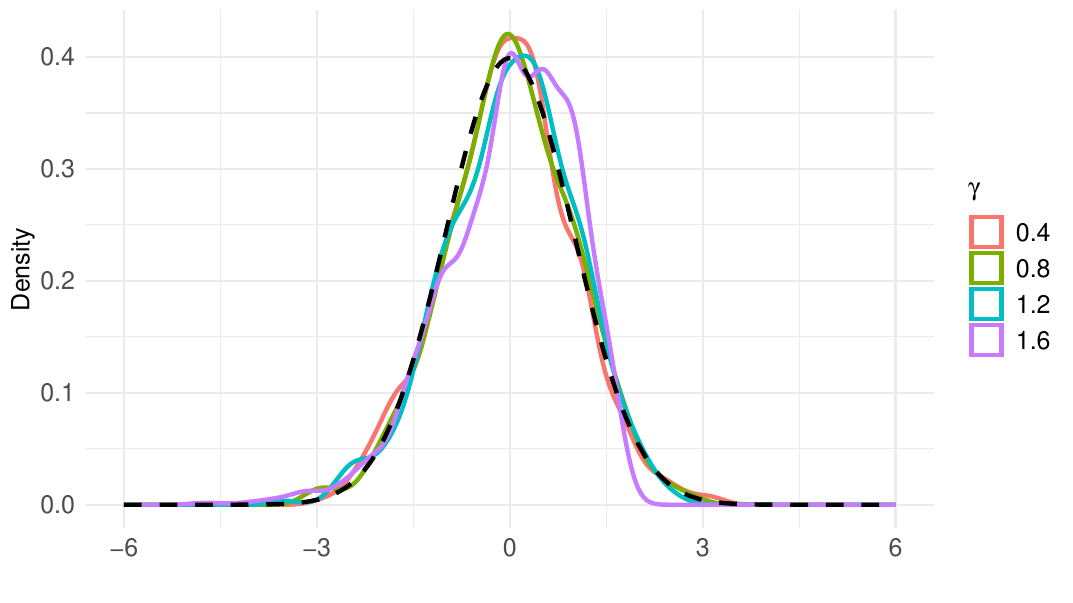}
  \end{minipage}}\hfill
  \subfigure[$n=10000$]{
  \begin{minipage}[t]{0.48\linewidth}
\centering
    \includegraphics[width=\linewidth]{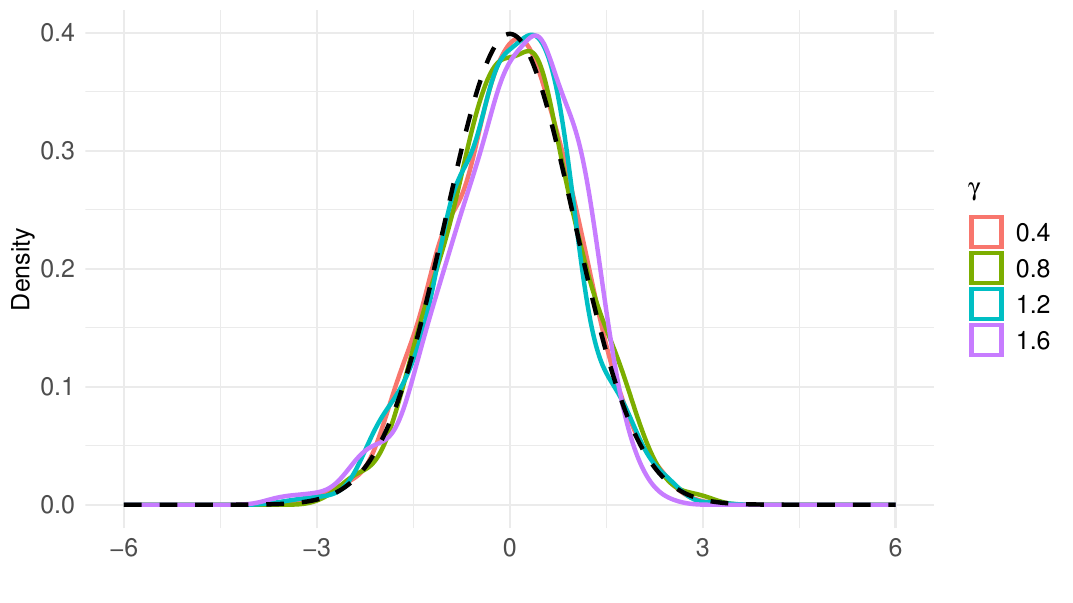}
  \end{minipage}}
  \caption{Density of the second-order estimator $\hat{\eta}$. The black dashed line is the standard normal density.}
  \label{fig:simu-eta}
\end{figure}

We first illustrate the asymptotic normality of the estimators considered in the paper, including the plug-in estimator $\hat{\psi}_1$ from Section~\ref{sec:ate-plugin}, the estimator $\hat{\tau}$ exploiting effect homogeneity from Section~\ref{sec:homogeneity}, and the second-order estimator $\hat{\eta}$ from Section~\ref{sec:cova-distribution}. The data-generating process is the same as in Section~\ref{sec:simulation}. For each repetition $m \in [M]$, we compute a standardized version of the estimator. For example, the standardized plug-in estimator is $(\hat{\psi}_1^{(m)}-\psi_1)/\widehat{\mathrm{sd}}(\hat{\psi}_1)$, where $\widehat{\mathrm{sd}}(\hat{\psi}_1)$ denotes the empirical standard deviation of $\hat{\psi}_1$ across the $M$ repetitions. Throughout, we center the estimators at the true parameter values, rather than at their Monte Carlo means. We then plot the estimated densities of the standardized estimators for different choices of $\gamma$, recalling that $d=\lfloor n^\gamma\rfloor$. We focus on $n=1000$ and $n=10{,}000$. The results are summarized in Figures~\ref{fig:simu-psi1}--\ref{fig:simu-eta}.

In Figure~\ref{fig:simu-psi1}, when $\gamma \in \{0.3, 0.5, 0.7\}$, the estimated densities of the standardized plug-in estimator closely match the standard normal distribution, consistent with the asymptotic normality established in Theorem~\ref{thm:plugin-lowdim}. When $\gamma=0.9$, the density remains qualitatively bell-shaped, but it is no longer centered at zero because we center the estimator using the true parameter $\psi_1$ rather than its Monte Carlo mean. In this high-dimensional regime the plug-in estimator is biased, so centering at $\psi_1$ induces a visible shift in the standardized distribution.

In Figures~\ref{fig:simu-tau}--\ref{fig:simu-eta}, when $\gamma \in \{0.4, 0.8, 1.2\}$, the estimated densities of the standardized estimators $\hat{\tau}$ and $\hat{\eta}$ closely match the standard normal distribution. Since these estimators can attain faster convergence rates under the additional structure they exploit, the standardized distributions remain approximately normal even at $\gamma=1.2$. This contrasts with Figure~\ref{fig:simu-psi1}, where the plug-in estimator $\hat{\psi}_1$ is biased in the high-dimensional regime and its standardized distribution is noticeably shifted when $\gamma=0.9$. When $\gamma=1.6$, the densities deviate from normality for $n=1000$ but move closer to a standard normal shape for $n=10{,}000$, suggesting that a CLT may still hold in this more challenging regime. We leave a formal analysis of asymptotic normality for $\hat{\tau}$ and $\hat{\eta}$ as an interesting direction for future work.

\subsection{Bootstrap Variance Estimation}

\begin{figure}[t]
  \centering
  \subfigure[$n=1000$]{
  \begin{minipage}[t]{0.48\linewidth}
\centering
    \includegraphics[width=\linewidth]{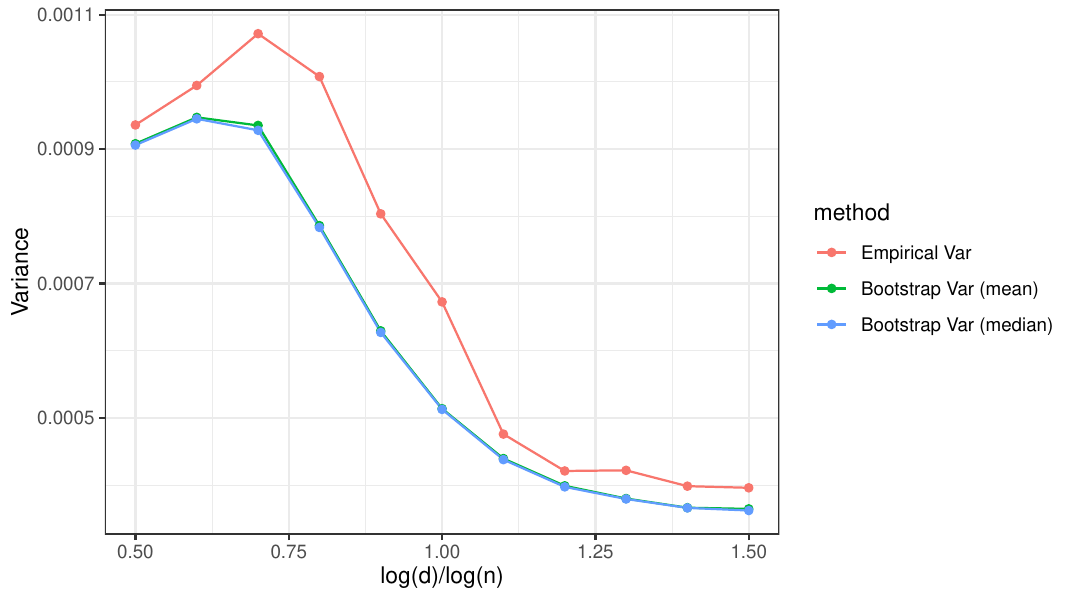}
  \end{minipage}}\hfill
  \subfigure[$n=10000$]{
  \begin{minipage}[t]{0.48\linewidth}
\centering
    \includegraphics[width=\linewidth]{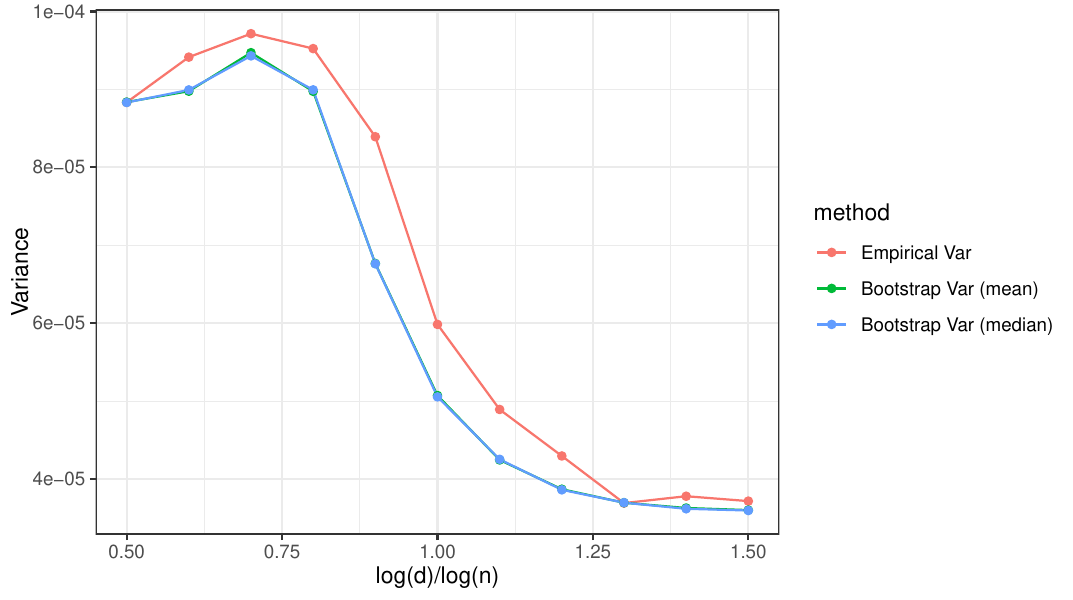}
  \end{minipage}}
  \caption{Comparison between estimated variance from bootstrap with the empirical variance}
  \label{fig:simu-var}
\end{figure}

We then examine the accuracy of bootstrap-based variance estimation in our high-dimensional setting with discrete covariates. We focus on estimating the variance of the plug-in estimator $\hat{\psi}_1$. The data-generating process is the same as in Section~\ref{sec:simulation}. For each replication $m \in [M]$, we estimate the variance of $\hat{\psi}_1^{(m)}$ using the nonparametric bootstrap with $B=500$ resamples. We then compare the mean and median of the resulting bootstrap variance estimates to the empirical variance of $\hat{\psi}_1^{(m)}$ across $M=500$ replications, for each $\gamma \in \{0.5,0.6,\ldots,1.5\}$. The results are summarized in Figure~\ref{fig:simu-var}.

In Figure~\ref{fig:simu-var}, the mean and median of the bootstrap variance estimates are generally smaller than the empirical variance of $\hat{\psi}_1$ across replications, which serves as a consistent estimate of the true sampling variance. This pattern holds for both $n=1000$ and $n=10{,}000$. Overall, these results suggest that, in our simulated high-dimensional discrete-covariate setting, the nonparametric bootstrap can systematically underestimate the variance of $\hat{\psi}_1$. When $d$ is large relative to $n$, many categories have $0$ or very small counts. The sampling distribution of plug-in estimators is then strongly influenced by whether a stratum appears at all and whether it has both treated/control observations. A bootstrap resample is drawn from the observed data, so it cannot recreate “new” strata that were absent in the original sample and tends to under-represent the variability coming from these appearance/disappearance events.
Developing more robust variance estimators in this regime is an interesting direction for future work.

\section{Proof of Main Results}

\subsection{Proof of Proposition \ref{prop:ests-equivalence}}
\begin{proof}
    We will prove three estimators for $\psi_1 = \ME[Y^1] = \sum_{k=1}^d {p}_k {\mu}_{1k}$ are equal and similar results hold for $\ME[Y^0] = \sum_{k=1}^d {p}_k {\mu}_{0k}$. First consider the regression estimator
    \[
    \begin{aligned}
     \hat{\psi}_{1,\text{reg}} = &\, \frac{1}{n} \sum_{i=1}^n \hat{\mu}_{1X_i} \\
     =&\, \frac{1}{n} \sum_{i=1}^n \sum_{k=1}^d I(X_i = k) \hat{\mu}_{1X_i} \\
     =&\, \frac{1}{n} \sum_{i=1}^n \sum_{k=1}^d I(X_i = k) \hat{\mu}_{1 k} \\
     = &\, \sum_{k=1}^d \left(\frac{1}{n} \sum_{i=1}^n I(X_i = k) \right) \hat{\mu}_{1 k} \\
     = &\, \sum_{k=1}^d \hat{p}_k\hat{\mu}_{1 k} = \hat{\psi},
    \end{aligned}
    \]
    where the second equation follows from the fact $\sum_{k=1}^d I(X_i=k)=1$. In fact, given $n$ samples, we have $\sum_{k=1}^d I(X_i=k)=\sum_{k: \hat{p}_k >0} I(X_i=k)=1$. With this in mind, for the inverse probability weighting estimator we have 
    \[
    \begin{aligned}
     \hat{\psi}_{1,\text{ipw}} = &\, \frac{1}{n}\sum_{i=1}^n \frac{A_i Y_i}{\hat{\pi}_{X_i}} \\
     = &\, \frac{1}{n}\sum_{i=1}^n \sum_{k:\hat{p}_k>0} I(X_i=k)\frac{A_i Y_i}{\hat{\pi}_{X_i}} \\
     = &\, \frac{1}{n}\sum_{i=1}^n \sum_{k:\hat{p}_k>0} I(X_i=k) \frac{A_i Y_i}{\hat{\pi}_{k}} \\ 
     = &\, \sum_{k:\hat{p}_k>0} \frac{1}{n}\sum_{i=1}^n  \frac{I(X_i=k)A_i Y_i}{\hat{\pi}_{k}} \\ 
     = &\, \sum_{k:\hat{p}_k>0} \hat{p}_k \frac{1}{n}\sum_{i=1}^n  \frac{I(X_i=k)A_i Y_i}{\hat{p}_k\hat{\pi}_{k}} \\
     = &\, \sum_{k:\hat{p}_k>0} \hat{p}_k \frac{\hat{q}_{1k}}{\hat{w}_k} \\
     = &\, \sum_{k:\hat{p}_k>0} \hat{p}_k \hat{\mu}_{1k} \\
     = &\, \sum_{k=1}^d \hat{p}_k \hat{\mu}_{1k} = \hat{\psi},
    \end{aligned}
    \]
    where $\hat{q}_{1k} = \frac{1}{n}\sum_{i=1}^n I(X_i=k)A_i Y_i$ and $\hat{q}_{1k}/\hat{w}_k = \hat{\mu}_{1k}$ by definition in \eqref{eq:mle-nuisance}. Note that all the equations still hold when some categories have no treated samples (i.e. $\hat{w}_k=0$) since we define $0/0=0$ whenever it appears. The last equation holds since the categories with $\hat{p}_k=0$ do not contribute to the estimation. Finally, we consider the doubly robust estimator. Note that
    \[
    \begin{aligned}
        &\, \MP_n \left [ \frac{A\hat{\mu}_{1X}}{\hat{\pi}_X} \right] \\
        =&\, \frac{1}{n} \sum_{i=1}^n \sum_{k:\hat{p}_k>0} I(X_i=k) \frac{A_i \hat{\mu}_{1X_i}}{\hat{\pi}_{X_i}} \\
        =&\, \sum_{k:\hat{p}_k>0} \left( \frac{\sum_{i=1}^n A_i I(X_i=k)}{n \hat{p}_k \hat{\pi}_k} \right) \hat{p}_k \hat{\mu}_{1k}.
    \end{aligned}
    \]
    Note that $n \hat{p}_k \hat{\pi}_k = \sum_{i=1}^n A_i I(X_i=k)$. If $n \hat{p}_k \hat{\pi}_k>0$ then
    \[
    \frac{\sum_{i=1}^n A_i I(X_i=k)}{n \hat{p}_k \hat{\pi}_k} = 1.
    \]
    If $n \hat{p}_k \hat{\pi}_k>0$ by our definition on $0/0=0$
    \[
    \frac{\sum_{i=1}^n A_i I(X_i=k)}{n \hat{p}_k \hat{\pi}_k} = 0, \,\hat{\mu}_{1k} = \frac{\hat{q}_{1k}}{\hat{w}_k} = 0
    \]
    and hence we always have
    \[
    \left( \frac{\sum_{i=1}^n A_i I(X_i=k)}{n \hat{p}_k \hat{\pi}_k} \right) \hat{p}_k \hat{\mu}_{1k} = \hat{p}_k \hat{\mu}_{1k}.
    \]
    We conclude
    \[
    \MP_n \left [ \frac{A\hat{\mu}_{1X}}{\hat{\pi}_X} \right] = \sum_{k: \hat{p}_k >0} \hat{p}_k \hat{\mu}_{1k} = \sum_{k=1}^d \hat{p}_k \hat{\mu}_{1k}.
    \]
    This together with 
    \[
    \MP_n \left[\frac{AY}{\hat{\pi}_X} \right] = \hat{\psi}
    \]
    shows
    \[
    \hat{\psi}_{1,\text{dr}} = \MP_n \left[ \frac{A(Y-\hat{\mu}_{1X})}{\hat{\pi}_X} + \hat{\mu}_{1X}\right] = \hat{\psi}.
    \]
\end{proof}

\subsection{Proof of Proposition \ref{thm:ate-exact-bias}}
\begin{proof}
    We focus on the bias of $\hat{\psi}_1$ and $\hat{\psi}_0$ can be similarly analyzed. We may rewrite the plug-in estimator as
    \[
    \hat{\psi}_1 = \sum_{k=1}^d \frac{\hat{q}_{1k}\hat{p}_k}{\hat{w}_k} I(\hat{w}_k > 0)
    \]
    to emphasize the definition $0/0=0$ in each term. On the event $\{ n\hat{w}_k = 0\}$, we automatically have $\hat{q}_{1k} = 0$. On the event $\{ n\hat{w}_k > 0\}$, 
    \[
    \ME[n\hat{q}_{1k} \mid\bX^n, \bA^n] = n \hat{w}_k \mu_{1k} = n\frac{q_{1k}}{w_{k}}\hat{w}_k
    \]
    since $n\hat{q}_{1k} \mid \bX^n, \bA^n \sim \text{B}(n\hat{w}_k, \mu_{1k})$. Hence we have
    \[
    \ME\left[ \frac{\hat{q}_{1k}\hat{p}_k}{\hat{w}_k} I(\hat{w}_k > 0) \right] = \ME\left[ \frac{\hat{p}_k}{\hat{w}_k} I(\hat{w}_k > 0) \ME(\hat{q}_{1k}\mid\bX^n, \bA^n) \right] = \ME[\hat{p}_k I(\hat{w}_k>0)]\frac{q_{1k}}{w_k}.
    \]
    The bias of each individual term is 
    \[
    \begin{aligned}
        &\, \ME\left[ \frac{\hat{q}_{1k}\hat{p}_k}{\hat{w}_k} I(\hat{w}_k > 0) \right]-\frac{q_{1k}p_k}{w_k} \\
        = &\, \ME[\hat{p}_k I(\hat{w}_k>0)] \frac{q_{1k}}{w_k} - \frac{q_{1k}p_k}{w_k} \\
        = &\, -\ME[\hat{p}_k I(\hat{w}_k=0)]\frac{q_{1k}}{w_{k}}
    \end{aligned}  
    \]
    Further note that $n\hat{w}_k \mid \bX^n \sim \text{B}(n\hat{p}_k, \pi_k)$, 
    \[
    \begin{aligned}
     &\, \ME[\hat{p}_k I(\hat{w}_k=0)] \\
     = &\,  \ME[\hat{p}_k I(n\hat{w}_k=0)] \\
     =&\, \ME[\hat{p}_k \MP(n\hat{w}_k =0\mid\bX^n)] \\
     = &\, \ME[\hat{p}_k (1-\pi_k)^{n\hat{p}_k}].
    \end{aligned}
    \]
    We then use the fact $\ME[Vc^V] = npc(1-p+cp)^{n-1} $ for $V \sim \text{Bin}(n,p)$ and obtain
    \[
    \begin{aligned}
      &\,  \ME[\hat{p}_k (1-\pi_k)^{n\hat{p}_k}] \\
      = &\, \frac{1}{n}\ME[n\hat{p}_k (1-\pi_k)^{n\hat{p}_k}] \\
      = &\, \frac{1}{n} \times n p_k (1-\pi_k) [1-p_k+(1-\pi_k)p_k]^{n-1} \\
      = &\,  p_k (1-\pi_k)(1-p_k \pi_k)^{n-1} 
    \end{aligned}
    \]
    Hence the bias for individual terms is
    \[
    \ME\left[ \frac{\hat{q}_{1k}\hat{p}_k}{\hat{w}_k} I(\hat{w}_k > 0) \right]-\frac{q_{1k}p_k}{w_k}= -\mu_{1k}p_k (1-\pi_k)(1-p_k \pi_k)^{n-1}
    \]
    and we conclude 
    \[
    \ME[\hat{\psi}_1 - \psi_1] = -\sum_{k=1}^d \mu_{1k}p_k (1-\pi_k)(1-p_k \pi_k)^{n-1}.
    \]
    Similarly (or by symmetry) for $\hat{\psi}_0$ we have
    \[
    \ME[\hat{\psi}_0 - \psi_0] = -\sum_{k=1}^d \mu_{0k}p_k \pi_k(1-p_k+p_k \pi_k)^{n-1}
    \]
\end{proof}

\subsection{Proof of Proposition \ref{prop:ate-worst-bias}}
\begin{proof}
Recall that we use bold letters $(\bp, \bpi, \bmu_1)$ to denote the vectors $(p_k, \pi_k, \mu_{1k})_{k=1}^d$. We have
\[
\begin{aligned}
    \sup_{\MP \in \mathcal{D}(\epsilon)}|\ME_{\MP}[\hat{\psi}_1 - \psi_1]| =&\, \sup_{\bp, \bpi, \bmu_1} \sum_{k=1}^d \mu_{1k}p_k (1-\pi_k)(1-p_k \pi_k)^{n-1} \\
    =&\, \sup_{ \bp, \bpi} \sum_{k=1}^d p_k (1-\pi_k)(1-p_k \pi_k)^{n-1} \\
    = &\, \sup_{ \bp } \sum_{k=1}^d p_k (1-\epsilon)(1-\epsilon p_k )^{n-1}, 
\end{aligned}
\]
where the first equation follows from $\mu_{1k} \leq 1$, and the second follows from $\pi_k \geq \epsilon$. The upper bound then follows from
\[
\max_{p\ge 0} p(1-\epsilon p)^{n-1} = \frac{1}{\epsilon n}\left(1-\epsilon\cdot \frac{1}{\epsilon n}\right)^{n-1} < \frac{1}{\epsilon n},
\]
where the maximizer $p^\star = 1/(\epsilon n)$ follows from simple differentiation. The lower bound follows from the enumeration of three different scenarios: 
\begin{enumerate}
    \item If $n\epsilon < 1$, choosing $\bp = (1,0,\cdots,0)$ gives the lower bound
    \begin{align*}
        (1-\epsilon)\cdot (1-\epsilon)^{n-1}  \ge (1-\epsilon) \left(1-\frac{1}{n}\right)^{n-1} \ge \frac{1-\epsilon}{e}. 
    \end{align*}
    \item If $1\le n\epsilon<d-1$, choosing $\bp = (1/(n\epsilon), \cdots, 1/(n\epsilon), 1-\lfloor n\epsilon \rfloor / (n\epsilon), 0, \cdots, 0)$ gives the lower bound
    \begin{align*}
        \lfloor n\epsilon \rfloor \cdot \frac{1-\epsilon}{n\epsilon}\left(1-\frac{1}{n}\right)^{n-1} \ge \frac{n\epsilon}{2}\cdot \frac{1-\epsilon}{en\epsilon} \ge \frac{1-\epsilon}{2e}.
    \end{align*}
    \item If $n\epsilon \ge d-1$, choosing $\bp = (1/(n\epsilon), \cdots, 1/(n\epsilon), 1-(d-1)/(n\epsilon))$ gives the lower bound
    \begin{align*}
        (d-1)\cdot \frac{1-\epsilon}{n\epsilon}\left(1-\frac{1}{n}\right)^{n-1} \ge \frac{1-\epsilon}{e}\cdot \frac{d-1}{n\epsilon}.
    \end{align*}
\end{enumerate}

\end{proof}

\subsection{Proof of Theorem \ref{thm:ate-var}}
\begin{proof}
\begin{equation}\label{eq:var-psi1}
    \Var(\hat{\psi}_1) = \ME[\Var(\hat{\psi}_1 \mid \bX^n, \bA^n)] + \Var(\ME[\hat{\psi}_1\mid\bX^n, \bA^n]).
\end{equation}
We start with the first term. By the conditional independence of $\hat{q}_{11}, \dots, \hat{q}_{1d}$ given $\bX^n, \bA^n$, we have
\[
\begin{aligned}
    &\, \Var(\hat{\psi}_1 \mid \bX^n, \bA^n)\\
    = &\, \Var \left( \sum_{k=1}^d \frac{\hat{q}_{1k}\hat{p}_k I(\hat{w}_k>0)}{\hat{w}_k} \mid \bX^n, \bA^n \right) \\
    =&\, \sum_{k=1}^d \frac{\hat{p}_k^2 I(\hat{w}_k >0)}{\hat{w}_k^2} \Var(\hat{q}_{1k}\mid\bX^n, \bA^n) \\
    = &\, \sum_{k=1}^d \frac{\hat{p}_k^2 I(\hat{w}_k >0)}{\hat{w}_k^2} \frac{\hat{w}_k \mu_{1k}(1-\mu_{1k})}{n} \\
    = &\, \frac{1}{n} \sum_{k=1}^d  \frac{\hat{p}_k^2 I(\hat{w}_k >0)}{\hat{w}_k}\mu_{1k}(1-\mu_{1k})
\end{aligned}
\]
where we use the fact $n \hat{q}_{1k} \mid \bX^n, \bA^n \sim \text{Bin}(n\hat{w}_k, \mu_{1k})$. We need the following lemma to proceed the analysis.

\begin{lemma}\label{lemma:bino-inverse}
        (Lemma A.2 in \cite{devroye2013probabilistic})If $X \sim B(n,p)$, then
        \[
        \mathbb{E}\left\{\frac{I(X>0)}{X}\right\} \leq \frac{2}{(n+1) p}.
        \]
    \end{lemma}
Note that $n\hat{w}_k \sim \text{Bin}(n\hat{p}_k, \pi_k) \mid \bX^n$ and apply Lemma \ref{lemma:bino-inverse}, we have
\[
\begin{aligned}
    &\, \ME[\Var(\hat{\psi}_1 \mid \bX^n, \bA^n)] \\
    = &\, \frac{1}{n} \sum_{k=1}^d \ME \left[ \frac{\hat{p}_k^2 I(\hat{w}_k >0)}{\hat{w}_k} \right]\mu_{1k}(1-\mu_{1k}) \\
    = &\, \sum_{k=1}^d \ME \left\{ \hat{p}_k^2 \ME \left[ \frac{ I(\hat{w}_k >0)}{n\hat{w}_k} \mid \bX^n\right] \right\}\mu_{1k}(1-\mu_{1k}) \\
    \leq &\, \sum_{k=1}^d \ME \left\{ \hat{p}_k^2 \frac{2}{(n \hat{p}_k +1)\pi_k} \right\}\mu_{1k}(1-\mu_{1k})\\
    \leq &\,  \frac{2}{n} \sum_{k=1}^d \ME \left\{  \frac{\hat{p}_k}{\pi_k} \right\}\mu_{1k}(1-\mu_{1k}) \\
    \leq &\, \frac{1}{2n\epsilon} \sum_{k=1}^d \ME[\hat{p}_k] \\
    = &\, \frac{1}{2n\epsilon},
\end{aligned}
\]
where we use the fact $\mu_{1k}(1-\mu_{1k}) \leq 1/4$. Then we evaluate the second term in \eqref{eq:var-psi1}
\[
    \ME[\hat{\psi}_1 \mid \bX^n, \bA^n] = \sum_{k=1}^d \ME \left[ \frac{\hat{q}_{1k}\hat{p}_k}{\hat{w}_k}I(\hat{w}_k >0)  \mid \bX^n, \bA^n \right] = \sum_{k=1}^d  \hat{p}_kI(\hat{w}_k >0) \mu_{1k}.
\]
By further using the property of conditional variance we have
\begin{equation}\label{eq:var-psi1_cond}
    \Var(\ME[\hat{\psi}_1 \mid \bX^n, \bA^n]) =  \ME\left[\Var\left(\sum_{k=1}^d  \hat{p}_kI(\hat{w}_k >0) \mu_{1k}  \mid \bX^n \right)\right] + \Var\left(\ME\left[\sum_{k=1}^d  \hat{p}_kI(\hat{w}_k >0) \mu_{1k}  \mid \bX^n \right]\right)
\end{equation}
We analyze the first term in \eqref{eq:var-psi1_cond}. By the conditional independence of $n\hat{w}_1 , \dots, n \hat{w}_d$ given $\bX^n$ we have
\[
\begin{aligned}
    &\, \Var\left(\sum_{k=1}^d  \hat{p}_kI(\hat{w}_k >0) \mu_{1k} \mid \bX^n \right)\\
    =&\, \sum_{k=1}^d \hat{p}_k^2 \mu_{1k}^2 \Var(I(\hat{w}_k > 0)\mid\bX^n)\\
    \leq &\, \sum_{k=1}^d \hat{p}_k^2 \mu_{1k}^2 (1-\pi_k)^{n\hat{p}_k} \\
    \leq &\, \sum_{k=1}^d \hat{p}_k^2 \mu_{1k}^2 \frac{1}{\pi_k(1+n\hat{p}_k)} \\
    \leq &\, \frac{1}{\epsilon n} \sum_{k=1}^d \hat{p}_k \\
    = &\, \frac{1}{\epsilon n},
\end{aligned}
\]
where we use the fact
\[
(1-x)^n \leq \frac{1}{x(1+n)}, \, \forall \, 0\leq x \leq 1.
\]
Hence we have
\[
\ME\left[\Var\left(\sum_{k=1}^d  \hat{p}_kI(\hat{w}_k >0) \mu_{1k} \mid \bX^n \right)\right] \leq \frac{1}{\epsilon n}.
\]
For the second term in \eqref{eq:var-psi1_cond}, we have
\[
\begin{aligned}
    &\, \ME\left[\sum_{k=1}^d  \hat{p}_kI(\hat{w}_k >0) \mu_{1k}  \mid \bX^n \right] \\
    = &\, \sum_{k=1}^d \hat{p}_k \mu_{1k}(1-\MP(n\hat{w}_k=0\mid\bX^n)) \\
    = &\, \sum_{k=1}^d \mu_{1k} \hat{p}_k [1-(1-\pi_k)^{n\hat{p}_k}]
\end{aligned}
\]
So we need to evaluate 
\[
\begin{aligned}
    &\, \Var\left(\ME\left[\sum_{k=1}^d  \hat{p}_kI(\hat{w}_k >0) \mu_{1k} \mid \bX^n \right]\right) \\
    = &\, \Var \left( \sum_{k=1}^d \mu_{1k} \hat{p}_k [1-(1-\pi_k)^{n\hat{p}_k}]\right) \\
    = &\,  \sum_{k=1}^d \mu_{1k}^2 \Var\left(\hat{p}_k [1-(1-\pi_k)^{n\hat{p}_k}]\right) + \sum_{i\neq j} \mu_{1i}\mu_{1j}\operatorname{Cov}\left(\hat{p}_i [1-(1-\pi_i)^{n\hat{p}_i}], \hat{p}_j [1-(1-\pi_j)^{n\hat{p}_j}]\right).
\end{aligned}
\]
We note that
\[
\ME[\hat{p}_k [1-(1-\pi_k)^{n\hat{p}_k}]] = p_k - p_k(1-\pi_k)(1-\pi_k p_k)^{n-1}
\]
and 
\[
\ME\left\{\hat{p}_k^2 [1-(1-\pi_k)^{n\hat{p}_k}]^2\right\} \leq \ME[\hat{p}_k^2] = \frac{p_k(1-p_k)}{n} + p_k^2.
\]
This yields the following bound
\[
\begin{aligned}
    &\, \Var(\hat{p}_k [1-(1-\pi_k)^{n\hat{p}_k}]) \\
    \leq &\, \frac{p_k(1-p_k)}{n} + 2p_k^2 (1-\pi_k)(1-\pi_k p_k)^{n-1} \\
    \leq &\, \frac{p_k(1-p_k)}{n} + 2p_k^2 (1-\pi_k) \frac{1}{n\pi_k p_k} \\
    \leq &\, \frac{p_k(1-p_k)}{n} + \frac{2p_k}{n \epsilon} \\
    \leq  &\, \left( 1 + \frac{2}{\epsilon} \right) \frac{p_k}{n}.
\end{aligned}
\]
So we have 
\[
\sum_{k=1}^d \mu_{1k}^2 \Var\left(\hat{p}_k [1-(1-\pi_k)^{n\hat{p}_k}]\right) \leq \left( 1 + \frac{2}{\epsilon} \right) \frac{1}{n}. 
\]
For the covariance part, the computations are more involved. We need to compute
\[
\begin{aligned}
    &\, \ME\left\{ \hat{p}_i[1-(1-\pi_i)^{n\hat{p}_i}] \hat{p}_j [1-(1-\pi_j)^{n\hat{p}_j}] \right\} \\
    = &\, \ME[\hat{p}_i \hat{p}_j] -\ME[\hat{p}_i\hat{p}_j (1-\pi_i)^{n\hat{p}_i}]  -\ME[\hat{p}_i\hat{p}_j (1-\pi_j)^{n\hat{p}_j}]  + \ME[\hat{p}_i\hat{p}_j (1-\pi_i)^{n\hat{p}_i} (1-\pi_j)^{n\hat{p}_j}]
\end{aligned}
\]
For three-dimensional multinomial distribution $(X_1, X_2, X_3)\sim$ Multinomial$(n, p_1, p_2, p_3)$ the probability generating function is
\[
\ME[z_1^{X_1}z_2^{X_2}z_3^{X_3}] = (p_1z_1 + p_2z_2 + p_3z_3)^n.
\]
From this formula and differentiation we have
\[
\ME[X_1X_2 z_1^{X_1}z_2^{X_2}] = n(n-1)p_1 z_1 p_2 z_2 (p_1z_1 + p_2z_2 + p_3)^{n-2}.
\]
We can obtain the expectations appearing in the covariance as:
\[
\ME[\hat{p}_i \hat{p}_j] = \frac{n-1}{n}p_i p_j
\]
\[
\ME[\hat{p}_i\hat{p}_j (1-\pi_i)^{n\hat{p}_i}] = \frac{n-1}{n}p_i p_j (1-\pi_i)(1-\pi_i p_i)^{n-2}
\]
\[
\ME[\hat{p}_i\hat{p}_j (1-\pi_i)^{n\hat{p}_i} (1-\pi_j)^{n\hat{p}_j}] =  \frac{n-1}{n} p_i p_j (1-\pi_i)(1-\pi_j)(1-\pi_i p_i - \pi_j p_j)^{n-2}
\]
Now the covariance of pair $(i,j)$ is
\[
\begin{aligned}
    &\, \operatorname{Cov}\left(\hat{p}_i [1-(1-\pi_i)^{n\hat{p}_i}], \hat{p}_j [1-(1-\pi_j)^{n\hat{p}_j}]\right)\\
    = &\, \ME\left[ \hat{p}_i[1-(1-\pi_i)^{n\hat{p}_i}] \hat{p}_j [1-(1-\pi_j)^{n\hat{p}_j}] \right] - \ME \left[ \hat{p}_i[1-(1-\pi_i)^{n\hat{p}_i}] \right] \ME \left[ \hat{p}_j [1-(1-\pi_j)^{n\hat{p}_j}]\right] \\
    = &\, \frac{n-1}{n}p_i p_j - \frac{n-1}{n}p_i p_j (1-\pi_i)(1-\pi_i p_i)^{n-2} - \frac{n-1}{n}p_i p_j (1-\pi_j)(1-\pi_j p_j)^{n-2} \\
    &\, +\frac{n-1}{n} p_i p_j (1-\pi_i)(1-\pi_j)(1-\pi_i p_i - \pi_j p_j)^{n-2} \\
    &\, - [p_i - p_i(1-\pi_i)(1-\pi_i p_i)^{n-1}][p_j - p_j(1-\pi_j)(1-\pi_j p_j)^{n-1}] \\
    =&\, -\frac{1}{n}p_ip_j + p_i p_j (1-\pi_i)(1-\pi_ip_i)^{n-2}\left( \frac{1}{n} - \pi_i p_i \right) + p_i p_j (1-\pi_j)(1-\pi_jp_j)^{n-2}\left( \frac{1}{n} - \pi_j p_j \right) \\
    &\, + \frac{1}{n}p_i p_j (1-\pi_i)(1-\pi_j)\left[(n-1)(1-\pi_ip_i-\pi_jp_j)^{n-2} - n(1-\pi_ip_i)^{n-1}(1-\pi_jp_j)^{n-1} \right]
\end{aligned}
\]
where in the last equation we combine the four terms in $\ME\left[ \hat{p}_i[1-(1-\pi_i)^{n\hat{p}_i}] \hat{p}_j [1-(1-\pi_j)^{n\hat{p}_j}] \right]$ and $\ME \left[ \hat{p}_i[1-(1-\pi_i)^{n\hat{p}_i}] \right] \ME \left[ \hat{p}_j [1-(1-\pi_j)^{n\hat{p}_j}]\right]$ correspondingly. We proceed as taking summations
\[
\begin{aligned}
  &\,\bigg| \sum_{i \neq j} \mu_{1i} \mu_{1j} \operatorname{Cov}\left(\hat{p}_i [1-(1-\pi_i)^{n\hat{p}_i}], \hat{p}_j [1-(1-\pi_j)^{n\hat{p}_j}]\right) \bigg| \\
  \leq &\,\sum_{i \neq j} \frac{1}{n}p_i p_j + 2 \sum_{i \neq j}p_i p_j (1-\pi_i)(1-\pi_i p_i)^{n-2}\bigg| \frac{1}{n}-\pi_i p_i \bigg| \\
  &\, + \sum_{i \neq j} \frac{1}{n}p_i p_j (1-\pi_i)(1-\pi_j)\big |(n-1)(1-\pi_ip_i-\pi_jp_j)^{n-2} - n(1-\pi_ip_i)^{n-1}(1-\pi_jp_j)^{n-1} \big |
\end{aligned}
\]
For each individual term in the expression above, we have
\[
\sum_{i \neq j} \frac{1}{n}p_i p_j = \sum_{i} \frac{1}{n}p_i (1-p_i) \leq \frac{1}{n}.
\]
\[
\begin{aligned}
    &\, \sum_{i \neq j}p_i p_j (1-\pi_i)(1-\pi_i p_i)^{n-2}\bigg| \frac{1}{n}-\pi_i p_i\bigg|\\
    \leq &\, \sum_{i \neq j} \frac{1}{n}  p_i p_j + \sum_{i \neq j} p_i p_j (1-\pi_i)(1-\pi_i p_i)^{n-2} \pi_i p_i\\
    \leq &\, \frac{1}{n} + \sum_{i \neq j} p_i p_j (1-\pi_i)\frac{1}{\pi_i p_i (n-1)} \pi_i p_i \\
    \leq &\, \frac{1}{n}+ \frac{1}{n-1}.
\end{aligned}
\]
For the last term we need an auxiliary lemma as follows:
\begin{lemma}\label{lemma:bound-var}
    For $n \geq 3$, consider the function
    \[
    f_2(x,y) = (n-1)(1-x-y)^{n-2} - n(1-x)^{n-1}(1-y)^{n-1}
    \]
    defined on the triangle $D = \{ x \geq 0, y \geq 0, x+y\leq 1\}$. Then we have
    \[
    |f_2(x,y)| \leq 1.
    \]
\end{lemma}
For the last term we invoke Lemma \ref{lemma:bound-var} and have
\[
\begin{aligned}
    &\, \sum_{i \neq j} \frac{1}{n}p_i p_j (1-\pi_i)(1-\pi_j)\big |(n-1)(1-\pi_ip_i-\pi_jp_j)^{n-2} - n(1-\pi_ip_i)^{n-1}(1-\pi_jp_j)^{n-1} \big |\\
    \leq &\, \sum_{i \neq j} \frac{1}{n}p_i p_j \leq \frac{1}{n}.
\end{aligned}
\]
We thus showed
\[
\big| \sum_{i \neq j} \mu_{1i} \mu_{1j} \operatorname{Cov}\left(\hat{p}_i [1-(1-\pi_i)^{n\hat{p}_i}], \hat{p}_j [1-(1-\pi_j)^{n\hat{p}_j}]\right) \big| \leq \frac{4}{n} + \frac{2}{n-1}.
\]
\[
\Var\left(\ME\left[\sum_{k=1}^d  \hat{p}_kI(\hat{w}_k >0) \mu_{1k} \bigg | \bX^n \right]\right) \leq \left( 1 + \frac{2}{\epsilon} \right) \frac{1}{n} + \frac{4}{n} + \frac{2}{n-1}.
\]
\[
\Var(\ME[\hat{\psi}_1 \mid \bX^n, \bA^n]) \leq \left( 1 + \frac{3}{\epsilon} \right) \frac{1}{n} + \frac{4}{n} + \frac{2}{n-1}.
\]
\[
\Var(\hat{\psi}_1) \leq \left( 1 + \frac{7}{2\epsilon} \right) \frac{1}{n} + \frac{4}{n} + \frac{2}{n-1}.
\]
\end{proof}

\subsection{Proof of Theorem \ref{thm:ate-minimax}}
\begin{proof}
The proof strategy follows from the moment matching method commonly used in theoretical computer science literature \citep{wu2016minimax, wu2019chebyshev, jiao2015minimax}. We first define a ``relaxed" model class with proportion vector $\bp = (p_1, p_2,\dots)$ being \textit{approximately} a probability vector. Mathematically, for $\delta>0$ define
    \begin{equation}\label{eq:model-class-relaxed}
    \mathcal{D}(\epsilon, \delta) = \left\{ \left|\sum_{k=1}^d p_k -1\right| \leq \delta, p_k \in  \left[0,1\right], \pi_k \in [\epsilon, 1-\epsilon], \forall k \in [d] \right\},
    \end{equation}
    where we relax the assumption $\sum_k p_k=1$ to $\left|\sum_k p_k -1\right| \leq \delta$ so that we can set $p_k $'s to be random variables in our construction and use the method of fuzzy hypotheses \citep{Tsybakov2009,le2012asymptotic}. Under $\MP \in  \mathcal{D}(\epsilon, \delta)$, we assume $n\hat{p}_k \sim \text{Poi}(np_k)$ and $(n\hat{p}_1, n\hat{p}_2, \dots)$ are independent, i.e. we again rely on a Poisson sampling model to prove our results. The treatment assignment $A$ and outcome $Y$ have the same distribution as in Section \ref{sec:DGP} conditioned on the category each sample falls into. Define the minimax lower bound over $\mathcal{D}(\epsilon, \delta)$ as
    \[
    \tilde{R}^*(d,n,\delta) := \inf_{\hat{\psi}_1} \sup_{\MP \in \mathcal{D}(\epsilon, \delta)} \ME_{\MP}\left[\left(\hat{\psi}_1 - \psi_1\right)^2\right]
    \]
    The following lemma allows us to relate ${R}^*(d,n)$ with $\tilde{R}^*(d,n,\delta)$.
    \begin{lemma}\label{lemma:poission-minimax}
        For any $\delta \in [0,1/3)$,
        \[
        {R}^*(d,n/2) \geq \frac{1}{2} \tilde{R}^*(d,n,\delta) - \exp(-n/50)-\delta^2.
        \]
    \end{lemma}
    We then present an auxiliary result characterizing the prior distribution on the parameters in the method of fuzzy hypotheses. In the following proof of this section, we will set $\mu_{0k}=0$ and abbreviate $\mu_{1k}$ as $\mu_k$, which is different from the mean of outcome in $k$-th category $\ME[Y|X=k]$ as in Section \ref{sec:cova-distribution}. 
    \begin{lemma}\label{lemma:prior-dist}
        There exists constants $c, c' >0$ such that for any constants $c_1, c_2, c_3>0$ satisfying $\frac{c c_1 c_3}{c_2^2} \leq 1$ and $d \leq c_3 n \log n$, there exist two distributions $\mu_0 $ and $ \mu_1$ on $(p, \pi, \mu )$ satisfying the following properties:
        \begin{enumerate}
            \item $\mu_0$ a.s. and $\mu_1$ a.s. 
            \[
            0 \leq p \leq \frac{c_1 \log n}{n}, \, \epsilon \leq \pi \leq 1-\epsilon, \, 0 \leq \mu \leq 1;
            \]
            \item For $i, j ,k \geq 0$ and $i+j+k \leq 3K$ with $K = c_2 \log n$,
            \[
            \ME_{\mu_0}\left[ p^i (p \pi)^j (p \pi \mu)^k \right] = \ME_{\mu_1} \left[ p^i (p \pi)^j (p \pi \mu)^k \right].
            \]
            \item 
            \[
            a_0 : =\ME_{\mu_0}[p] =  \ME_{\mu_1}[p] \leq \frac{1}{d}.
            \]
            \item 
            \[
            |\ME_{\mu_0}[p \mu] - \ME_{\mu_1}[p \mu]| \geq \frac{c_4}{n \log n},
            \]
            where $c_4 = \frac{cc'c_1}{c_2^2}$. 
        \end{enumerate}
    \end{lemma}
    The proof of Lemma \ref{lemma:poission-minimax} and \ref{lemma:prior-dist} are provided in Section \ref{sec:proof-auxiliary}. Note that in the construction of $\mu_i$, $\pi$ and $ \mu$ are actually functions of $p$. Under null hypothesis $P$, let $(p_k, \pi_k, \mu_k)$ i.i.d. $\sim \mu_0, 1\leq k \leq d$. We add one more category with $p_{d+1}=1-d a_0, \pi_{d+1}=\epsilon, \mu_{d+1}=0$. Under the alternative hypothesis $P'$, let $(p_k', \pi_k', \mu_k')$ i.i.d. $\sim \mu_1, 1\leq k \leq d$ and $p_{d+1}'=1-d a_0, \pi_{d+1}'=\epsilon, \mu_{d+1}'=0$. Obviously, adding one more category will not affect the final rate. The sufficient statistics for $\psi_1$ are 
    \[
    \#\{i:X_i=k,A_i=1,Y_i=1 \}, \#\{i:X_i=k,A_i=1,Y_i=0 \}, \#\{i:X_i=k,A_i=0 \}, 1\leq k \leq d+1.
    \]
    By the property of Poisson distribution these counting statistics are independent under the Poisson sampling model and under the null hypothesis $P$ (given $\{p_1, \dots, p_d\}$, which is equivalent to conditioning on $\{(p_k, \pi_k, \mu_k ), 1\leq k \leq d\}$ since $\pi, \mu$ are functions of $p$)
    \[
    N_{k11} : = \#\{i:X_i=k,A_i=1,Y_i=1 \} \sim \text{Poi} (np_k \pi_k\mu_{k})  ,
    \]
    \[
    N_{k10} : = \#\{i:X_i=k,A_i=1,Y_i=0 \} \sim \text{Poi} (np_k \pi_k(1-\mu_{k})) ,
    \]
    \[
    N_{k0} := \#\{i:X_i=k,A_i=0 \} \sim \text{Poi} (np_k (1-\pi_k)) .
    \]
    Similarly under the alternative hypothesis $P'$,
    \[
    N_{k11}' : = \#\{i:X_i=k,A_i=1,Y_i=1 \} \sim \text{Poi} (np_k' \pi_k'\mu_{k}') ,
    \]
    \[
    N_{k10}' : = \#\{i:X_i=k,A_i=1,Y_i=0 \} \sim \text{Poi} (np_k' \pi_k'(1-\mu_{k}'))  ,
    \]
    \[
    N_{k0}' := \#\{i:X_i=k,A_i=0 \} \sim \text{Poi} (np_k' (1-\pi_k')) .
    \]
    Denote $\bN_k = (N_{k11}, N_{k10}, N_{k0})$ as the sufficient statistics in $k$-th category and $\bN=(\bN_1,\dots, \bN_{d}, \bN_{d+1})$ as the collection of these sufficient statistics. Define $\bN'$ similarly for the alternative hypothesis. The total variation distance between the marginal distribution of $\bN$ and $\bN'$ (marginalize over the distribution of $\{(p_k, \pi_k, \mu_k ), 1\leq k \leq d\}$ and $\{(p_k', \pi_k', \mu_k' ), 1\leq k \leq d\}$) can be bounded by triangle inequality as (since $\bN_k$ depends only on $(p_k, \pi_k, \mu_k )$ for each $k$ and $\{(p_k, \pi_k, \mu_k ), 1\leq k \leq d\}$ are independent, $(\bN_1,\dots, \bN_d, \bN_{d+1})$ are also independent)
    \[
     \text{TV}(\bN,\bN') \leq\, \sum_{k=1}^d \text{TV}(\bN_k,\bN_k').
    \]
    Note that the marginal distributions of components of $\bN_k$ are not independent since $N_{k11}, N_{k10}, N_{k0}$ all depend on $(p_k, \pi_k, \mu_k )$. Conditioned on $(p_k, \pi_k, \mu_k )$, $(N_{k11}, N_{k10}, N_{k0})$ are conditionally independent with each being a Poisson distribution. By definition of total variation distance,
    \[
    \begin{aligned}
        &\, \text{TV}(\bN_k,\bN_k')\\
        = &\, \frac{1}{2} \sum_{i,j,\ell=0}^{\infty}  \big |  \MP (N_{k11}=i, N_{k10}=j, N_{k0}=\ell)- \MP (N_{k11}'=i, N_{k10}'=j, N_{k0}'=\ell)\big |.
    \end{aligned}
    \]
    Conditioning on $(p_k, \pi_k, \mu_k )$ we have
    \[
    \begin{aligned}
        &\,\MP (N_{k11}=i, N_{k10}=j, N_{k0}=\ell)\\
        = & \, \ME\left [ \exp(-np_k \pi_k \mu_k) \frac{(np_k \pi_k \mu_k)^i}{i!} \exp(-np_k \pi_k(1- \mu_k)) \frac{[np_k \pi_k(1- \mu_k)]^j}{j!} \right.\\
        &\, \left.\exp(-np_k (1-\pi_k)) \frac{[np_k(1-\pi_k)]^{\ell}}{\ell !} \right ] \\
        = &\, \ME \left[ \exp(-n p_k) \frac{(np_k \pi_k \mu_k)^i[np_k \pi_k(1- \mu_k)]^j[np_k(1-\pi_k)]^{\ell}}{i!j!\ell!} \right] \\
        = & \, \ME \left[ \sum_{t=0}^{\infty} \frac{(-np_k)^t(np_k \pi_k \mu_k)^i[np_k \pi_k(1- \mu_k)]^j[np_k(1-\pi_k)]^{\ell}}{i!j!\ell!t!}  \right].
    \end{aligned}
    \]
    Hence the total variation distance can be written as
    \[
    \begin{aligned}
    &\, \text{TV}(\bN_k,\bN_k')  \\
    = &\, \frac{1}{2} \sum_{i,j,\ell} \left | \ME \left[ \sum_{t=0}^{\infty} \frac{(-np_k)^t(np_k \pi_k \mu_k)^i[np_k \pi_k(1- \mu_k)]^j[np_k(1-\pi_k)]^{\ell}}{i!j!\ell!t!}  \right] \right. \\
    &\, \left.- \ME \left[ \sum_{t=0}^{\infty} \frac{(-np_k')^t(np_k' \pi_k' \mu_k')^i[np_k' \pi_k'(1- \mu_k')]^j[np_k'(1-\pi_k')]^{\ell}}{i!j!\ell!t!}  \right] \right |.
    \end{aligned}
    \]
    Note that 
    \[
    (np_k)^t(np_k \pi_k \mu_k)^i[np_k \pi_k(1- \mu_k)]^j[np_k(1-\pi_k)]^{\ell}
    \]
    is a polynomial of $(p_k, p_k \pi_k, p_k\pi_k \mu_k)$ with degree $ t+i+j+\ell$. By moment matching property $1$ in Theorem \ref{lemma:prior-dist}, we have
    \[
    \begin{aligned}
    &\, \text{TV}(\bN_k,\bN_k')  \\
    \leq &\, \frac{1}{2} \left\{ \sum_{i+j+\ell+t > 3K}  \ME \left[ \frac{(np_k)^t(np_k \pi_k \mu_k)^i[np_k \pi_k(1- \mu_k)]^j[np_k(1-\pi_k)]^{\ell}}{i!j!\ell!t!}  \right] \right.  \\
    &\, \left. + \ME \left[  \frac{(np_k')^t(np_k' \pi_k' \mu_k')^i[np_k' \pi_k'(1- \mu_k')]^j[np_k'(1-\pi_k')]^{\ell}}{i!j!\ell!t!}  \right] \right\} .
    \end{aligned}
    \]
    Note that
    \[
    \begin{aligned}
        &\,\frac{(np_k)^t(np_k \pi_k \mu_k)^i[np_k \pi_k(1- \mu_k)]^j[np_k(1-\pi_k)]^{\ell}}{i!j!\ell!t!} \\
        = &\, \exp(2np_k) \exp(-np_k) \frac{(np_k)^t}{t!} \exp(-np_k \pi_k \mu_k) \frac{(np_k \pi_k \mu_k)^i}{i!} \\
        &\, \exp(-np_k \pi_k(1- \mu_k)) \frac{[np_k \pi_k(1- \mu_k)]^j}{j!} \exp(-np_k (1-\pi_k)) \frac{[np_k(1-\pi_k)]^{\ell}}{\ell !}  \\
        = & \, \exp(2np_k)\MP(V_1 = t, V_2 = i, V_3 = j, V_4 = \ell \mid p_k, \pi_k, \mu_k),
    \end{aligned}
    \]
    where conditioning on $(p_k, \pi_k, \mu_k)$, 
    \[
    \begin{aligned}
        V_1 \sim &\, \text{Poi}(n p_k) ,\\
        V_2 \sim &\, \text{Poi}(n p_k \pi_k \mu_k) ,\\
        V_3 \sim &\, \text{Poi}(n p_k \pi_k (1-\mu_k))  ,\\
        V_4 \sim &\, \text{Poi}(n p_k (1-\pi_k)) ,
    \end{aligned}
    \]
    and $(V_1, V_2, V_3, V_4)$ are independent. Further let $V= V_1 + V_2 + V_3 + V_4 \sim \text{Poi}(2n p_k)$ given $(p_k, \pi_k, \mu_k)$. Thus we have
    \[
    \begin{aligned}
    &\, \sum_{i+j+\ell+t > 3K}  \ME \left[ \frac{(np_k)^t(np_k \pi_k \mu_k)^i[np_k \pi_k(1- \mu_k)]^j[np_k(1-\pi_k)]^{\ell}}{i!j!\ell!t!}  \right]\\
    = &\, \ME[\exp(2np_k) \MP(V_1+V_2+V_3+V_4 > 3K \mid p_k, \pi_k,\mu_k)] \\
    = &\, \ME[\exp(2np_k) \MP(V > 3K \mid p_k, \pi_k,\mu_k)] \\
    = &\, \ME\left[\exp(2np_k) \sum_{\ell > 3K} \exp(-2np_k) \frac{(2np_k)^{\ell}}{\ell!}\right] \\
    = &\, \ME \left[\sum_{\ell > 3K}  \frac{(2np_k)^{\ell}}{\ell!} \right] \\
    \leq &\, \sum_{\ell > 3K}  \frac{(2c_1 \log n)^{\ell}}{\ell!}\\
    = & \, \exp(2c_1 \log n) \sum_{\ell > 3K} \exp(-2c_1 \log n) \frac{(2c_1 \log n)^{\ell}}{\ell!} \\
    = &\, n^{2c_1} \MP(W>3K),
    \end{aligned}
    \]
    where $W \sim \text{Poi}(2c_1 \log n)$. Apply the following Chernoff bound: For $L>eM$
    \begin{equation}\label{eq:chernoff-poi}
    \MP(\text{Poi}(M)>L) \leq \exp(-M) \left(\frac{eM}{L} \right)^L.        
    \end{equation}
    When $3K > 2e c_1 \log n$, we have
    \[
    \MP(W>3K) \leq \exp(-2c_1 \log n) \left(\frac{2ec_1}{3c_2} \right)^{3c_2 \log n},
    \]
    \[
    n^{2c_1} \MP(W>3K) \leq \left(\frac{2ec_1}{3c_2} \right)^{3c_2 \log n} = n^{3c_2 \log \left(\frac{2ec_1}{3c_2} \right)} \leq n^{-2} \leq \frac{c_3 \log n}{nd}
    \]
    as long as we choose constant $c_1, c_2$ satisfying $3c_2 \log \left(\frac{3c_2}{2ec_1} \right) \geq 2$ and $c_3$ is the constant in Lemma \ref{lemma:prior-dist}. Thus the total variation distance is bounded as (similar inequality holds under the alternative hypothesis) 
    \[
    \text{TV}(\bN_k,\bN_k') \leq \frac{c_3 \log n}{nd}
    \]
    \[
    \text{TV}(\bN,\bN') \leq \frac{c_3 \log n}{n}.
    \]
    The functional separation is
    \[
    \psi_1(P) - \psi_1(P') = \sum_{k=1}^d (p_k \mu_k - p_k' \mu_k'),
    \]
    \[
    \left |\ME\left[\psi_1(P) - \psi_1(P') \right] \right | = d | \ME[p\mu - p' \mu'] | \geq \frac{c_4 d}{n \log n}:= q.
    \]
    Consider the following events:
    \[
    E=\left\{\bigg|  \sum_{k=1}^d p_k -da_0\bigg|\leq \delta, \bigg | \sum_{k=1}^d p_k \mu_k  - d\ME_{\mu_0}[p \mu]\bigg|\leq q/4 \right\},
    \]
    \[
    E'=\left\{\bigg|  \sum_{k=1}^d p_k' -da_0\bigg|\leq \delta, \bigg | \sum_{k=1}^d p_k' \mu_k'  -d \ME_{\mu_1}[p' \mu']\bigg|\leq q/4 \right\},
    \]
    By Chebyshev's inequality, we have
    \[
    \begin{aligned}
        \MP(E^c) \leq &\, \MP\left(\bigg|  \sum_{k=1}^d p_k -da_0\bigg|> \delta\right) \\
        &\, + \MP \left(\bigg | \sum_{k=1}^d p_k \mu_k  - d\ME_{\mu_0}[p \mu]\bigg| > q/4 \right) \\
        \leq &\, \frac{d \operatorname{Var}_{\mu_0}(p)}{ \delta^2} + \frac{16 d \operatorname{Var}_{\mu_0}(p \mu)}{q^2} \\
        \leq &\, \frac{c_1^2 d (\log n)^2}{ \delta^2 n^2} + \frac{16 c_1^2 d (\log n)^2}{n^2 q^2} \\
        = &\, \frac{1}{16} + \frac{16c_1^2 (\log n)^4}{ c_4^2 d},
    \end{aligned}
    \]
    where the third inequality follows from the bound $|p\mu| \leq |p|\leq c_1 \log n/n $ and the last equation follows by setting $\delta = \frac{4 c_1 \sqrt{d} \log n }{n} $. Similarly, we have
    \[
    \MP(E'^c) \leq \frac{1}{16} + \frac{16c_1^2 (\log n)^4}{ c_4^2 d}. 
    \]
    We put the following prior distributions induced by $\{(p_k, \pi_k, \mu_k), 1 \leq k \leq d\}$ and $\{(p_k', \pi_k', \mu_k'), 1 \leq k \leq d\}$ on $P$ and $P'$, respectively:
    \[
    \pi \stackrel{d}{=} P \mid E, \pi' \stackrel{d}{=} P' \mid E'.
    \]
    Note that under $\pi, \pi'$,
    \[
    |\psi_1(P)-\psi_1(P')|\geq q/2.
    \]
    By triangle inequality, the total variation distance of the sufficient counting statistics $\bN$ and $\bN'$ under two priors is bounded by
    \[
    \begin{aligned}
         \text{TV}\left( {\bN \mid E}, {\bN' \mid E'} \right) \leq &\, \text{TV}\left( {\bN \mid E}, {\bN} \right) + \text{TV}({\bN},{\bN'} ) + \text{TV}\left( {\bN' \mid E'}, {\bN'} \right)\\
         \leq &\, \MP(E^c) + \MP(E'^c) + \text{TV}({\bN},{\bN'} )\\
         \leq &\, \frac{1}{8} + \frac{32c_1^2 (\log n)^4}{ c_4^2 d} + \frac{c_3 \log n}{n} .
    \end{aligned}
    \]
    By method of fuzzy hypotheses (Section 2.7.4 of \cite{Tsybakov2009}) we conclude
    \[
    \begin{aligned}
    &\, \tilde{R}^*(d,n,\delta) \\
    \geq &\, \frac{q^2}{32} \left( 1 -  \text{TV}\left( {\bN \mid E}, {\bN' \mid E'} \right) \right).\\
    \geq &\, \frac{q^2}{32} \left( \frac{7}{8} - \frac{32c_1^2 (\log n)^4}{ c_4^2 d} - \frac{c_3 \log n}{n} \right).
    \end{aligned}
    \]
    Hence in the regime $d \gtrsim (\log n)^4$, we have 
    \[
    \tilde{R}^*(d,n,\delta) \gtrsim \frac{d^2}{(n \log n)^2}.
    \]
    By Lemma \ref{lemma:poission-minimax} we conclude
    \[
    \begin{aligned}
    {R}^*(d,n) \gtrsim &\,  \frac{d^2}{(n \log n)^2} - \exp(-n/50)- \frac{d (\log n)^2}{n^2}\\
    \gtrsim &\, \frac{d^2}{(n \log n)^2}.
    \end{aligned}
    \]
    The overall requirements on the constants are
    \[
    \frac{c c_1 c_3}{c_2^2} \leq 1, \, 3c_2 \log \left(\frac{3c_2}{2ec_1} \right) \geq 2.
    \]
    Clearly for $c, c_3 >0$, one can choose $c_1$ sufficiently small and $c_2$ sufficiently large to satisfy these conditions.
    
    The lower bound $1/n$ can be proved using a two-point method. 
Without loss of generality, assume that $n \geq 8$. Under the null $P$ and alternative hypothesis $P'$, set
\[
p_k = \frac{1}{d}, \pi_k = \frac{1}{2}, \mu_{0k}=0,  \mu_{1k}(P) = \frac{1}{2}, \mu_{1k}(P') = \frac{1}{2} + \delta,
\]
with $\delta = 1/\sqrt{n}$, i.e., the covariate distribution and propensity score are the same. By Le Cam's two-point method we have
\begin{equation}\label{eq:two-point}
\inf_{\hat{\psi}_1} \sup_{\MP \in \mathcal{D} (\epsilon)}  \ME_{\MP}\left[\left(\hat{\psi}_1 - \psi_1\right)^2\right] \geq \frac{1}{4}(\psi_1(P)-\psi_1(P'))^2 \exp (-n D(P \| P')).
\end{equation}
Note that the functional separation is
\[
|\psi_1(P)-\psi_1(P')| = \delta = \frac{1}{\sqrt{n}}.
\]
Under the null hypothesis $P$, we have
\[
P(X=k, A=a, Y=y) = \left\{
\begin{array}{rl}
\frac{1}{4d} & \text{if } a=1, y=1,\\
\frac{1}{4d} & \text{if } a=1, y=0,\\
\frac{1}{2d} & \text{if } a=0, y=0.
\end{array} \right.
\]
Under the alternative hypothesis $P'$, we have
\[
P'(X=k, A=a, Y=y) = \left\{
\begin{array}{rl}
\frac{1}{2d}\left( \frac{1}{2}+\delta \right) & \text{if } a=1, y=1,\\
\frac{1}{2d} \left( \frac{1}{2}-\delta \right) & \text{if } a=1, y=0,\\
\frac{1}{2d} & \text{if } a=0, y=0.
\end{array} \right.
\]
The K-L divergence between $P$ and $P'$ is
\[
D(P \| P') = \sum_{k=1}^d \left[- \frac{1}{4d} \log(1+2\delta) - \frac{1}{4d} \log(1-2\delta) \right] = -\frac{1}{4} \log (1-4\delta^2).
\]
Using the inequality 
\[
\log(1-x) \geq -2x, x \in [0,1/2],
\]
we have
\[
D(P \| P') =-\frac{1}{4} \log \left(1-\frac{4}{n} \right) \leq \frac{2}{n}.
\]
Plug in the functional separation and bound on K-L divergence into \eqref{eq:two-point}, we conclude
\[
\inf_{\hat{\psi}_1} \sup_{\MP \in \mathcal{D} (\epsilon)}  \ME_{\MP}\left[\left(\hat{\psi}_1 - \psi_1\right)^2\right] \gtrsim \frac{1}{n}.
\]
    
\end{proof}

\subsection{Proof of Lemma \ref{lemma:nocollision}}
\begin{proof}
By the definition of $\hat{t}_k = I(\hat{p}_k >0, 0 < \hat{\pi}_k < 1)$ we have
\[
\begin{aligned}
\mathbb{P}\left(\sum_{j=1}^d \widehat{t}_j=0\right) & =\mathbb{P}\left(\widehat{w}_k \in\left\{0, \widehat{p}_k\right\} \forall k \text { with } \widehat{p}_k>0\right) \\
& =\mathbb{E}\left\{\mathbb{P}\left(\widehat{w}_k \in\left\{0, \widehat{p}_k\right\} \forall k \text { with } \widehat{p}_k>0 \mid \bX^n\right)\right\} \\
& =\mathbb{E}\left\{\prod_{k: \hat{p}_k > 0} \mathbb{P}\left(\widehat{w}_k \in\left\{0, \widehat{p}_k\right\} \mid \bX^n\right)\right\} \\
& =\mathbb{E}\left(\prod_{k=1}^d\left[\left\{\left(1-\pi_k\right)^{n \widehat{p}_k}+\pi_k^{n \widehat{p}_k}\right\} I\left(\widehat{p}_k>0\right)+I\left(\widehat{p}_k=0\right)\right]\right) \\
& =\mathbb{E}\left[\prod_{k=1}^d\left\{\left(1-\pi_k\right)^{n \widehat{p}_k}+\pi_k^{n \widehat{p}_k}-I\left(\widehat{p}_k=0\right)\right\}\right]
\end{aligned}
\]
The proof relies on the poissonization technique to bound the above expectation of the product. Poissonization allows us to replace $\left(n \widehat{p}_1, \ldots, n \widehat{p}_d\right) \sim \operatorname{Multinomial}\left(n, p_1, \ldots, p_d\right)$ with $n\hat{p}_k \sim \text{Poisson}(n p_k)$ and $n\hat{p}_1, \dots, n\hat{p}_d $ are independent. The following lemma connects the expectation in the multinomial case with that in the independent Poisson case.

\begin{lemma}[Theorem 5.10 in \cite{mitzenmacher2017probability}]\label{lemma:poissonization}

Let $\bX^n \in \mathbb{R}^d \sim \text{Multinomial}\left(n, p_1, \ldots, p_d\right)$, $ \bY^n \in \mathbb{R}^d$ and $Y_{i}^n \sim \text{Poisson}(np_i)$ and $Y_{1}^n, \dots, Y_{d}^n$ are independent. Consider a non-negative function $f(x_1, \dots, x_d)$, if $\ME[f(X_1^n,\dots, X_d^n)]$ is monotonely non-increasing with $n$, then
\[
\ME[f(X_1^n,\dots, X_d^n)] \leq 2 \ME[f(Y_1^n,\dots, Y_d^n)].
\]
\end{lemma}
The proof is left as an exercise in \cite{mitzenmacher2017probability} and we include it in Appendix \ref{sec:proof-poissonization}. Let 
\[
\begin{aligned}
    a_n =&\,  \mathbb{E}\left[\prod_{k=1}^d\left\{\left(1-\pi_k\right)^{n \widehat{p}_k}+\pi_k^{n \widehat{p}_k}-I\left(\widehat{p}_k=0\right)\right\}\right] \\
    =&\, \mathbb{E}\left[\prod_{k=1}^d\left\{\left(1-\pi_k\right)^{\sum_{i=1}^n I(X_i=k)}+\pi_k^{\sum_{i=1}^n I(X_i=k)}-I\left(\sum_{i=1}^n I(X_i=k) = 0\right)\right\}\right].
\end{aligned}
\]

First we verify the monotonicity of $a_n$. We claim
\[
\begin{aligned}
    &\, \left(1-\pi_k\right)^{\sum_{i=1}^{n+1} I(X_i=k)}+\pi_k^{\sum_{i=1}^{n+1} I(X_i=k)}-I\left(\sum_{i=1}^{n+1} I(X_i=k) = 0\right) \\
    \leq &\, \left(1-\pi_k\right)^{\sum_{i=1}^{n} I(X_i=k)}+\pi_k^{\sum_{i=1}^{n} I(X_i=k)}-I\left(\sum_{i=1}^{n} I(X_i=k) = 0\right).
\end{aligned}
\]
On the event $\left\{\sum_{i=1}^{n+1 } I(X_i=k) = 0\right\}$ then $\sum_{i=1}^{n } I(X_i=k) = 0$  and
\[
\left(1-\pi_k\right)^{\sum_{i=1}^{n+1} I(X_i=k)}+\pi_k^{\sum_{i=1}^{n+1} I(X_i=k)}-I\left(\sum_{i=1}^{n+1} I(X_i=k) = 0\right) = 1
\]
\[
\left(1-\pi_k\right)^{\sum_{i=1}^{n} I(X_i=k)}+\pi_k^{\sum_{i=1}^{n} I(X_i=k)}-I\left(\sum_{i=1}^{n} I(X_i=k) = 0\right) = 1.
\]
On the event $\left\{\sum_{i=1}^{n} I(X_i=k) = 0, X_{n+1} = k\right\}$,
\[
\left(1-\pi_k\right)^{\sum_{i=1}^{n+1} I(X_i=k)}+\pi_k^{\sum_{i=1}^{n+1} I(X_i=k)}-I\left(\sum_{i=1}^{n+1} I(X_i=k) = 0\right) = 1-\pi_k + \pi_k = 1
\]
\[
\left(1-\pi_k\right)^{\sum_{i=1}^{n} I(X_i=k)}+\pi_k^{\sum_{i=1}^{n} I(X_i=k)}-I\left(\sum_{i=1}^{n} I(X_i=k) = 0\right) = 1.
\]
On the event $\left\{\sum_{i=1}^n I(X_i=k) > 0 \right\}$,
\[
\begin{aligned}
    &\, \left(1-\pi_k\right)^{\sum_{i=1}^{n+1} I(X_i=k)}+\pi_k^{\sum_{i=1}^{n+1} I(X_i=k)}-I\left(\sum_{i=1}^{n+1} I(X_i=k) = 0\right) \\
    = &\, \left(1-\pi_k\right)^{\sum_{i=1}^{n+1} I(X_i=k)}+\pi_k^{\sum_{i=1}^{n+1} I(X_i=k)} \\
    \leq &\, \left(1-\pi_k\right)^{\sum_{i=1}^{n}  I(X_i=k)}+\pi_k^{\sum_{i=1}^{n} I(X_i=k)} \\
    = &\, \left(1-\pi_k\right)^{\sum_{i=1}^{n} I(X_i=k)}+\pi_k^{\sum_{i=1}^{n} I(X_i=k)}-I\left(\sum_{i=1}^{n} I(X_i=k) = 0\right).
\end{aligned}
\]
Hence the claim is verified and $a_n$ is non-increasing. Apply Lemma \ref{lemma:poissonization} we can now assume $n\hat{p}_k \sim \text{Poisson}(np_k)$ and $n\hat{p}_1,\dots, n\hat{p}_d$ are independent (with an additional factor 2)
\[
\begin{aligned}
    &\, \mathbb{P}\left(\sum_{j=1}^d \widehat{t}_j=0\right) \\
    \leq &\, 2 \mathbb{E}\left[\prod_{k=1}^d\left\{\left(1-\pi_k\right)^{n \widehat{p}_k}+\pi_k^{n \widehat{p}_k}-I\left(\widehat{p}_k=0\right)\right\}\right] \\
    =&\, 2 \prod_{k=1}^d \ME\left[\left(1-\pi_k\right)^{n \widehat{p}_k}+\pi_k^{n \widehat{p}_k}-I\left(\widehat{p}_k=0\right)\right] \\
    =&\, 2 \prod_{k=1}^d \left[ \exp(-n\pi_k p_k) + \exp(-n(1-\pi_k)p_k) -\exp(-np_k)  \right] \\
    \leq &\, 2 \prod_{k=1}^d \left[ \exp(-\epsilon n p_k) + \exp(-(1-\epsilon)n p_k) -\exp(-np_k)  \right].
\end{aligned}
\]
The first equation follows from independence and second one follows from probability generating function of Poisson distribution and the last inequality follows since $f_3(x) = \exp(-cx) + \exp(-c(1-x))$ (for constant $c>0$) is decreasing on $[0,1/2]$ and increasing on $[1/2, 1]$. For simplicity define 
\[
Z_k =  \exp(-\epsilon n p_k) + \exp(-(1-\epsilon)n p_k) -\exp(-np_k) .
\]
1. In the first case $\epsilon n p_k \geq 2\log 2$ we have (recall $0 < \epsilon < 1/2$)
\[
Z_k \leq 2\exp(-\epsilon n p_k) \leq \exp\left( -\frac{\epsilon n p_k}{2} \right).
\]
2. In the second case $n p_k \leq \epsilon$, we use the following inequalities for $t \geq 0$:
\[
\begin{aligned}
& \exp (-t) \leq 1-t+t^2 / 2, \\
& \exp (-t) \geq 1-t+t^2 / 2-t^3 / 6
\end{aligned}
\]
and obtain
\[
\begin{aligned}
    Z_k \leq &\, 1-\epsilon n p_k + \frac{\epsilon^2 n^2 p_k^2}{2} + 1-(1-\epsilon)n p_k + \frac{(1-\epsilon)^2 n^2 p_k^2}{2} \\
    &\, - \left( 1-n p_k + \frac{n^2 p_k^2}{2} - \frac{n^3 p_k^3}{6} \right) \\
    = &\, 1-\epsilon(1-\epsilon)n^2 p_k^2 + \frac{n^3 p_k^3}{6} \\
    \leq & \, 1- \left( \frac{5}{6}\epsilon -\epsilon^2 \right)n^2 p_k^2\\
    \leq & \, 1-\frac{\epsilon n^2 p_k^2}{3} \\
    \leq &\, \exp \left( -\frac{\epsilon n^2 p_k^2}{3} \right).
\end{aligned}
\]
3. In the last case $\epsilon < n p_k < \frac{2 \log 2}{\epsilon}$, since $f_4(x) = \exp(-\epsilon x) + \exp(-(1-\epsilon)x) -\exp(-x)$ is monotonely non-increasing on $[0,+\infty]$ (one can check this by taking the first-order derivative easily), we have
\[
Z_k \leq \exp(-\epsilon^2) + \exp(-\epsilon(1-\epsilon)) - \exp(-\epsilon) = \exp (-C_1(\epsilon))
\]
where $C_1(\epsilon) = -\log(\exp(-\epsilon^2) + \exp(-\epsilon(1-\epsilon)) - \exp(-\epsilon)) \in(0, \epsilon^2)$. Note that 
\[
1 > \frac{\epsilon n p_k}{2 \log 2}.
\]
We have
\[
\exp (-C_1(\epsilon)) \leq \exp (-C_2(\epsilon) n p_k)
\]
where 
\[
C_2(\epsilon) = \frac{C_1(\epsilon)\epsilon}{2 \log 2} < \frac{\epsilon ^3}{2 \log 2} < \frac{\epsilon }{2} 
\]
Hence we can combine the first case with the third case as
\[
Z_k \leq \exp (-C_2(\epsilon) n p_k).
\]
Let $I_1$ include the index $k \in [d]$ in the first and third case, $I_2$ include the index $k$ in the second case. Denote $S_1 = \sum_{i\in I_1} p_i, T_2 = \sum_{i \in I_2}p_i^2$. We thus have
\[
\prod_{k=1}^d Z_k \leq \exp(-C_2(\epsilon)nS_1) \exp \left(-\frac{\epsilon n^2 T_2}{3} \right).
\]
In the case $n \geq d$ we have
\[
1-S_1 = \sum_{i \in I_2} p_i \leq \frac{\epsilon |I_2|}{n} \leq \frac{\epsilon d}{n} \leq \epsilon,
\]
i.e. $S_1 > 1-\epsilon$ and we have
\[
\prod_{k=1}^d Z_k \leq \exp(-C_2(\epsilon)(1-\epsilon)n) \leq \exp\left(- \frac{C_2(\epsilon)n}{2 } \right).
\]
In the case $n < d$, if $S_1 \geq 1/2$ then we also have
\[
\prod_{k=1}^d Z_k \leq \exp\left(- \frac{C_2(\epsilon)n}{2 } \right).
\]
In the case $S_1 < 1/2$, then by Cauchy-Schwarz inequality we have
\[
\frac{1}{4} \leq (1-S_1)^2 = \left(\sum_{i \in I_2} p_i \right)^2 \leq |I_2|T_2 \leq d T_2,
\]
i.e. $T_2 \geq \frac{1}{4d}$, this further implies
\[
\prod_{k=1}^d Z_k \leq \exp \left(- \frac{\epsilon n^2}{12d} \right).
\]
So we conclude that
\[
\prod_{k=1}^d Z_k \leq \max \left\{ \exp\left(- \frac{C_2(\epsilon)n}{2 } \right),  \exp \left(- \frac{\epsilon n^2}{12d} \right) \right \} \leq \exp\left(-C (\epsilon) \frac{n^2}{n \vee d} \right)
\]
where 
\[
C(\epsilon) = \min \left( \frac{C_2(\epsilon)}{2}, \frac{\epsilon}{12} \right).
\]
\end{proof}

\subsection{Proof of Theorem \ref{thm:homogeneity}}

\begin{proof}
    We first bound the bias. By conditioning on $\bX^n,\bA^n$ and noting $\ME[\hat{\mu}_{1k} \mid \bX^n, \bA^n] = \mu_{1k} I(\hat{p}_k \hat{\pi}_k >0)$, $\ME[\hat{\mu}_{0k} \mid \bX^n, \bA^n] = \mu_{0k} I(\hat{p}_k (1-\hat{\pi}_k) >0)$ we have
    \[
    \ME\left[\hat{\tau}  \mid\bX^n, \bA^n\right] = \frac{\sum_{k=1}^d \hat{t}_k \tau_k}{\sum_{k=1}^d \hat{t}_k} I\left(\sum_{j=1}^d \hat{t}_j>0\right).
    \]
    Note that 
    \[
    \psi = \psi I\left(\sum_{j=1}^d \hat{t}_j>0\right) + \psi I\left(\sum_{j=1}^d \hat{t}_j=0\right) =  \frac{\sum_{k=1}^d \hat{t}_k \psi}{\sum_{k=1}^d \hat{t}_k} I\left(\sum_{j=1}^d \hat{t}_j>0\right) + \psi I\left(\sum_{j=1}^d \hat{t}_j=0\right).
    \]
    Hence the bias can be expressed as 
    \[
    \begin{aligned}
        &\,\ME\left[\hat{\tau}-\psi\right] \\
        =&\, \ME \left[ \frac{\sum_{k=1}^d \hat{t}_k (\tau_k-\psi)}{\sum_{k=1}^d \hat{t}_k} I\left(\sum_{j=1}^d \hat{t}_j>0\right) \right] - \psi \MP\left(\sum_{j=1}^d \hat{t}_j=0\right)
    \end{aligned}
    \]
    Using $\sigma_n$ in \eqref{eq:homogeneity} to parameterize the effect homogeneity and Lemma \ref{lemma:nocollision}, we have 
    \[
    \begin{aligned}
        &\,\left |\ME\left[\hat{\tau}-\psi\right] \right| \\
        \leq &\, \ME \left[ \frac{\sum_{k=1}^d \hat{t}_k |\tau_k-\psi|}{\sum_{k=1}^d \hat{t}_k} I\left(\sum_{j=1}^d \hat{t}_j>0\right) \right]+ |\psi| \MP\left(\sum_{j=1}^d \hat{t}_j=0\right) \\
        \leq &\, \sigma_n + 2 \exp\left(-C (\epsilon) \frac{n^2}{n \vee d} \right).
    \end{aligned}
    \]
    We then consider the variance. By the property of conditional variance, we have
    \begin{equation}\label{eq:var-homogeneity}
        \operatorname{Var}(\hat{\tau}) = \operatorname{Var}\left(\ME\left[\hat{\tau} \mid\bX^n, \bA^n\right] \right) + \ME \left[ \operatorname{Var}\left(\hat{\tau} \mid \bX^n, \bA^n\right) \right].
    \end{equation}
    Rewrite 
    \[
    \ME\left[\hat{\tau} \mid\bX^n, \bA^n\right] = \frac{\sum_{k=1}^d \hat{t}_k \tau_k}{\sum_{k=1}^d \hat{t}_k} I\left(\sum_{j=1}^d \hat{t}_j>0\right) = \frac{\sum_{k=1}^d \hat{t}_k (\tau_k-\psi)}{\sum_{k=1}^d \hat{t}_k} I\left(\sum_{j=1}^d \hat{t}_j>0\right) + \psi I\left(\sum_{j=1}^d \hat{t}_j>0\right).
    \]
    For a bounded random variable $X$ satisfying $m \leq X \leq M$ we have $\operatorname{Var}(X) \leq (M-m)^2/4$. Note that
    \[
    \left| \frac{\sum_{k=1}^d \hat{t}_k (\tau_k-\psi)}{\sum_{k=1}^d \hat{t}_k} I\left(\sum_{j=1}^d \hat{t}_j>0\right) \right| \leq \sigma_n
    \]
    and we have
    \[
    \operatorname{Var}\left(\frac{\sum_{k=1}^d \hat{t}_k (\tau_k-\psi)}{\sum_{k=1}^d \hat{t}_k} I\left(\sum_{j=1}^d \hat{t}_j>0\right) \right) \leq \sigma_n^2. 
    \]
    Further notice
    \[
    \operatorname{Var}\left( \psi I\left(\sum_{j=1}^d \hat{t}_j>0\right) \right) \leq \psi^2 \MP \left(\sum_{j=1}^d \hat{t}_j=0 \right) \leq  2 \exp\left(-C (\epsilon) \frac{n^2}{n \vee d} \right).
    \]
    The covariance can be expressed as
    \[
    \begin{aligned}
        &\, \operatorname{Cov} \left( \frac{\sum_{k=1}^d \hat{t}_k (\tau_k-\psi)}{\sum_{k=1}^d \hat{t}_k} I\left(\sum_{j=1}^d \hat{t}_j>0\right), \psi  I\left(\sum_{j=1}^d \hat{t}_j>0\right)\right) \\
        = &\, \psi \ME \left[\frac{\sum_{k=1}^d \hat{t}_k (\tau_k-\psi)}{\sum_{k=1}^d \hat{t}_k} I\left(\sum_{j=1}^d \hat{t}_j>0\right) \right] -\psi \ME \left[\frac{\sum_{k=1}^d \hat{t}_k (\tau_k-\psi)}{\sum_{k=1}^d \hat{t}_k} I\left(\sum_{j=1}^d \hat{t}_j>0\right) \right] \MP \left(\sum_{j=1}^d \hat{t}_j>0 \right) \\
        =&\, \psi \ME \left[\frac{\sum_{k=1}^d \hat{t}_k (\tau_k-\psi)}{\sum_{k=1}^d \hat{t}_k} I\left(\sum_{j=1}^d \hat{t}_j>0\right) \right] \MP \left(\sum_{j=1}^d \hat{t}_j=0 \right)
    \end{aligned}
    \]
    and hence we can bound the covariance as follows:
    \[
    \begin{aligned}
        &\, \left|\operatorname{Cov} \left( \frac{\sum_{k=1}^d \hat{t}_k (\tau_k-\psi)}{\sum_{k=1}^d \hat{t}_k} I\left(\sum_{j=1}^d \hat{t}_j>0\right), \psi  I\left(\sum_{j=1}^d \hat{t}_j>0\right)\right) \right|\\
        \leq &\, 2\sigma_n \exp\left(-C (\epsilon) \frac{n^2}{n \vee d} \right).
    \end{aligned}
    \]
    Combining these three terms above we have
    \begin{equation}\label{eq:var-cond-exp}
    \operatorname{Var}\left(\ME\left[\hat{\tau}\mid\bX^n, \bA^n\right] \right) \leq \sigma_n^2 + 2 \exp\left(-C (\epsilon) \frac{n^2}{n \vee d} \right) + 4\sigma_n \exp\left(-C (\epsilon) \frac{n^2}{n \vee d} \right).
    \end{equation}
    To complete the proof we need to bound 
    \[
    \ME \left[ \operatorname{Var}\left(\hat{\tau} \mid \bX^n, \bA^n\right) \right].
    \]
    Note that conditioning on $\bX^n,\bA^n$ the estimated $ \{\hat{\tau}_k, 1\leq k \leq d \}$ are independent with 
    \[
    \operatorname{Var}(\hat{\tau}_k\mid\bX^n, \bA^n) = \frac{\mu_{1k}(1-\mu_{1k})}{n\hat{p}_k \hat{\pi}_k} I(\hat{p}_k \hat{\pi}_k>0) +  \frac{\mu_{0k}(1-\mu_{0k})}{n\hat{p}_k (1-\hat{\pi}_k)}I\left(\hat{p}_k (1-\hat{\pi}_k)>0\right).
    \]
    We have
    \[
    \begin{aligned}
        \operatorname{Var}(\hat{\tau}\mid\bX^n, \bA^n) = &\, \frac{I\left(\sum_{j=1}^d \hat{t}_j > 0\right)}{\left(\sum_{j=1}^d \hat{t}_j\right)^2} \sum_{k=1}^d \hat{p}_k^2 I(0 < \hat{\pi}_k < 1) \left( \frac{\mu_{1k}(1-\mu_{1k})}{n\hat{p}_k \hat{\pi}_k} +\frac{\mu_{0k}(1-\mu_{0k})}{n\hat{p}_k (1-\hat{\pi}_k)} \right)\\
        \leq &\, \frac{1}{4}\frac{I\left(\sum_{j=1}^d \hat{t}_j > 0\right)}{\left(\sum_{j=1}^d \hat{t}_j\right)^2} \sum_{k=1}^d \hat{p}_k^2 I(0 < \hat{\pi}_k < 1) \left( \frac{1}{n\hat{p}_k \hat{\pi}_k} +\frac{1}{n\hat{p}_k (1-\hat{\pi}_k)} \right)
    \end{aligned}
    \]
    Denote $\bI = (I_1,\dots, I_d)$ with $I_k = I(0 < \hat{\pi}_k < 1)$ and
    further note that
    \[
    \begin{aligned}
        &\,\frac{1}{4} \ME \left[ \frac{I\left(\sum_{j=1}^d \hat{t}_j > 0\right)}{\left(\sum_{j=1}^d \hat{t}_j\right)^2} \sum_{k=1}^d   \frac{\hat{p}_k^2 I(0 < \hat{\pi}_k < 1)}{n\hat{p}_k \hat{\pi}_k}  \right]\\
        = &\, \frac{1}{4} \ME \left[ \frac{I\left(\sum_{j=1}^d \hat{t}_j > 0\right)}{\left(\sum_{j=1}^d \hat{t}_j\right)^2} \sum_{k=1}^d  \hat{p}_k^2 I(0 < \hat{\pi}_k < 1) \ME\left(\frac{1}{n\hat{p}_k \hat{\pi}_k} \mid \bX^n, \bI \right)  \right] \\
        = &\, \frac{1}{4} \ME \left[ \frac{I\left(\sum_{j=1}^d \hat{t}_j > 0\right)}{\left(\sum_{j=1}^d \hat{t}_j\right)^2} \sum_{k=1}^d  \hat{p}_k^2 I(0 < \hat{\pi}_k < 1) \ME\left(\frac{1}{n\hat{p}_k \hat{\pi}_k} \mid \bX^n, I_k \right)  \right],
    \end{aligned}
    \]
    where the last equation follows from the independence of treatment assignment within each category after covariates $\bX^n$ are sampled. It is easy to see from Lemma \ref{lemma:bino-inverse} that for a Binomial random variable $V \sim \text{B}(n,p)$ we have
    \[
    \ME\left[\frac{1}{V} \mid 0 < V < n\right] \leq \frac{2}{(n+1)p\left[1-p^n-(1-p)^n\right]}.
    \]
    Thus we have (assume $n \geq 2$)
    \[
    \begin{aligned}
        &\,\hat{p}_k^2 I(0 < \hat{\pi}_k < 1) \ME\left(\frac{1}{n\hat{p}_k \hat{\pi}_k} \mid \bX^n, I_k \right)\\
        \leq &\, \frac{2\hat{p}_k^2 I(0 < \hat{\pi}_k < 1)}{(n\hat{p}_k+1)\pi_k\left[1-\pi_k^n-(1-\pi_k)^n\right]}\\
        \leq &\, \frac{2\hat{p}_k^2 I(0 < \hat{\pi}_k < 1)}{(n\hat{p}_k+1)\epsilon\left[1-\epsilon^n-(1-\epsilon)^n\right]}\\
        \leq &\, \frac{2\hat{p}_k^2 I(0 < \hat{\pi}_k < 1)}{(n\hat{p}_k+1)\epsilon\left[1-\epsilon^2-(1-\epsilon)^2\right]}\\
        = &\, \frac{\hat{p}_k^2 I(0 < \hat{\pi}_k < 1)}{(n\hat{p}_k+1)\epsilon^2(1-\epsilon)} \\
        \leq &\, \frac{\hat{p}_k I(0 < \hat{\pi}_k < 1)}{\epsilon^2(1-\epsilon)n}.
    \end{aligned}
    \]
    Sum these terms up, we have
    \[
    \frac{1}{4} \ME \left[ \frac{I\left(\sum_{j=1}^d \hat{t}_j > 0\right)}{\left(\sum_{j=1}^d \hat{t}_j\right)^2} \sum_{k=1}^d  \hat{p}_k^2 I(0 < \hat{\pi}_k < 1) \ME\left(\frac{1}{n\hat{p}_k \hat{\pi}_k} \mid \bX^n, I_k \right)  \right] \leq \frac{1}{4\epsilon^2(1-\epsilon)} \ME \left[ \frac{I\left(\sum_{j=1}^d \hat{t}_j > 0\right)}{n\sum_{j=1}^d \hat{t}_j} \right]
    \]
    Similarly, we can show
    \[
    \frac{1}{4} \ME \left[ \frac{I\left(\sum_{j=1}^d \hat{t}_j > 0\right)}{\left(\sum_{j=1}^d \hat{t}_j\right)^2} \sum_{k=1}^d   \frac{\hat{p}_k^2 I(0 < \hat{\pi}_k < 1)}{n\hat{p}_k (1-\hat{\pi}_k)}  \right]\leq \frac{1}{4\epsilon(1-\epsilon)^2} \ME \left[ \frac{I\left(\sum_{j=1}^d \hat{t}_j > 0\right)}{n\sum_{j=1}^d \hat{t}_j} \right]
    \]
    So we have the following bound on the expected conditional variance:
    \[
    \ME\left[\operatorname{Var}(\hat{\tau}\mid\bX^n, \bA^n)\right]\leq \frac{1}{4\epsilon^2(1-\epsilon)^2} \ME \left[ \frac{I\left(\sum_{j=1}^d \hat{t}_j > 0\right)}{n\sum_{j=1}^d \hat{t}_j} \right]  
    \]
    Define $h(\hat{\bp}, \hat{\pi})$
    \[
    h(\hat{\bp}, \hat{\boldsymbol{\pi}}) =
    \begin{cases}
    1 &\quad\text{if} \sum_{j=1}^d \hat{t}_j =0\\
    \frac{1}{\sum_j n\hat{p}_j I(0 < \hat{\pi}_j < 1)} &\quad\text{if} \sum_{j=1}^d \hat{t}_j >0 . \\
    \end{cases}
    \]
    Note that $\sum_j n\hat{p}_j I(0 < \hat{\pi}_j < 1)$ is the number of subjects in categories with both treated and untreated units, thus it is an integer and will not decrease as we collect more samples. As a result,
    $\ME[h(\hat{\bp}, \hat{\boldsymbol{\pi}})]$ is non-increasing in $n$ and by Lemma \ref{lemma:poissonization} we have
    \begin{equation}\label{eq:bound-the-inverse}
         \begin{aligned}
        &\, \ME \left[ \frac{I\left(\sum_{j=1}^d \hat{t}_j > 0\right)}{n\sum_{j=1}^d \hat{t}_j} \right]\\
        \leq &\, \ME[h(\hat{\bp}, \hat{\boldsymbol{\pi}})]\\
        \leq &\,2\ME_{\text{Poi}}[h(\hat{\bp}, \hat{\boldsymbol{\pi}})] \\
        =&\, 2\MP\left(\sum_j n\hat{p}_j I(0 < \hat{\pi}_j < 1) = 0\right) + 2 \ME_{\text{Poi}} \left[ \frac{I\left(\sum_j n\hat{p}_j I(0 < \hat{\pi}_j < 1)\geq 1\right)}{\sum_j n\hat{p}_j I(0 < \hat{\pi}_j < 1)}\right]\\
        \leq &\, 2 \exp\left(-C(\epsilon)\frac{n^2}{n \vee d} \right) + 4\ME_{\text{Poi}} \left[ \frac{1}{1+\sum_j n\hat{p}_j I(0 < \hat{\pi}_j < 1)} \right],
    \end{aligned}
    \end{equation}
    where we use the fact that $\sum_j n\hat{p}_j I(0 < \hat{\pi}_j < 1)$ is an integer and the last inequality follows from Lemma \ref{lemma:nocollision} and the inequality
    \[
    \frac{I(x \geq 1)}{x} \leq \frac{2}{1+x}, \, x\geq 0.
    \]
    $\ME_{\text{Poi}}$ means the components of $(n\hat{p}_1, \dots, n\hat{p}_d)$ are independent and $n\hat{p}_k \sim \text{Poi}(np_k)$. Let $W = \sum_j n\hat{p}_j I(0 < \hat{\pi}_j < 1)$ and we will derive a tail bound for $1/(1+W)$ by bounding its MGF $\ME_{\text{Poi}}[\exp(-c_1W)]$ for an absolute constant $c_1>0$ (e.g. can be taken as $1/2$). Note that given $\hat{p}_k$, $n\hat{p}_k \hat{\pi}_k$ is a binomial variable and 
    \[
    \MP(0 < \hat{\pi}_k < 1 \mid \hat{p}_k) = 1-\pi_k^{n\hat{p}_k} - (1-\pi_k)^{n\hat{p}_k}
    \]
    when $\hat{p}_k>0$. Denote $q_k(j):= 1-\pi_k^{j} - (1-\pi_k)^{j}$. We have
    \[
    \begin{aligned}
        &\,\ME_{\text{Poi}}[\exp(-c_1n \hat{p}_k I(0<\hat{\pi}_k <1))] \\
        \leq &\, \exp(-np_k)(1+np_k) +\sum_{j=2}^\infty \ME_{\text{Poi}}[\exp(-c_1n \hat{p}_k I(0<\hat{\pi}_k <1))\mid n\hat{p}_k = j] \exp(-np_k) \frac{(np_k)^j}{j!}\\
        =&\, \exp(-np_k)(1+np_k) +\sum_{j=2}^\infty \left[ \exp(-c_1j)q_k(j)+1-q_k(j) \right] \exp(-np_k) \frac{(np_k)^j}{j!}.
    \end{aligned}
    \]
    For $c_2>0$ another constant to be fixed, we divide into two cases:

    \textbf{Case 1}: $np_k \leq c_2$, we have
    \[
    \begin{aligned}
        &\,\ME_{\text{Poi}}[\exp(-c_1n \hat{p}_k I(0<\hat{\pi}_k <1))]\\
        \leq &\, 1-\sum_{j=2}^\infty q_k(j)(1-\exp(-c_1j))\exp(-np_k)\frac{(np_k)^j}{j!}\\
        \leq &\, 1-q_k(2)(1-\exp(-2c_1))\exp(-np_k)\frac{(np_k)^2}{2} \\
        \leq &\, 1-\epsilon (1-\epsilon)(1-\exp(-2c_1))\exp(-c_2){(np_k)^2} \\
        \leq &\, \exp(-c_3 n^2p_k^2),
    \end{aligned}
    \]
    where $c_3 = \epsilon (1-\epsilon)(1-\exp(-2c_1))\exp(-c_2)$.

    \textbf{Case 2}: $np_k > c_2$, we have
    \[
    \begin{aligned}
        &\,\ME_{\text{Poi}}[\exp(-c_1n \hat{p}_k I(0<\hat{\pi}_k <1))]\\
        \leq &\, \exp(-np_k)(1+np_k)  +\sum_{j=2}^\infty \exp(-c_1j)\exp(-np_k) \frac{(np_k)^j}{j!} \\
        &\,+\sum_{j=2}^\infty (1-q_k(j))(1-\exp(-c_1j)) \exp(-np_k) \frac{(np_k)^j}{j!}\\
    \end{aligned}
    \]
    The summation of first two terms is equal to 
    \[
    \begin{aligned}
        &\,\exp(-np_k)(1+np_k)  +\exp(-np_k) \left [\exp\{\exp(-c_1)np_k\}-1-\exp(-c_1)np_k \right]\\
        =&\, \exp\{-(1-\exp(-c_1))np_k\} +np_k \exp(-np_k)(1-\exp(-c_1)).
    \end{aligned}
    \]
    The third term can be bounded as
    \[
    \begin{aligned}
        &\,\sum_{j=2}^\infty (1-q_k(j))(1-\exp(-c_1j)) \exp(-np_k) \frac{(np_k)^j}{j!}\\
        \leq &\, \sum_{j=2}^\infty (1-q_k(j)) \exp(-np_k) \frac{(np_k)^j}{j!}\\
        = &\, \sum_{j=2}^\infty (\pi_k^j + (1-\pi_k)^j) \exp(-np_k) \frac{(np_k)^j}{j!} \\
        \leq &\,\sum_{j=2}^\infty (\epsilon^j + (1-\epsilon)^j) \exp(-np_k) \frac{(np_k)^j}{j!} \\
        \leq &\,\sum_{j=0}^\infty (\epsilon^j + (1-\epsilon)^j) \exp(-np_k) \frac{(np_k)^j}{j!} \\   
        =&\, \exp(-(1-\epsilon)np_k) + \exp(-\epsilon np_k).
    \end{aligned}
    \]
    Hence we have
    \[
    \begin{aligned}
    &\,\ME_{\text{Poi}}[\exp(-c_1n \hat{p}_k I(0<\hat{\pi}_k <1))]\\
    \leq &\, \exp\{-(1-\exp(-c_1))np_k\} +np_k \exp(-np_k)(1-\exp(-c_1)) + \exp(-(1-\epsilon)np_k) + \exp(-\epsilon np_k)\\
    \leq &\, \exp(-c_4np_k)
    \end{aligned}
    \]
    where we take $c_2>0$ sufficiently large and $c_4>0$ sufficiently small (both depend on $\epsilon$) such that 
    \[
    \exp\{-(1-\exp(-c_1))x\} +x \exp(-x)(1-\exp(-c_1)) + \exp(-(1-\epsilon)x) + \exp(-\epsilon x)\leq \exp(-c_4x) , \, x\geq c_2.
    \]
    Denote the indices corresponding to case 1 as $I_1$ and case 2 as $I_2$. Now we have
    \[
    \begin{aligned}
    &\,\ME_{\text{Poi}}[\exp(-c_1W)]\\
    = &\, \prod_{k=1}^d \ME_{\text{Poi}}[\exp(-c_1n\hat{p}_k I(0 < \hat{\pi}_k < 1))]\\
    \leq &\, \exp \left\{-\left(c_3 \sum_{k \in I_1}n^2p_k^2 + c_4 \sum_{k \in I_2}np_k \right)\right\}\\
    = &\, \exp \left\{-\left(c_3 n^2T_1 + c_4 nS_2 \right)\right\},
    \end{aligned}
    \]
    where $T_1 = \sum_{k \in I_1}p_k^2, S_2 = \sum_{k \in I_2}p_k$. In the case $d \leq n/(2c_2)$, we have
    \[
    1-S_2 = \sum_{k \in I_1}p_k \leq \frac{c_2|I_1|}{n} \leq \frac{c_2d}{n} \leq \frac{1}{2}.
    \]
    So $S_2 \geq 1/2$ and we have
    \[
    \ME_{\text{Poi}}[\exp(-c_1W)]
     \leq \exp(-c_4n/2)
    \]
    In the case $d > n/(2c_2)$, if $S_2 \geq 1/2$, the above inequality still holds. If $S_2 < 1/2$, by Cauchy-Schwarz inequality we have
    \[
    \frac{1}{4} \leq (1-S_2)^2=\left(\sum_{k \in I_1}p_k \right)^2 \leq |I_1|T_1 \leq dT_1.
    \]
    Hence $T_1 \geq \frac{1}{4d}$ and 
    \[
    \ME_{\text{Poi}}[\exp(-c_1W)] \leq \exp\left(-c_3\frac{n^2}{4d} \right)
    \]
    Thus we always have
    \[
    \ME_{\text{Poi}}[\exp(-c_1W)] \leq \max\left\{\exp(-c_4n/2), \exp\left(-c_3\frac{n^2}{4d} \right)\right\} \leq \exp \left(-c_5 \frac{n^2}{n \vee d}\right).
    \]
    Finally for a small constant $c_6>0$ such that $c_1c_6 < c_5$ we have
    \[
    \begin{aligned}
        &\,\MP_{\text{Poi}}\left(W\leq \frac{c_6n^2}{n \vee d} \right)\\
        \leq &\, \exp\left(\frac{c_1c_6n^2}{n\vee d}\right)
        \ME_{\text{Poi}}[\exp(-c_1W)] \\
        \leq &\, \exp\left(\frac{c_1c_6n^2}{n\vee d}-\frac{c_5n^2}{n\vee d}\right) \\
        \leq &\, \exp\left(\frac{-c_7n^2}{n\vee d}\right)
    \end{aligned}
    \]
    Hence we have the following bound 
    \[
    \ME_{\text{Poi}}\left[\frac{1}{1+W} \right] \leq \frac{1}{c_6} \frac{n \vee d}{n^2} + \MP_{\text{Poi}}\left(W\leq \frac{c_6n^2}{n \vee d}\right) \\
    \leq \frac{1}{c_6} \frac{n \vee d}{n^2} + \exp\left(\frac{-c_7n^2}{n\vee d}\right).
    \]
    Plug into \eqref{eq:bound-the-inverse}, we conclude 
    \[
    \ME \left[ \frac{I\left(\sum_{j=1}^d \hat{t}_j > 0\right)}{n\sum_{j=1}^d \hat{t}_j} \right] \lesssim \exp\left(-C'(\epsilon)\frac{n^2}{n\vee d}\right) + \frac{n \vee d}{n^2},
    \]
    which is also a bound on expected conditional variance. Combining this bound with \eqref{eq:var-cond-exp} we have
    \[
    \operatorname{Var}(\hat{\tau}) \lesssim \sigma_n^2 + \exp\left(-C'(\epsilon)\frac{n^2}{n\vee d}\right) + \frac{n \vee d}{n^2}.
    \]
    Thus the mean squared error of $\hat{\tau}$ can be bounded as
    \[
    \ME[(\hat{\tau}-\psi)^2] \lesssim \sigma_n^2 + \exp\left(-C'(\epsilon)\frac{n^2}{n\vee d}\right) + \frac{n \vee d}{n^2}.
    \]
\end{proof}

\subsection{Proof of Theorem \ref{thm:homo-minimax}}

\begin{proof}
The parametric rate $1/n$  can be similarly achieved by using the construction in proof of Theorem \ref{thm:ate-minimax}.
We focus on proving the $d/n^2$ lower bound when $d \lesssim n^2$.

By Lemma \ref{lemma:poission-minimax}, we can again focus on a Poissonized experiment where the total sample size $N \sim \mathrm{Poi}(n)$. Since $p_k = 1/d$, this is equivalent to having independent category counts
\[
N_k \sim \mathrm{Poi}(\lambda),
\qquad
\lambda := \frac{n}{d},
\qquad
1 \le k \le d.
\]
We then define the null and alternative distributions.
Fix $q \in (0,q_0]$, where
\[
q_0 := \min\left\{\frac{1}{8},\, 8\Bigl(\frac{1}{2}-\epsilon\Bigr)^2\right\},
\qquad
r := \sqrt{\frac{q}{8}}.
\]
Define the null model $P_0 \in \mathcal{H}_d(\epsilon)$ by
\[
p_k = \frac{1}{d},
\qquad
\pi_k = \frac{1}{2},
\qquad
\mu_{0k} = \mu_{1k} = \frac{1}{2},
\qquad
1 \le k \le d.
\]
Then $\psi(P_0)=0$.
For the alternative, let $\xi_1,\dots,\xi_d$ be i.i.d.\ Rademacher random variables. Given $\boldsymbol{\xi}=(\xi_1,\dots,\xi_d)$, define the model $P_\xi$ by
\[
p_k = \frac{1}{d},
\qquad
\pi_k = \frac{1}{2} + r \xi_k,
\]
\[
\mu_{0k} = \frac{1}{2} - r \xi_k - \frac{q}{4},
\qquad
\mu_{1k} = \frac{1}{2} - r \xi_k + \frac{q}{4},
\qquad
1 \le k \le d.
\]
Since $r \le \frac{1}{2}-\epsilon$ and $r + q/4 \le 1/2$, every realization $P_\xi$ belongs to $\mathcal{H}_d(\epsilon)$. Moreover,
\[
\mu_{1k} - \mu_{0k} = \frac{q}{2}
\qquad \text{for all } k \in [d],
\]
so that
\[
\psi(P_\xi) = \frac{q}{2}
\qquad \text{almost surely}.
\]

With this construction, we will see that one observation from a category is completely uninformative on average.
Order the four cells as
\[
(A,Y)\in\{(1,1),(1,0),(0,1),(0,0)\},
\]
and write
\[
u_0 := \left(\frac{1}{4},\frac{1}{4},\frac{1}{4},\frac{1}{4}\right).
\]
Under the null model $P_0$, the cell probabilities are exactly $u_0$.

Under the alternative $P_\xi$, for category $k$ a direct calculation gives the following cell probabilities
\[
u_{\xi_k}
=
\left(
\frac{1}{4}+2r^3\xi_k,\ 
\frac{1}{4}+(r-2r^3)\xi_k,\ 
\frac{1}{4}+(-r+2r^3)\xi_k,\ 
\frac{1}{4}-2r^3\xi_k
\right).
\]
Equivalently,
\[
u_{\xi_k} = u_0 + \xi_k \Delta,
\qquad
\Delta := \bigl(2r^3,\ r-2r^3,\ -r+2r^3,\ -2r^3\bigr).
\]
Hence
\[
\frac{1}{2}u_{+1} + \frac{1}{2}u_{-1} = u_0.
\]
Therefore, after averaging over the latent sign $\xi_k$, the law of \emph{one} observation from category $k$ is exactly the same under the alternative prior and the null. Thus categories observed zero or one time contain no information for distinguishing the null from the alternative; only repeated observations within a category can separate them.

We then compute the $\chi^2$ divergence between the null and marginalized alternative $P_1 = \ME[P_{\xi}]$.
Under the null $P_0$, the four cell counts in category $k$,
\[
\bN_k = (N_{k,11},N_{k,10},N_{k,01},N_{k,00}),
\]
are independent with law
\[
Q_0 := \bigotimes_{j=1}^4 \mathrm{Poi}(\lambda/4).
\]
Under the alternative prior, conditional on $\xi_k=\pm 1$, the category law is
\[
Q_\pm := \bigotimes_{j=1}^4 \mathrm{Poi}\!\bigl(\lambda(u_{0j}\pm \Delta_j)\bigr),
\]
and after averaging over $\xi_k$,
\[
Q_1 := \frac{1}{2}Q_+ + \frac{1}{2}Q_-.
\]
Let
\[
s_q := \|\Delta\|_2^2.
\]
A direct computation gives
\[
s_q
=
(2r^3)^2 + (r-2r^3)^2 + (-r+2r^3)^2 + (-2r^3)^2
=
\frac{q}{4} - \frac{q^2}{8} + \frac{q^3}{32}
\le \frac{q}{4}.
\]
We now compute $\chi^2(Q_1 \|Q_0)$. Since
\[
u_{0j}\pm \Delta_j = \frac{1}{4}(1\pm 4\Delta_j),
\]
the likelihood ratios satisfy
\[
\frac{dQ_\pm}{dQ_0}(\bN_k)
=
\prod_{j=1}^4 (1\pm 4\Delta_j)^{N_{k,j}},
\]
because $\sum_{j=1}^4 \Delta_j = 0$.
Therefore
\[
\mathbb{E}_{Q_0}\!\left[\left(\frac{dQ_{\pm}}{dQ_0}\right)^2\right]
=
\prod_{j=1}^4
\exp\!\left(
\frac{\lambda}{4}\bigl[(1\pm 4\Delta_j)^2 - 1\bigr]
\right)
=
\exp(4\lambda s_q).
\]
Similarly,
\[
\mathbb{E}_{Q_0}\!\left[
\frac{dQ_+}{dQ_0}\frac{dQ_-}{dQ_0}
\right]
=
\prod_{j=1}^4
\exp\!\left(
\frac{\lambda}{4}\bigl[(1+4\Delta_j)(1-4\Delta_j)-1\bigr]
\right)
=
\exp(-4\lambda s_q).
\]
Hence
\[
1+\chi^2(Q_1\|Q_0)
=
\mathbb{E}_{Q_0}\!\left[\left(\frac{dQ_1}{dQ_0}\right)^2\right]
=
\frac{1}{4}
\left(
e^{4\lambda s_q}+2e^{-4\lambda s_q}+e^{4\lambda s_q}
\right)
=
\cosh(4\lambda s_q).
\]
Since $4\lambda s_q \le \lambda q$, and later we will choose $q$ so that $\lambda q \le 1$, it follows that
\[
\chi^2(Q_1\|Q_0)
=
\cosh(4\lambda s_q)-1
\leq C' \lambda^2 s_q^2 \leq 
C(\lambda q)^2
\]
for some universal constant $C>0$.

Because the categories are independent under Poissonization and the latent signs $\xi_k$ are independent, the cell counts' laws satisfy
\[
{P}_0 = Q_0^{\otimes d}.
\qquad
{P}_1 = Q_1^{\otimes d},
\]
Therefore
\[
1+\chi^2({P}_1\|{P}_0)
=
\bigl(1+\chi^2(Q_1\|Q_0)\bigr)^d
\le
\exp\!\bigl(Cd\lambda^2 q^2\bigr)
=
\exp\!\left(C\frac{n^2q^2}{d}\right).
\]
Now choose
\[
q^2 = c_0\frac{d}{n^2},
\]
where $c_0>0$ is a sufficiently small constant depending only on $\epsilon$, chosen so that $q\le q_0$ and $\lambda q \le 1$. Then the $\chi^2$ divergence is bounded away from infinity. By method of fuzzy hypothesis,
we conclude that
\[
\inf_{\hat{\psi}} \sup_{\MP \in \mathcal{H}(\epsilon)}\ME_{\MP}\left[\left(\hat{\psi} - \psi\right)^2\right] \gtrsim
q^2
\asymp
\frac{d}{n^2}.
\]
\end{proof}

\subsection{Proof of Theorem \ref{thm:ate-minimax-cova}}
\begin{proof}
By the same argument in Lemma \ref{lemma:poission-minimax}, we only need to prove the result under the Poisson sampling model, where $(n \hat{p}_1,\dots,n \hat{p}_d) \text{ i.i.d} \sim \text{Poi}(\lambda)$ with $\lambda = \frac{n}{d}$.
Let $\nu_0, \nu_1$ be two prior distributions on propensity score $\pi$ defined on $[\epsilon, 1-\epsilon]$ satisfying 
\begin{itemize}
    \item $\ME_{\nu_0}[\pi^j] = \ME_{\nu_1}[\pi^j], 1\leq j \leq 2L, L = \lfloor \frac{1}{\beta} \rfloor +1$.
    \item $\ME_{\nu_0}\left[\frac{1}{\pi} \right] - \ME_{\nu_1}\left[\frac{1}{\pi} \right] = c_1(\beta, \epsilon)>0.$
\end{itemize}
The existence of $\nu_0, \nu_1$ follows from the duality between moment matching and best polynomial approximation. In fact, one can show 
\begin{equation}\label{eq:const-1/x}
c_1(\beta, \epsilon) = 2 \ME_{2L} \left( \frac{1}{x}; [\epsilon, 1-\epsilon] \right),    
\end{equation}
where $E_n(f;S)$ is the best polynomial (with order no greater than $n$) approximation error of $f$ on the interval $S$. Since $L$ depends on $\beta$, the RHS of \eqref{eq:const-1/x} is a constant only depending on $\beta,\epsilon$. 

Under null hypothesis $P$, let $ (\pi_1, \dots, \pi_d)$ i.i.d. $\sim \nu_0$ and $\mu_k = \epsilon /\pi_k$. The sufficient statistics for $\psi_1$ (conditioning on $(\pi_1, \dots, \pi_d)$) are
\[
N_{k11} : = \#\{i:X_i=k,A_i=1,Y_i=1 \} \sim \text{Poi} (\epsilon \lambda)  ,
\]
\[
N_{k10} : = \#\{i:X_i=k,A_i=1,Y_i=0 \} \sim \text{Poi} \left(\lambda(\pi_k-\epsilon)\right) ,
\]
\[
N_{k0} := \#\{i:X_i=k,A_i=0 \} \sim \text{Poi} \left(\lambda (1-\pi_k)\right) 
\]
with $N_{k11}, N_{k10}, N_{k0} $ conditionally independent. Similarly under the alternative hypothesis $P'$, let $ (\pi_1', \dots, \pi_d')$ i.i.d. $\sim \nu_1$ and $\mu_k' = \epsilon /\pi_k'$. The sufficient statistics for $\psi_1$ (conditioning on $(\pi_1', \dots, \pi_d')$) are
\[
N_{k11}' : = \#\{i:X_i=k,A_i=1,Y_i=1 \} \sim \text{Poi} (\epsilon \lambda)  ,
\]
\[
N_{k10}' : = \#\{i:X_i=k,A_i=1,Y_i=0 \} \sim \text{Poi} \left(\lambda(\pi_k'-\epsilon)\right) ,
\]
\[
N_{k0}' := \#\{i:X_i=k,A_i=0 \} \sim \text{Poi} \left(\lambda (1-\pi_k')\right) .
\]
Using the same notation as in proof of Theorem \ref{thm:ate-minimax}, the total variation distance between the sufficient statistics is 
\[
\begin{aligned}
    &\, \text{TV}(\bN_k,\bN_k')\\
    = &\, \frac{1}{2} \sum_{i,j,\ell=0}^{\infty}  \big |  \MP (N_{k11}=i, N_{k10}=j, N_{k0}=\ell)- \MP (N_{k11}'=i, N_{k10}'=j, N_{k0}'=\ell)\big |.
\end{aligned}
\]
Conditioning on $\pi_k$ we have
\[
\begin{aligned}
    &\,\MP (N_{k11}=i, N_{k10}=j, N_{k0}=\ell)\\
    = & \, \ME\left [ \exp(-\epsilon \lambda ) \frac{(\epsilon \lambda)^i}{i!} \exp(-\lambda(\pi_k-\epsilon)) \frac{[\lambda(\pi_k-\epsilon)]^j}{j!} \right.\\
    &\, \left.\exp(-\lambda (1-\pi_k)) \frac{[\lambda(1-\pi_k)]^{\ell}}{\ell !} \right ] \\
    = &\, \ME \left[ \exp(-\lambda) \frac{(\epsilon \lambda)^i[\lambda(\pi_k-\epsilon)]^j[\lambda(1-\pi_k)]^{\ell}}{i!j!\ell!} \right]. \\
\end{aligned}
\]
Hence the total variation distance can be written as
    \[
    \begin{aligned}
    &\, \text{TV}(\bN_k,\bN_k')  \\
    = &\, \frac{1}{2} \sum_{i,j,\ell} \left | \ME \left[ \exp(-\lambda) \frac{(\epsilon \lambda)^i[\lambda(\pi_k-\epsilon)]^j[\lambda(1-\pi_k)]^{\ell}}{i!j!\ell!} \right] \right. \\
    &\, \left.- \ME \left[ \exp(-\lambda) \frac{(\epsilon \lambda)^i[\lambda(\pi_k'-\epsilon)]^j[\lambda(1-\pi_k')]^{\ell}}{i!j!\ell!} \right] \right | \\
    = &\,  \frac{1}{2} \sum_{j,\ell} \left | \ME \left [ \exp(-(1-\epsilon)\lambda) \frac{[\lambda(\pi_k-\epsilon)]^j[\lambda(1-\pi_k)]^{\ell}}{j!\ell!}\right]\right. \\
    &\, \left. -\ME \left [ \exp(-(1-\epsilon)\lambda) \frac{[\lambda(\pi_k'-\epsilon)]^j[\lambda(1-\pi_k')]^{\ell}}{j!\ell!}\right] \right|,
    \end{aligned}
    \]
where the last equation follows from $\sum_i \exp(-\epsilon \lambda) \frac{(\epsilon \lambda)^i}{i!} = 1$. Since $\ME[\pi_k^j] = \ME[\pi_k'^j], 1 \leq j \leq 2L$, we have 
\[
\begin{aligned}
    &\, \text{TV}(\bN_k,\bN_k')  \\
    = &\,  \frac{1}{2} \sum_{j+\ell>2L} \left | \ME \left [ \exp(-(1-\epsilon)\lambda) \frac{[\lambda(\pi_k-\epsilon)]^j[\lambda(1-\pi_k)]^{\ell}}{j!\ell!}\right]\right. \\
    &\, \left. -\ME \left [ \exp(-(1-\epsilon)\lambda) \frac{[\lambda(\pi_k'-\epsilon)]^j[\lambda(1-\pi_k')]^{\ell}}{j!\ell!}\right] \right| \\
    \leq &\, \frac{1}{2} \left \{ \sum_{j+\ell>2L} \ME \left [ \exp(-(1-\epsilon)\lambda) \frac{[\lambda(\pi_k-\epsilon)]^j[\lambda(1-\pi_k)]^{\ell}}{j!\ell!}\right] \right. \\
    &\,  \left. +\ME \left [ \exp(-(1-\epsilon)\lambda) \frac{[\lambda(\pi_k'-\epsilon)]^j[\lambda(1-\pi_k')]^{\ell}}{j!\ell!}\right] \right \}.  \\
\end{aligned}
\]
Since given $\pi_k$, $N_{k11} \sim  \text{Poi}(\lambda(\pi_k-\epsilon)), N_{k10} \sim  \text{Poi}(\lambda(1-\pi_k)), N_{k11}+N_{k10} \sim \text{Poi}(\lambda(1-\epsilon))$, we have
\[
\begin{aligned}
 &\,\sum_{j+\ell>2L} \ME \left [ \exp(-(1-\epsilon)\lambda) \frac{[\lambda(\pi_k-\epsilon)]^j[\lambda(1-\pi_k)]^{\ell}}{j!\ell!}\right]\\
 = &\, \ME[\MP(N_{k11}+N_{k10}>2L|\pi_k)] \\
 \leq&\,  \exp(-\lambda(1-\epsilon)) \left( \frac{e\lambda(1-\epsilon)}{2L} \right)^{2L},
\end{aligned}
\]
where the last inequality follows from Chernoff bound \eqref{eq:chernoff-poi}. Hence in the regime $n \lesssim d^{1-\beta}$ with the choice $L=\lfloor 1/\beta \rfloor +1$, the total variation distance is bounded as 
\[
\text{TV}(\bN_k,\bN_k')  \leq \left( \frac{en(1-\epsilon)}{2dL} \right)^{2L} \lesssim d^{-2\beta L} = o(1/d^2).
\]
Since $(\bN_1, \dots, \bN_d)$ are independent, we have
\[
\text{TV}(\bN,\bN') \leq \sum_{k=1}^d \text{TV}(\bN_k,\bN_k') = o(1/d)  \rightarrow 0.
\]
The functional separation between the null and the alternative is
\[
\psi_1(P) - \psi_1(P') = \sum_{k=1}^d \frac{\epsilon}{d} \left( \frac{1}{\pi_k} - \frac{1}{\pi_k'} \right)
\]
with expectation
\[
\ME[\psi_1(P) - \psi_1(P') ] = \epsilon \ME \left[ \frac{1}{\pi_k} - \frac{1}{\pi_k'} \right]= \epsilon c_1( \beta, \epsilon) := c_2(\beta, \epsilon) = c_2.
\]
Define two events
\[
E =\left \{ \left| \frac{1}{d}\sum_{k=1}^d \frac{1}{\pi_k} - \ME\left[\frac{1}{\pi_k}\right] \right| \leq \frac{c_2}{4} \right\},
\]
\[
E' =\left \{ \left| \frac{1}{d}\sum_{k=1}^d \frac{1}{\pi_k'} - \ME\left[\frac{1}{\pi_k'}\right] \right| \leq \frac{c_2}{4} \right\}.
\]
By Chebyshev's inequality, we have
\[
\MP(E^c) \leq \frac{16}{c_2^2d} \operatorname{Var}\left(\frac{1}{\pi_k} \right) \leq \frac{16}{\epsilon^2c_2^2d}.
\]
\[
\MP(E'^c) \leq \frac{16}{c_2^2d} \operatorname{Var}\left(\frac{1}{\pi_k'} \right) \leq \frac{16}{\epsilon^2c_2^2d}.
\]
We put the following prior distributions induced by $\{ \pi_k, 1 \leq k \leq d\}$ and $\{ \pi_k', 1 \leq k \leq d\}$ on $P$ and $P'$, respectively:
\[
\pi \stackrel{d}{=} P \mid E, \pi' \stackrel{d}{=} P' \mid E'.
\]
Note that under $\pi, \pi'$,
\[
|\psi_1(P)-\psi_1(P')|\geq c_2/2.
\]
By triangle inequality, the total variation distance of the sufficient counting statistics $\bN$ and $\bN'$ under two priors is bounded by
\[
\begin{aligned}
     \text{TV}\left( {\bN \mid E}, {\bN' \mid E'} \right) \leq &\, \text{TV}\left( {\bN \mid E}, {\bN} \right) + \text{TV}({\bN},{\bN'} ) + \text{TV}\left( {\bN' \mid E'}, {\bN'} \right)\\
     \leq &\, \MP(E^c) + \MP(E'^c) + \text{TV}({\bN},{\bN'} )\\
     \leq &\,  \frac{32}{ \epsilon^2 c_2^2 d} + \text{TV}({\bN},{\bN'} ) \rightarrow 0 .
\end{aligned}
\]
By method of fuzzy hypotheses, we conclude
\[
\begin{aligned}
&\, \inf_{\hat{\psi}_1} \sup_{\MP \in \mathcal{D}^U (\epsilon)}  \ME_{\MP}\left[\left(\hat{\psi}_1 - \psi_1\right)^2\right] \\
\geq &\, \frac{c_2^2}{32} \left( 1 -  \text{TV}\left( {\bN \mid E}, {\bN' \mid E'} \right) \right).\\
\geq &\, c(\beta, \epsilon),
\end{aligned}
\]
for some constant $c(\beta, \epsilon)$ when $d$ is large enough in the regime $n \lesssim d^{1-\beta}$.
\end{proof}

\subsection{Proof of Theorem \ref{thm:cova-dist-agnostic}}
\begin{proof}

Note that $\eta = \ME[(A-\pi_X)(Y-\mu_X)]$ The conditional bias is (Let $\bZ=(X,A,Y), \bZ_1=(X_1,A_1,Y_1), \bZ_2=(X_2,A_2,Y_2)$ be samples independent of $D$)
\[
\begin{aligned}
   &\, \ME[\hat{\eta}\mid D]-\eta \\
    =&\,  \ME[(A-\hat{\pi}_X)(Y-\hat{\mu}_X)\mid D] - \ME[(A-\pi_X)(Y-\mu_X)\mid D] \\
    &\,- \ME\left[\frac{(A_1-\hat{\pi}_{X_1})I(X_1=X_2)(Y_2 - \hat{\mu}(X_2))}{\hat{p}_{X_1}} \mid D\right]
\end{aligned}
\]
Note that
\[
\begin{aligned}
    &\, \ME[(A-\hat{\pi}_X)(Y-\hat{\mu}_X)\mid D] - \ME[(A-\pi_X)(Y-\mu_X)\mid D] \\
    = &\, \ME[(A-\hat{\pi}_X)(Y-\hat{\mu}_X)\mid D] - \ME[(A-\hat{\pi}_X)(Y-{\mu}_X)\mid D] \\
    &\, +\ME[(A-\hat{\pi}_X)(Y-{\mu}_X)\mid D] - \ME[(A-{\pi}_X)(Y-{\mu}_X)\mid D] \\
    = &\, \ME[(A-\hat{\pi}_X)(\mu_X-\hat{\mu}_X)\mid D] + \ME[(\pi_X-\hat{\pi}_X)(Y-{\mu}_X)\mid D] \\
    = &\, \ME[(\pi_X-\hat{\pi}_X)(\mu_X-\hat{\mu}_X)\mid D],
\end{aligned}
\]
the last equation follows from conditioning on $X$. By conditioning on $X_1$ we have
\[
\ME\left[ I(X_2=X_1)(Y_2 - \hat{\mu}(X_2))\mid D,X_1\right] = p_{X_1}(\mu_{X_1}-\hat{\mu}_{X_1}). 
\]
Hence condition on $\bZ_1$ we have
\[
\begin{aligned}
    &\, \ME\left[\frac{(A_1-\hat{\pi}_{X_1})I(X_1=X_2)(Y_2 - \hat{\mu}(X_2))}{\hat{p}_{X_1}} \mid D\right]\\
    = &\, \ME\left[ \frac{p_{X_1}}{\hat{p}_{X_1}}(A_1-\hat{\pi}_{X_1})(\mu_{X_1}-\hat{\mu}_{X_1}) \mid D\right] \\
    = &\, \ME\left[ \frac{p_{X_1}}{\hat{p}_{X_1}}(\pi_{X_1}-\hat{\pi}_{X_1})(\mu_{X_1}-\hat{\mu}_{X_1}) \mid D\right].
\end{aligned}
\]
Hence the conditional bias is 
\begin{equation}\label{eq:high-order-bias}
    \ME[\hat{\eta}\mid D]-\eta = \ME \left[ (\hat{\mu}_{X}-\mu_{X})(\hat{\pi}_{X} - \pi_X)\left (1-\frac{p_{X}}{\hat{p}_X} \right)  \mid D \right].
\end{equation}
By Cauchy-Schwarz inequality one can bound the conditional bias as
\[
|\ME[\hat{\eta}\mid D]-\eta| \leq  \|\hat{\mu}-\mu\|_2 \|\hat{\pi}-\pi\|_2 \max_{k} \left| 1-\frac{p_k}{\hat{p}_k} \right|
\]
We next bound the variance. Let \[
h_2(\bZ_1,\bZ_2)  = (A_1 - \hat{\pi}_{X_1})(Y_1-\hat{\mu}_{X_1}) - \frac{(A_1-\hat{\pi}_{X_1})I(X_1=X_2)(Y_2-\hat{\mu}_{X_2})}{\hat{p}_{X_1}}
\]
and note that the estimator can be expressed as 
\[
\hat{\eta} = \MU_n[h_2(\bZ_1,\bZ_2)]. 
\]
We will use the variance of U-statistics (Lemma 6 in \cite{robins2009quadratic}) to bound the conditional variance of $\hat{\eta}$. Let
\[
h_1(\bZ_1) = \ME[h_2(\bZ_1,\bZ_2)\mid\bZ_1] = (A_1 - \hat{\pi}_{X_1})(Y_1-\hat{\mu}_{X_1}) - \frac{(A_1-\hat{\pi}_{X_1})p_{X_1}(\mu_{X_1}-\hat{\mu}_{X_1})}{\hat{p}_{X_1}}.
\]
\[
\begin{aligned}
    \sigma_1^2 =&\, \ME[h_1^2(\bZ_1)\mid D] \\
    \leq &\, 2 \left( \ME[(A_1 - \hat{\pi}_{X_1})^2(Y_1-\hat{\mu}_{X_1})^2\mid D] + \ME \left[ \frac{(A_1-\hat{\pi}_{X_1})^2p_{X_1}^2(\mu_{X_1}-\hat{\mu}_{X_1})^2}{\hat{p}_{X_1}^2} \mid D \right] \right) \\
    \lesssim &\, 1 + \ME \left[ \frac{p_{X}^2}{\hat{p}_X^2}  \mid D \right] \\
    \lesssim &\, 1 + \left(1+\max_{k} \left| 1-\frac{p_k}{\hat{p}_k} \right|\right)^2 \\
    \lesssim &\, \left(1+\max_{k} \left| 1-\frac{p_k}{\hat{p}_k} \right|\right)^2.
\end{aligned}
\]
\[
\begin{aligned}
    \sigma_2^2 =&\,  \ME[h_2^2(\bZ_1,\bZ_2)\mid D] \\
    \leq &\, 2 \left( \ME[(A_1 - \hat{\pi}_{X_1})^2(Y_1-\hat{\mu}_{X_1})^2\mid D] + \ME \left[ \frac{(A_1-\hat{\pi}_{X_1})^2I(X_1=X_2)(Y_2-\hat{\mu}_{X_2})^2}{\hat{p}_{X_1}^2} \mid  D \right] \right) \\
    \lesssim &\,  1 + \ME\left [ \frac{I(X_1=X_2)}{\hat{p}_{X_1}^2}  \mid  D \right] \\
    \lesssim &\, 1 + \sum_{k=1}^d \frac{p_k^2}{\hat{p}_k^2} \\
    \leq &\, 1 + d\left(1+\max_{k} \left| 1-\frac{p_k}{\hat{p}_k} \right|\right)^2 \\
    \lesssim&\,  d\left(1+\max_{k} \left| 1-\frac{p_k}{\hat{p}_k} \right|\right)^2.
\end{aligned}
\]
By Lemma 6 in \cite{robins2009quadratic} the conditional variance is upper bounded as 
\[
\operatorname{Var}(\hat{\eta}\mid D) \lesssim \frac{1}{n}\left(1+\max_{k} \left| 1-\frac{p_k}{\hat{p}_k} \right|\right)^2 + \frac{d}{n^2}\left(1+\max_{k} \left| 1-\frac{p_k}{\hat{p}_k} \right|\right)^2.
\]
\end{proof}

\subsection{Proof of Theorem \ref{thm:cova-dist}}
\begin{proof}
    
We will use the following equation for MSE 
\begin{equation}\label{eq:mse-decomposition}
    \begin{aligned}
    &\, \ME[(\hat{\eta}-\eta)^2] \\
    =&\, \ME \left\{ \ME[(\hat{\eta}-\eta)^2\mid D] \right\} \\
    = &\, \ME[\operatorname{Var}(\hat{\eta}\mid D)] + \ME\{(\ME[\hat{\eta}\mid D] - \eta)^2\}.
\end{aligned}
\end{equation}
For the conditional bias derived in \eqref{eq:high-order-bias} we have (using the property $(\ME[X])^2 \leq \ME[X^2]$)
\[
\begin{aligned}
    &\,\ME\{(\ME[\hat{\eta}\mid D] - \eta)^2\}\\
   \leq &\,  \ME \left \{ \ME \left [ (\hat{\mu}_{X}-\mu_{X})^2(\hat{\pi}_{X} - \pi_X)^2\left (1-\frac{p_{X}}{\hat{p}_X} \right)^2  \mid  D\right] \right\} \\
   = &\, \ME \left[ \sum_{k=1}^d p_k (\hat{\mu}_{k}-\mu_{k})^2(\hat{\pi}_{k} - \pi_k)^2\left (1-\frac{p_{k}}{\hat{p}_k} \right)^2\right] \\
   = &\, \sum_{k=1}^d \ME \left[ p_k (\hat{\mu}_{k}-\mu_{k})^2(\hat{\pi}_{k} - \pi_k)^2\left (1-\frac{p_{k}}{\hat{p}_k} \right)^2\right] \\
   \leq&\, \xi_n^2 \sum_{k=1}^d p_k\ME [(\hat{\pi}_k - \pi_k)^2],
\end{aligned}
\]
where in the last inequality we use the bound 
\[
\max_k \left |1-\frac{p_k}{\hat{p}_k} \right| \leq \xi_n. 
\]
and $(\hat{\mu}_k - \mu_k)^2 \leq 1$. A naive bound 
\[
\sum_{k=1}^d p_k\ME [(\hat{\pi}_k - \pi_k)^2] \leq 1
\]
holds for both empirical average estimator $\hat{\pi}_k, \hat{\mu}_k$ or simply letting $\hat{\pi}_k = \hat{\mu}_k =0$. For empirical average estimator $\hat{\pi}_k$ from a sample of size $n$ we can derive an alternative bound. 
By the property of conditional variance, we have
\[
\operatorname{Var}(\hat{\pi}_k) = \ME[\operatorname{Var}(\hat{\pi}_k \mid  \bX^n)] + \operatorname{Var}(\ME[\hat{\pi}_k \mid  \bX^n]).
\]
Recall $\ME[\hat{\pi}_k \mid  \bX^n] = \pi_k I(\hat{p}_k>0), \operatorname{Var}(\hat{\pi}_k \mid  \bX^n) = \frac{\pi_k(1-\pi_k)}{n\hat{p}_k} I(\hat{p}_k>0)$ we obtain
\[
\operatorname{Var}(\hat{\pi}_k) \leq \frac{1}{4} \ME\left[\frac{I(\hat{p}_k>0)}{n\hat{p}_k}\right] +  \pi_k^2 (1-p_k)^n \leq \frac{1}{2(n+1)p_k} + (1-p_k)^n,
\]
\[
(\ME[\hat{\pi}_k ] - \pi_k)^2 = \pi_k^2 (1-p_k)^{2n} \leq  (1-p_k)^{2n}.
\]
Thus we have
\[
\ME [(\hat{\pi}_k - \pi_k)^2] \leq \frac{1}{2(n+1)p_k} + (1-p_k)^n + (1-p_k)^{2n}
\]
\[
\sum_{k=1}^d p_k\ME [(\hat{\pi}_k - \pi_k)^2] \leq  \frac{d}{2(n+1)} + \sum_{k=1}^d p_k(1-p_k)^n +\sum_{k=1}^d p_k(1-p_k)^{2n} .
\]
Let $f_6(x) = x(1-x)^n, x\in [0,1]$ and $f_6'(x) = (1-nx)(1-x)^{n-1}$. Hence $f_6(x) \leq f_6(1/n) < 1/n$, which implies
\[
 \sum_{k=1}^d p_k(1-p_k)^n \leq \frac{d}{n}, \sum_{k=1}^d p_k(1-p_k)^{2n} \leq \frac{d}{2n}.
\]
We conclude that
\[
\sum_{k=1}^d p_k\ME [(\hat{\pi}_k - \pi_k)^2] \leq  \frac{2d}{n},
\]
Combining this bound the naive bound
\[
\sum_{k=1}^d p_k\ME [(\hat{\pi}_k - \pi_k)^2] \leq 1
\]
we have
\[
\ME\{(\ME[\hat{\eta}\mid D] - \eta)^2\} \lesssim  \xi_n^2 \frac{n \wedge d}{n}. 
\]
The bound on conditional variance in proof of Theorem \ref{thm:cova-dist-agnostic} can be reduced to
\[
\operatorname{Var}(\hat{\eta}\mid D) \lesssim \frac{1}{n} + \frac{d}{n^2}.
\]
The proof is completed by combining the bounds on conditional bias with variance as in \eqref{eq:mse-decomposition}.

\end{proof}

\subsection{Proof of Theorem \ref{thm:plugin-lowdim}}
\begin{proof}
    We will show the plug-in estimator of $\psi_1 = \sum_{k}p_k \mu_{1k}$ satisfies 
    \[
    \hat{\psi}_1 - \psi_1 = (\MP_n - \MP) [\varphi_1(\bZ)] + o_{\MP}(1/\sqrt{n}),
    \]
    where
    \[
    \varphi_1(\bZ) = \frac{A(Y-\mu_{1X})}{\pi_X}+\mu_{1X}.
    \]
    Similar argument can be applied to $\psi_0 = \sum_{k}p_k \mu_{0k}$. By Proposition \ref{prop:ests-equivalence} we will use the doubly robust form of $\hat{\psi} = \MP_n[\hat{\varphi}_1(\bZ)] $. Using the following decomposition
    \[
    \hat{\psi}_1 - \psi_1 = (\MP_n - \MP) [\varphi_1(\bZ)] + (\MP_n - \MP) [\hat{\varphi}_1(\bZ)-\varphi_1(\bZ)] + \MP [\hat{\varphi}_1(\bZ)-\varphi_1(\bZ)].
    \]
    \textbf{Step1: Bound $\MP [\hat{\varphi}_1(\bZ)-\varphi_1(\bZ)]$.} We first show $\MP [\hat{\varphi}_1(\bZ)-\varphi_1(\bZ)] = \int \hat{\varphi}_1(\bz)-\varphi_1(\bz) d \MP(\bz) = o_{\MP}(1/\sqrt{n})$.
    By direct calculations, one can show
    \[
    |\MP [\hat{\varphi}_1(\bZ)-\varphi_1(\bZ)] |= \bigg |\MP \left[ (\hat{\mu}_{1X}-\mu_{1X})\left(1-\frac{\pi_X}{\hat{\pi}_X} \right) \right] \bigg |\leq \|\hat{\mu}_{1X}-\mu_{1X}\|_2 \left \| 1-\frac{\pi_X}{\hat{\pi}_X}\right \|_2.
    \]
    Note that since $n\hat{p}_k \hat{\pi}_k \hat{\mu}_{1k}| \bX^n, \bA^n \sim \text{Bin}(n\hat{p}_k \hat{\pi}_k, \mu_{1k})$, we have
    \[
    \ME[\hat{\mu}_{1k}\mid\bX^n,\bA^n] = \mu_{1k} I(\hat{p}_k \hat{\pi}_k>0), \, \operatorname{Var}(\hat{\mu}_{1k}\mid\bX^n, \bA^n) = \frac{\mu_{1k}(1-\mu_{1k})I(\hat{p}_k \hat{\pi}_k>0)}{n\hat{p}_k \hat{\pi}_k}.
    \]
    The bias of $\hat{\mu}_{1k}$ is
    \[
    \ME[\hat{\mu}_{1k}-\mu_{1k}] = -\mu_{1k}\MP(\hat{p}_k \hat{\pi}_k=0).
    \]
    By conditioning on $\bX^n$, we have
    \[
    \MP(\hat{p}_k \hat{\pi}_k=0) = \ME[\MP(\hat{p}_k \hat{\pi}_k=0\mid \bX^n)] = \ME[(1-\pi_k)^{n\hat{p}_k}] = (1-p_k \pi_k)^n \leq (1-\epsilon p_k)^n,
    \]
    where we use the fact $\ME[c^V] = (1-p+pc)^n$ for $V \sim B(n,p)$. Similarly, we have
    \[
    \operatorname{Var}(\ME[\hat{\mu}_{1k}\mid\bX^n,\bA^n]) = \mu_{1k}^2 \MP(\hat{p}_k \hat{\pi}_k=0)\MP(\hat{p}_k \hat{\pi}_k>0) \leq (1-\epsilon p_k)^n.
    \]
    For the expected conditional variance,
    \[
    \ME \left[ \operatorname{Var}(\hat{\mu}_{1k}\mid\bX^n, \bA^n) \right] \leq \frac{1}{4}\ME \left[\frac{I(\hat{p}_k \hat{\pi}_k>0)}{n\hat{p}_k \hat{\pi}_k} \right]
    \]
    Note that $n\hat{p}_k \hat{\pi}_k \sim B(n,p_k \pi_k)$, by lemma \ref{lemma:bino-inverse} we have 
    \[
    \ME \left[ \operatorname{Var}(\hat{\mu}_{1k}\mid\bX^n, \bA^n) \right] \leq \frac{1}{2(n+1)p_k\pi_k}\leq \frac{1}{2\epsilon(n+1)p_k}.
    \]
    Thus the variance of $\hat{\mu}_{1k}$ is bounded by
    \[
    \operatorname{Var}(\hat{\mu}_{1k}) \leq (1-\epsilon p_k)^n + \frac{1}{2\epsilon(n+1)p_k}.
    \]
    We conclude 
    \[
    \ME[(\hat{\mu}_{1k}-\mu_{1k})^2] \leq (1-\epsilon p_k)^n + \frac{1}{2\epsilon(n+1)p_k}+ (1-\epsilon p_k)^{2n} = O(1/n)
    \]
    since when $d$ is fixed, $p_k$'s are considered as fixed. Hence
    \[
    \|\hat{\mu}_{1X}-\mu_{1X}\|_2^2 = \sum_{k}p_k (\hat{\mu}_{1k}-\mu_{1k})^2 = O_{\MP}(1/n),
    \]
    \[
    \|\hat{\mu}_{1X}-\mu_{1X}\|_2 = O_{\MP}(1/\sqrt{n}).
    \]
    Similarly $\ME[\hat{\pi}_k \mid \bX^n] = \pi_k I(\hat{p}_k>0), \operatorname{Var}(\hat{\pi}_k \mid \bX^n) = \frac{\pi_k(1-\pi_k)}{n\hat{p}_k} I(\hat{p}_k>0)$, we obtain
    \[
    \operatorname{Var}(\hat{\pi}_k) \leq \frac{1}{4} \ME\left[\frac{I(\hat{p}_k>0)}{n\hat{p}_k}\right] +  \pi_k^2 (1-p_k)^n \leq \frac{1}{2(n+1)p_k} + (1-p_k)^n,
    \]
    \[
    (\ME[\hat{\pi}_k ] - \pi_k)^2 = \pi_k^2 (1-p_k)^{2n} \leq  (1-p_k)^{2n}.
    \]
    Thus we have
    \[
    \ME [(\hat{\pi}_k - \pi_k)^2] \leq \frac{1}{2(n+1)p_k} + (1-p_k)^n + (1-p_k)^{2n} = O(1/n)
    \]
    \[
    (\hat{\pi}_k-\pi_k)^2 = O_{\MP}(1/n).
    \]
    Since $\pi_k \geq \epsilon$, this shows $1/\hat{\pi}_k = O_{\MP}(1)$ and 
    \[
    \frac{(\hat{\pi}_k-\pi_k)^2}{\hat{\pi}_k^2} =O_{\MP}(1/n).
    \]
    We conclude 
    \[
    \left \| 1-\frac{\pi_X}{\hat{\pi}_X}\right \|_2^2 = \sum_{k} p_k \frac{(\hat{\pi}_k-\pi_k)^2}{\hat{\pi}_k^2}=O_{\MP}(1/n),
    \]
    \[
    \left \| 1-\frac{\pi_X}{\hat{\pi}_X}\right \|_2=O_{\MP}(1/\sqrt{n}).
    \]
    This shows 
    \[
    \MP [\hat{\varphi}_1(\bZ)-\varphi_1(\bZ)] = O_{\MP}(1/n)=o_{\MP}(1/\sqrt{n}).
    \]
    \textbf{Step2: Bound $(\MP_n - \MP) [\hat{\varphi}_1(\bZ)-\varphi_1(\bZ)]$.}We then show 
    \[
    (\MP_n - \MP) [\hat{\varphi}_1(\bZ)-\varphi_1(\bZ)] = o_{\MP}(1/\sqrt{n}).
    \]
    Since $X$ is discrete we can write the nuisance functions $(\hat{\pi}, \hat{\mu}_1)$ as saturated linear models, i.e.
    \[
    {\pi}_x=\pi(x ; \aaa)=\aaa^{\top} \bw,
    \]
    where $\bw^{\top}=\{I(x=1), \ldots, I(x=d)\} \in\{0,1\}^{d}$ and 
    \[
    {\mu}_{1x}=\mu_1(x ; \bbb)=\bbb^{\top} \bw.
    \]
    Here the components of $\aaa$ and $\bbb$ are simply propensity scores and regression functions within different categories and $\|\bw\|_2=1$. Define the function class $\mathcal{F} = \{\varphi_{1}(\bz;\rrr), \rrr = (\aaa, \bbb), \alpha_k \in [\epsilon/2, 1], \beta_k \in [0, 1]\}$. Since propensity scores are lower bounded by $\epsilon /2$ for $\aaa \in \mathcal{F}$, one can show
    there is a constant $K$ that depends on $\epsilon$ such that 
    \[
    |\varphi_{1}(\bz;\rrr_1) - \varphi_{1}(\bz;\rrr_2)| \leq K\|\rrr_1 - \rrr_2\|_2.
    \]
    Since $d$ is fixed as constant, by example 19.7 in \cite{van2000asymptotic}, $\mathcal{F}$ is Donsker. Let 
    \[
    \tilde{\pi}_k = \hat{\pi}_k I(\hat{\pi}_k \geq \epsilon/2) + \frac{\epsilon}{2} I(\hat{\pi}_k <\epsilon/2), \, \tilde{\mu}_{1k} = \hat{\mu}_{1k},
    \]
    so that we truncate the estimated propensity score $\hat{\pi}_k$ to obtain $\tilde{\pi}_k$. Let 
    \[
    \tilde{\varphi}_1(\bZ) = \varphi_1(\bZ;\tilde{\rrr}) = \frac{A(Y-\tilde{\mu}_{1X})}{\tilde{\pi}_X}+\tilde{\mu}_{1X}.
    \]
    Clearly $\tilde{\varphi}_{1}(\cdot) = \varphi_1(\cdot, \tilde{\rrr}) \in \mathcal{F}$. We have
    \[
    \|\tilde{\rrr}- \rrr\|_2^2 =\|\tilde{\aaa} -\aaa\|_2^2+\|\tilde{\bbb}-\bbb\|_2^2 = \sum_{k=1}^d \left[(\tilde{\pi}_k-\pi_k)^2 + (\tilde{\mu}_{1k}-\mu_{1k})^2 \right] .
    \]
    Since we assume $\pi_k \geq \epsilon$, truncating $\hat{\pi}_k$ yields smaller error and hence $|\tilde{\pi}_k - \pi_k|\leq |\hat{\pi}_k - \pi_k|$. By the consistency of $\hat{\pi}_k$ and $\hat{\mu}_{1k}$ shown above, we have
    \[
    \|\tilde{\rrr}- \rrr\|_2^2 \leq \sum_{k=1}^d \left[(\hat{\pi}_k-\pi_k)^2 + (\hat{\mu}_{1k}-\mu_{1k})^2 \right] = o_{\MP}(1)
    \]
    again because $d$ is fixed. Thus
    \[
    \|\tilde{\varphi}_1-\varphi\|_2 \leq K \|\tilde{\rrr}- \rrr\| = o_{\MP}(1)
    \]
    and Lemma 19.24 in \cite{van2000asymptotic} shows 
    \begin{equation}\label{eq:varphitilde-neg}
        (\MP_n - \MP) [\tilde{\varphi}_1(\bZ)-\varphi_1(\bZ)] = o_{\MP}(1/\sqrt{n}).
    \end{equation}
    Now consider \[
    (\MP_n - \MP) [\hat{\varphi}_1(\bZ)-\tilde{\varphi}_1(\bZ)].
    \]
    By strong law of large numbers, we have 
    \[
    \hat{\pi}_k = \frac{\sum_{i=1}^n I(X_i=k,A_i=1)}{\sum_{i=1}^nI(X_i=k)} \rightarrow \pi_k
    \]
    for any $k \in [d]$ almost surely. Thus for almost every $\omega \in \Omega$ (sample space), one can find a $N(\omega,\epsilon,d) \in \bN^+$ such that for all $n \geq N(\omega, \epsilon,d)$, we have $|\hat{\pi}_k-\pi_k| \leq \epsilon/2$ for any $k \in [d]$. This together with $\pi_k \geq \epsilon$ shows $\hat{\pi}_k \geq \epsilon/2$ and hence $\hat{\pi}_k = \tilde{\pi}_k $ holds for all $k$ when $n \geq N(\omega,\epsilon,d)$. This implies 
    \[
    \tilde{\varphi}_1(\bz) \equiv \hat{\varphi}_1(\bz),
    \]
    \[
    (\MP_n - \MP) [\hat{\varphi}_1(\bZ)-\tilde{\varphi}_1(\bZ)]=0
    \]
    when $n \geq N(\omega,\epsilon,d)$. This clearly implies 
    \[
    \sqrt{n}(\MP_n - \MP) [\hat{\varphi}_1(\bZ)-\tilde{\varphi}_1(\bZ)] \rightarrow 0
    \]
    almost surely since the left-hand side is exactly $0$ when $n$ is large enough. 
    We conclude
    \begin{equation}\label{eq:diff-varphi}
    (\MP_n - \MP) [\hat{\varphi}_1(\bZ)-\tilde{\varphi}_1(\bZ)] = o_{\MP}(1/\sqrt{n}) 
    \end{equation}
    Combining \eqref{eq:varphitilde-neg} and \eqref{eq:diff-varphi}, we have
    \[
    (\MP_n - \MP) [\hat{\varphi}_1(\bZ)-\varphi_1(\bZ)] = o_{\MP}(1/\sqrt{n}).
    \]
\end{proof}

\section{Proof of Auxiliary Results}\label{sec:proof-auxiliary}

\subsection{Proof of Lemma \ref{lemma:bound-var}}
\begin{proof}
First consider the boundary. If $y=0$ then the function $f_2$ is reduced to
\[
g_1(x) = (n-1)(1-x)^{n-2} - n(1-x)^{n-1}, 0 \leq x \leq 1,
\]
\[
g_1'(x) = (n-1)(1-x)^{n-3}(2-nx).
\]
So $|g_1(x)|\leq \max (|g_1(0)|, |g_1(1)|, |g_1(2/n)|)$ and $g_1(0) = -1, g_1(1) = 0, g_1(2/n) = \left(\frac{n-2}{n}\right)^{n-2} \leq 1$. This shows $|g_1(x)| \leq 1$. On the boundary $x=0$ similar arguments hold. Now consider the boundary $x+y=1$, the function $f_2$ is reduced to
\[
|f_2(x,y)| = nx^{n-1}(1-x)^{n-1} \leq  \frac{n}{4^{n-1}} \leq 1
\]
Now we consider the interior of the triangle. By taking the partial derivative (or by noting that $x,y$ have symmetric roles in the function $f_2$) we see the maximizer of $|f_2(x,y)|$ must lie on the line $x=y$. So we define
\[
g_2(x) = (n-1)(1-2x)^{n-2} - n(1-x)^{2n-2}, 0 \leq x \leq 1/2
\]
and only need to show $|g_2(x)| \leq 1$. Let $X_0$ be the set of stationary points of $g_2$ on $[0,1/2]$. Any stationary point $x_0 \in X_0$ must satisfy $g_2'(x_0) = 0$, i.e.
\[
n(1-x_0)^{2n-3} = (n-2)(1-2x_0)^{n-3}.
\]
So for any $x_0 \in X_0$,
\[
g_2(x_0) = (n-1)(1-2x_0)^{n-2} - (n-2)(1-2x_0)^{n-3}(1-x_0) = (1-2x_0)^{n-3}(1-nx_0).
\]
We then define a new function to show $|g_2(x_0)| \leq 1$ for any $x_0 \in X_0$. Let 
\[
g_3(x) = (1-2x)^{n-3}(1-nx), 0 \leq x \leq 1/2
\]
\[
g_3'(x) = (n-2)(2nx-3)(1-2x)^{n-4}.
\]
So we have $|g_3(x)| \leq \max(|g_3(0)|,|g_3(1/2)|, |g_3(3/2n)|)$ and $g_3(0) = 1, g_3(1/2) = 0, |g_3(3/2n)| = |\frac{1}{2}(\frac{n-3}{n})^{n-3}| \leq 1/2$. This shows $|g_3(x)| \leq 1$ on $[0,1/2]$, which implies $|g_2(x_0)| \leq 1$ for any $x \in X_0$. Note that $g_2(0)=-1, g_2(1/2) = -n/4^{n-1}$ We conclude
\[
|g_2(x)| \leq \max\left( \max_{x_0 \in X_0 }|g_2(x)|, 1, n/4^{n-1} \right) \leq 1
\]
\end{proof}

\subsection{Proof of Lemma \ref{lemma:poission-minimax}}
\begin{proof}
    For any $\gamma > 0$, let $\hat{\psi}_1(n)$ be a near-minimax optimal estimator of $\psi_1(\bp)$ for fixed sample $n$, i.e.,
    \[
    \sup_{\MP \in \mathcal{D}(\epsilon)} \ME_{\MP}\left[\left(\hat{\psi}_1(n) - \psi_1(\bp)\right)^2\right] \leq \gamma + R^*(d,n).
    \]
    Note we emphasize the dependency of $\psi_1$ on $\bp=(p_1,\dots)$. Under a Poisson-sampling model $\MP \in \mathcal{D}(\epsilon, \delta)$, let $n' = \sum_{k}n\hat{p}_k \sim \text{Poi}(n\sum_{k} p_k)$ and construct an estimator $\tilde{\psi}_1 = \hat{\psi}_1(n')$. Note that conditioned on $n'=m, (n\hat{p}_1, \dots,) \sim \text{Multinomial}\left(m, \frac{\bp}{\sum_{k}p_k}\right) $. By definition, the model with probability vector $\frac{\bp}{\sum_{k}p_k} $ and the same propensity score, regression functions as $\MP$ is in $\mathcal{D}(\epsilon) $. Under the Poisson-sampling model we have
    \begin{equation*}
        \begin{aligned}
            &\, \ME\left[ \left( \tilde{\psi}_1 -\psi_1\left( \frac{\bp}{\sum_{k}p_k}\right) \right)^2 \right] \\
            = &\, \sum_{m=0}^{\infty} \ME\left[ \left( \tilde{\psi}_1 -\psi_1\left( \frac{\bp}{\sum_{k}p_k}\right) \right)^2   \mid  n'=m\right] \MP(n'=m)\\
            \leq & \, \sum_{m=0}^{\infty} R^*(d,m)\MP(n'=m) + \gamma.
        \end{aligned}
    \end{equation*}
    Since $n \mapsto R^*(d, n)$ is non-increasing for fixed $d$ and $R^*(d, n) \leq 1$ , we have
    \[
    \begin{aligned}
        &\,\ME\left[ \left( \tilde{\psi}_1 -\psi_1\left( \frac{\bp}{\sum_{k}p_k}\right) \right)^2 \right] \\
        \leq &\, \sum_{m \geq n / 2} R^*(d, m) \mathbb{P}\left[n^{\prime}=m\right]+ \mathbb{P}\left[n^{\prime} \leq \frac{n}{2}\right] + \gamma \\
        \leq &\, R^*(d, n/2) + \exp(-n/50)+\gamma.
    \end{aligned}
    \]
    The last inequality follows from Chernoff bound and $|\sum_{k}p_k-1|\leq \delta < 1/3$. The difference between $\psi_1\left( \frac{\bp}{\sum_{k}p_k}\right)$ and $\psi_1\left( {\bp}\right)$ is
    \[
    \begin{aligned}
        &\, \bigg |\psi_1\left( \frac{\bp}{\sum_{k}p_k}\right) - \psi_1\left( {\bp}\right) \bigg |\\
        = &\, \big|\sum_{k}p_k-1\big| \frac{\sum_k p_k \mu_{1k}}{\sum_k p_k} \\
        \leq &\, \delta.
    \end{aligned}
    \]
    The last inequality follows since $\frac{\sum_k p_k \mu_{1k}}{\sum_k p_k} $ is a weighted average of quantities bounded by 1. Thus we have
    \[
    \begin{aligned}
        &\, \frac{1}{2} \ME\left[ \left( \tilde{\psi}_1 -\psi_1\left( {\bp}\right) \right)^2 \right]\\
        \leq &\, \ME\left[ \left( \tilde{\psi}_1 -\psi_1\left( \frac{\bp}{\sum_{k}p_k}\right) \right)^2 \right] +\left(\psi_1\left( \frac{\bp}{\sum_{k}p_k}\right) - \psi_1\left( {\bp}\right) \right)^2 \\
        \leq &\, R^*(d, n/2) + \exp(-n/50)+\gamma + \delta^2.
    \end{aligned}
    \]
    Take supremum over $\MP \in D(\epsilon, \delta)$ and since $\gamma$ is arbitrary, we have
    \[
    \frac{1}{2}\tilde{R}^*(d, n, \delta) \leq R^*(d, n/2) + \exp(-n/50) + \delta^2.
    \]
\end{proof}

\subsection{Proof of Lemma \ref{lemma:prior-dist}}
\begin{proof}
We claim it suffices to find probability measures $\omega_0, \omega_1$ on $[c/K^2,1]$ such that
\begin{enumerate}
    \item 
    \begin{equation}\label{eq:matched-moments}
    \ME_{\omega_0}[X^l] = \ME_{\omega_1}[X^l] \text{ for all } l=-1,0,1,\dots, 3K.
    \end{equation}
    \item 
    \[
    \left | \ME_{\omega_0}\left[\frac{X}{X+c/K^2} \right]- \ME_{\omega_1}\left[\frac{X}{X+c/K^2} \right] \right| \geq c'.
    \]
\end{enumerate}
Note that here we use $X$ to denote a different random variable from the covariate in the maintexts. We first show the claim leads to the results in Lemma \ref{lemma:prior-dist}. With $\omega_i$ ($i=0,1$) we define $\tilde{\omega}_i$ supported on $\{0\} \cup [c/K^2,1]$ such that
\[
\tilde{\omega}_i (dx) = \frac{c}{K^2x}\omega_i(dx) + \left(1-\ME\left[\frac{c}{K^2X} \right] \right) \delta_0(dx).
\]
And for $X \sim \tilde{\omega}_i$, let $\nu_i$ be the distribution of $\frac{c_1 \log n}{n}X$. We let $p \sim \nu_i$, 
\[
\pi = \left\{
    \begin{array}{rl}
    \epsilon \left(1+\frac{c_1 \log n}{np} \frac{c}{K^2} \right) & \text{if } p>0,\\
    \epsilon & \text{if } p=0
    \end{array} \right.
\]
and $\mu = \epsilon /\pi$. Note that $\pi, \mu$ are both functions of $p$. We then verify the properties in Lemma \ref{lemma:prior-dist} with joint distribution $\mu_i$ defined above. 
\begin{enumerate}
    \item $\text{supp}(\nu_i) \subseteq \{0\} \cup \frac{c_1 \log n}{n}[c/K^2,1]$ implies the range of $p$.
    \item For $ i,j,k \geq 0$ and $1\leq i+ j +k \leq 3K$ we have
    \[
    \begin{aligned}
        &\, \ME_{\mu_0}[p^i (p \pi)^j (p \pi \mu)^k] \\
        = &\, \ME_{\mu_0}[p^i (p \pi)^j (p \pi \mu)^k I(p>0)] \\
        = &\, \ME_{\nu_0}\left[ p^{i+k} \epsilon^{j+k} \left( 
p+\frac{cc_1 \log n}{nK^2} \right)^j I(p>0) \right] \\
        = &\, \ME_{\tilde{\omega}_0} \left[ \left(\frac{c_1 \log n}{n} \right)^{i+j+k} \epsilon^{j+k} X^{i+k} \left(X + \frac{c}{K^2} \right)^j I(X>0) \right] \\
        = &\, \int \left(\frac{c_1 \log n}{n} \right)^{i+j+k} \epsilon^{j+k} \frac{c}{K^2} x^{i+k-1} \left(x+ \frac{c}{K^2} \right)^j \omega_0 (dx) \\
        = &\, \int \left(\frac{c_1 \log n}{n} \right)^{i+j+k} \epsilon^{j+k} \frac{c}{K^2} x^{i+k-1} \left(x+ \frac{c}{K^2} \right)^j \omega_1 (dx) \\
        = &\, \ME_{\mu_1}[p^i (p \pi)^j (p \pi \mu)^k].
    \end{aligned}
    \]
    The first four equations follow from the definition of distributions above and the fifth follows from \eqref{eq:matched-moments}. 
    \item 
    \[
    \begin{aligned}
        & \, \ME_{\mu_i}[p]\\
        = &\, \int_{p>0} p \nu_i(dp) \\
        = &\, \int_{x>0}\frac{c_1 \log n}{n} x \tilde{\omega}_i (dx)\\
        = &\, \int_{x>0}\frac{c_1 \log n}{n} \frac{c}{K^2} \tilde{\omega}_i (dx)\\
        = &\, \frac{c c_1}{c_2^2} \frac{1}{n \log n} \leq \frac{c c_1 c_3}{c_2^2 d} \leq \frac{1}{d}
    \end{aligned}
    \]
    as long as $c c_1 c_3/c_2^2 \leq 1$.
    \item 
    \[
    \begin{aligned}
        &\, |\ME_{\mu_0}[p \mu] - \ME_{\mu_1}[p \mu]|\\
        = &\, \left | \ME_{\nu_0}\left[ \frac{p^2}{p+ \frac{c_1 \log n}{n}\frac{c}{K^2} } \right] -  \ME_{\nu_1}\left[ \frac{p^2}{p+ \frac{c_1 \log n}{n}\frac{c}{K^2} } \right]\right| \\
        =&\, \frac{c_1 \log n}{n} \left | \ME_{\tilde{\omega}_0}\left[ \frac{X^2}{X+ \frac{c}{K^2} } \right] -  \ME_{\tilde{\omega}_1}\left[ \frac{X^2}{X+ \frac{c}{K^2} } \right]\right|\\
        = &\, \frac{c c_1 \log n}{n K^2} \left | \ME_{{\omega}_0}\left[ \frac{X}{X+ \frac{c}{K^2} } \right] -  \ME_{{\omega}_1}\left[ \frac{X}{X+ \frac{c}{K^2} } \right]\right| \\
        \geq &\, \frac{c c' c_1}{c_2^2 n \log n}.
    \end{aligned}
    \]
    And we define $c_4 = cc' c_1/c_2^2$.
\end{enumerate}
    Thus we only need to prove the claim. By duality of polynomial approximation and moment matching (see e.g., Lemma 19 in \cite{han2020optimal}) it suffices to prove 
    \[
    \inf_{P \in \text{span}\{1/x, 1, x, \dots, x^{3K}\} } \max_{x \in [c/K^2,1]} \left | \frac{x}{x+c/K^2}-P(x) \right| \geq c'/2.
    \]
    Let 
    \begin{equation}\label{eq:poly-approximation}
    E_n(f; S) = \min_{\text{deg}(P) \leq n} \max_{x \in S} |f(x)-P(x)|
    \end{equation}
    be the best polynomial approximation error (with polynomial's degree smaller than $n$) of $f$ on a set $S$. We have the following lemma.
    \begin{lemma}\label{lemma:poly-approximation}
        There exists $c_0, c' >0$ such that 
        \[
        \liminf_{n \rightarrow \infty} \inf_{\alpha \in \mathbb{R}} E_n \left(\frac{x}{x+c_0n^{-2}} + \frac{\alpha}{x}; [c_0 n^{-2},1] \right) \geq c' >0.
        \]
    \end{lemma}
    Let $c = c_0/9$, we then have for any $P \in \text{span}\{1/x, 1, x,\dots, x^{3K} \}$, write $P(x) = \frac{\alpha_{-1}}{x}+ \sum_{k=0}^{3K}\alpha_k x^k$ and let $P_{\geq}(x) = \sum_{k=0}^{3K}\alpha_k x^k$
    \[
    \begin{aligned}
        & \,\max_{x \in [c/K^2,1]} \left| \frac{x}{x+c/K^2} -P(x) \right| \\
        \geq &\, \inf_{\alpha \in \mathbb{R}} \max_{x \in [c/K^2,1]} \left| \frac{x}{x+c/K^2} + \frac{\alpha}{x} -P_{\geq}(x) \right| \\
        = &\, \inf_{\alpha \in \mathbb{R}} \max_{x \in [c_0/(3K)^2,1]} \left| \frac{x}{x+c_0/(3K)^2} + \frac{\alpha}{x} -P_{\geq}(x) \right| \\
        \geq &\, \inf_{\alpha \in \mathbb{R}} E_{3K} \left( \frac{x}{x+c_0/(3K)^2} + \frac{\alpha}{x}; [c_0 / (3K)^2, 1] \right) \\
        \geq &\, \frac{c'}{2},
    \end{aligned}
    \]
    as $K = c_2 \log n$ is large enough.

    We then prove Lemma \ref{lemma:poly-approximation}. By Theorem 7.2.4 in \cite{ditzian2012moduli} (and use the notation there), we have
    \begin{equation}\label{eq:moduli-smooth}
    \omega_{\varphi}^1 \left(f,\frac{1}{n} \right)_{\infty} \leq \frac{M}{n}\sum_{\ell=0}^n E_{\ell}(f;[0,1]),        
    \end{equation}
    where
    \[
    \omega_{\varphi}^1 \left(f,\frac{1}{n} \right)_{\infty} = \sup_{0 < h \leq 1/n} \sup_{x,x+h\varphi(x) \in [0,1]} |f(x+h\varphi(x))-f(x)|, \, \varphi(x) = \sqrt{x(1-x)}
    \]
    and $M$ is an absolute constant. Let $f_{\alpha}(x) = \frac{x}{x+c_0 n^{-2}} + \frac{\alpha}{x}$ and $\tilde{f}_{\alpha}(x) = f_{\alpha}\left(c_0n^{-2} + (1-c_0n^{-2}) x\right), x\in[0,1]$. Denote $m = \lceil n/\sqrt{c_0} \rceil $. Here $c_0$ is a constant to be chosen. Consider two cases of $\alpha$ separately:
    
    \textbf{Case 1:} $|\alpha| \leq c_0 / n^2$. By \eqref{eq:moduli-smooth} we have
    \[
    \begin{aligned}
    &\,E_n(f_{\alpha}; [c_0n^{-2},1]) \\
    = &\,E_n(\tilde{f}_{\alpha}; [0,1]) \\
    \geq &\, \frac{1}{m-n+1} \sum_{\ell =n}^m \ME_{\ell} (\tilde{f}_{\alpha}; [0,1]) \\
    \geq &\, \frac{1}{m} \sum_{\ell =n}^m \ME_{\ell} \left(\tilde{f}_{\alpha}; [0,1] \right) \\
    \geq &\, \frac{1}{M} \omega_{\varphi}^1 \left(\tilde{f}_{\alpha},\frac{1}{m} \right)_{\infty} - \frac{1}{m} \sum_{\ell =0}^{n-1} \ME_{\ell} \left(\tilde{f}_{\alpha}; [0,1] \right) \\
    \geq &\, \frac{1}{M} \omega_{\varphi}^1 \left(\tilde{f}_{\alpha},\frac{1}{m} \right)_{\infty} - \frac{n}{m}  \ME_{0} \left(\tilde{f}_{\alpha}; [0,1] \right) \\
    \geq &\, \frac{1}{M} \omega_{\varphi}^1 \left(\tilde{f}_{\alpha},\frac{1}{m} \right)_{\infty} - 2 \sqrt{c_0},
    \end{aligned}
    \]
    where the first and the fourth inequality follow from the monotonity of $E_n$, the third inequality applies \eqref{eq:moduli-smooth} and the last one follows from $|\tilde{f}_{\alpha}| \leq 2$. For the first term, we have
    \[
    \begin{aligned}
        &\, \omega_{\varphi}^1 \left(\tilde{f}_{\alpha},\frac{1}{m} \right)_{\infty} \\
        \geq &\, \sup_{t_1, t_2 \in [3,4]} \left | \tilde{f}_{\alpha}\left(\frac{t_1}{m^2} \right) -\tilde{f}_{\alpha}\left(\frac{t_2}{m^2} \right)  \right| \\
        = &\, \sup_{t_1, t_2 \in [3,4]} \left | {f}_{\alpha}\left(\frac{c_0}{n^2} + \left(1-\frac{c_0}{n^2} \right) \frac{t_1}{m^2} \right) -{f}_{\alpha}\left(\frac{c_0}{n^2} + \left(1-\frac{c_0}{n^2} \right) \frac{t_2}{m^2} \right)  \right|, \\
    \end{aligned}
    \]
    where the inequality follows from the fact that when $m$ is sufficiently large,
    \[
    \frac{t_2-t_1}{m^2} \leq \frac{1}{m} \sqrt{\frac{t_2}{m} \left(1-\frac{t_2}{m} \right)}
    \]
    holds for all $t_1 < t_2$ and $t_1,t_2 \in [3,4]$. Note that as $m \rightarrow \infty$
    \[
    1+ \left(\frac{n^2}{c_0}-1 \right)\frac{3}{m^2} \rightarrow 4,
    \]
    \[
    1+ \left(\frac{n^2}{c_0}-1 \right)\frac{4}{m^2} \rightarrow 5.
    \]
    Thus we have
    \[
    \begin{aligned}
        & \, \sup_{t_1, t_2 \in [3,4]} \left | {f}_{\alpha}\left(\frac{c_0}{n^2} + \left(1-\frac{c_0}{n^2} \right) \frac{t_1}{m^2} \right) -{f}_{\alpha}\left(\frac{c_0}{n^2} + \left(1-\frac{c_0}{n^2} \right) \frac{t_2}{m^2} \right)  \right|\\
        \geq &\, \sup_{t_1, t_2 \in [3.5, 5]}  \left | {f}_{\alpha}\left(\frac{c_0}{n^2} t_1 \right) - {f}_{\alpha}\left(\frac{c_0}{n^2} t_2 \right)\right| \\
        \geq &\, \sup_{t \in[3.5,4]} \left | {f}_{\alpha}\left(\frac{c_0}{n^2} \left(t+ \frac{1}{2} \right) \right) - {f}_{\alpha}\left(\frac{c_0}{n^2} t \right)\right| \\
        = &\, \sup_{t \in[3.5,4]} \left | \frac{t+1/2}{t+3/2} - \frac{t}{t+1} + \frac{\alpha n^2}{c_0} \left(\frac{1}{t+1/2} -\frac{1}{t} \right) \right |\\
        \geq &\, \inf_{\beta} \sup_{t \in[3.5,4]} \left | \frac{t+1/2}{t+3/2} - \frac{t}{t+1} + \beta \left(\frac{1}{t+1/2} -\frac{1}{t} \right) \right |\\
        \geq &\, a_1,
    \end{aligned}
    \]
    where $a_1$ is a positive constant independent of $\alpha, n, c_0$ since $\frac{t+1/2}{t+3/2} - \frac{t}{t+1}$ and $\frac{1}{t+1/2} -\frac{1}{t}$ are linearly independent. Hence we have
    \[
    E_n(f_{\alpha}; [c_0n^{-2},1]) \geq \frac{a_1}{M}-2\sqrt{c_0}, \, |\alpha|\leq \frac{c_0}{n^2}.
    \]
    \textbf{Case 2:} $|\alpha| >{c_0}/{n^2}$, similar to the proof in Case 1 we have
    \[
    \begin{aligned}
        &\, E_n(f_{\alpha}; [c_0n^{-2},1]) \\
        \geq &\, \frac{1}{M} \omega_{\varphi}^1 \left(\tilde{f}_{\alpha},\frac{1}{m} \right)_{\infty} - \sqrt{c_o}  \ME_{0} \left(\tilde{f}_{\alpha}; [0,1] \right) \\
        \geq &\, \frac{1}{M} \omega_{\varphi}^1 \left(\tilde{f}_{\alpha},\frac{1}{m} \right)_{\infty} - \sqrt{c_o} \left(1+ \frac{|\alpha|n^2}{c_0} \right). 
    \end{aligned}
    \]
    The second inequality follows from $|f_{\alpha}| \leq \left(1 + \frac{|\alpha|n^2}{c_0} \right)$. And for the first term by the same argument in Case 1,
    \[
    \begin{aligned}
        &\, \omega_{\varphi}^1 \left(\tilde{f}_{\alpha},\frac{1}{m} \right)_{\infty}\\
        \geq &\, \sup_{t \in[3.5,4]} \left | {f}_{\alpha}\left(\frac{c_0}{n^2} \left(t+ \frac{1}{2} \right) \right) - {f}_{\alpha}\left(\frac{c_0}{n^2} t \right)\right| \\
        = &\, \sup_{t \in[3.5,4]} \left | \frac{t+1/2}{t+3/2} - \frac{t}{t+1} + \frac{\alpha n^2}{c_0} \left(\frac{1}{t+1/2} -\frac{1}{t} \right) \right |\\
        \geq &\, \frac{\alpha n^2}{c_0} \inf_{\beta \in \mathbb{R}} \sup_{t \in[3.5,4]} \left | \beta \left(\frac{t+1/2}{t+3/2} - \frac{t}{t+1}\right) + \frac{1}{t+1/2} -\frac{1}{t}  \right |\\
        \geq & \, \frac{\alpha n^2}{c_0} a_2,
    \end{aligned}
    \]
    where $a_2 $ is a positive constant independent of $\alpha, n, c_0$ again since $\frac{t+1/2}{t+3/2} - \frac{t}{t+1}$ and $\frac{1}{t+1/2} -\frac{1}{t}$ are linearly independent. Hence we have for $|\alpha|> c_0/n^2$ (choose $c_0$ such that $\frac{a_2}{M}-\sqrt{c_0}>0$)
    \[
    \begin{aligned}
       &\, E_n(f_{\alpha}; [c_0n^{-2},1])  \\
       \geq &\, \frac{\alpha n^2}{c_0 M}a_2 -\sqrt{c_o} \left(1+ \frac{|\alpha|n^2}{c_0} \right) \\
       \geq &\, \inf_{|\alpha|> c_0/n^2} \frac{|\alpha|n^2}{c_0}  \left(\frac{a_2}{M}-\sqrt{c_0} \right) -\sqrt{c_0} \\
       \geq & \, \frac{a_2}{M}-2\sqrt{c_0}.
    \end{aligned}
    \]
    Combining Case 1 and Case 2, we conclude for any $\alpha \in \mathbb{R}$
    \[
    E_n(f_{\alpha}; [c_0n^{-2},1]) \geq \min \left\{  \frac{a_1}{M} - 2\sqrt{c_0}, \frac{a_2}{M} - 2\sqrt{c_0}\right\}.
    \]
    Choosing $c_0>0$ sufficiently small completes the proof.
\end{proof}

\subsection{Proof of Lemma \ref{lemma:poissonization}}\label{sec:proof-poissonization}
\begin{proof}
It is known that the conditional distribution of $Y^n$ given the summation of its components is multinomial:
\begin{equation}\label{eq:poisson-multinomial}
    \bY^n \mid \sum_{k=1}^d Y_k^n = m \stackrel{d}{=} \bX^m
\end{equation}
Hence we have
\begin{equation*}
    \begin{aligned}
        & \, \ME[f(Y_1^n,\dots, Y_d^n)]\\
        \geq&\,  \sum_{m=0}^n \ME\left[f(Y_1^n,\dots, Y_d^n) | \sum_{k=1}^d Y_k^n = m\right] \MP\left(\sum_{k=1}^d Y_k^n = m\right) \\
        = &\, \sum_{m=0}^n \ME\left[f(X_1^m,\dots, X_d^m)\right] \MP\left(\sum_{k=1}^d Y_k^n = m\right) \\
        \geq &\,  \ME\left[f(X_1^n,\dots, X_d^n)\right] \MP \left(\sum_{k=1}^d Y_k^n \leq n\right) \\
        \geq  &\, \frac{1}{2} \ME\left[f(X_1^n,\dots, X_d^n)\right] 
    \end{aligned}
\end{equation*}
where the first inequality follows from truncating the summation, the first equation follows from \eqref{eq:poisson-multinomial}, the second inequality follows from monotonicity of $\ME[f(X_1^n,\dots, X_d^n)]$ and the last inequality follows from $\sum_{k=1}^d Y_k^n \sim \text{Poisson}(n)$ and $\MP(\sum_{k=1}^d Y_k^n \leq n) \geq 1/2$.
\end{proof}

\end{document}